\input ./style/arxiv-general.cfg
\documentclass[aap,MSNbibl,seceqn,dvips]{arximspdf}
\makeatletter
   \@ifpackageloaded{graphicx}{}{\usepackage{graphicx}}
\makeatother
\usepackage{mathrsfs}
\doi{10.1214/14-AAP1065} % kopijuoti is laisko
\volume{25}
\issue{5}
\pubyear{2015}
\firstpage{2959}
\lastpage{3006}
\docsubty{FLA}

\makeatletter
\newcommand{\rrvert}{\vert}
\newcommand{\rrVert}{\Vert}
\newcommand{\llvert}{\vert}
\newcommand{\llVert}{\Vert}
\def\nicefrac{\frac}
\newcommand{\eqref}[1]{(\ref{#1})}
\newtheorem{stat}{Statement}[section]
\newtheorem{proposition}[stat]{Proposition}
\newtheorem{corollary}[stat]{Corollary}
\newtheorem{theorem}[stat]{Theorem}
\newtheorem{lemma}[stat]{Lemma}
\newtheorem{lemmaa}{Lemma}[section]
\newproclaim{definition}[stat]{Definition}
\newproclaim{remark}[stat]{Remark}
\newproclaim{OP}{Open Problems}
\newproclaim{example}[stat]{Example}
\newproclaim{nota}[stat]{Notation}
\newcommand{\Z}{\mathbf{Z}}
\newcommand{\R}{\mathbf{R}}

\newcommand{\lip}{\operatorname{Lip}_\sigma}

\renewcommand{\P}{\mathrm{P}}
\newcommand{\E}{\mathrm{E}}

\newcommand{\1}{\mathbf{1}}
\renewcommand{\d}{\mathrm{d}}
\renewcommand{\Re}{\operatorname{Re}}
\makeatother

\begin{document}
\begin{frontmatter}

\title{Semi-discrete semi-linear parabolic SPDEs}
\runtitle{Semi-discrete semi-linear parabolic SPDEs}
\begin{aug}
\author[A]{\fnms{Nicos} \snm{Georgiou}\thanksref{T1}\ead[label=e1]{N.Georgiou@sussex.ac.uk}},
\author[B]{\fnms{Mathew} \snm{Joseph}\thanksref{T1,T3}\ead[label=e2]{m.joseph@shef.ac.uk}\ead[label=u2,url]{http://mathew-joseph.staff.shef.ac.uk/}},
\author[C]{\fnms{Davar} \snm{Khoshnevisan}\thanksref{T2,T3}\ead[label=e3]{davar@math.utah.edu}\ead[label=u3,url]{http://www.math.utah.edu/\textasciitilde davar/}}
\and
\author[D]{\fnms{Shang-Yuan} \snm{Shiu}\corref{}\thanksref{T4}\ead[label=e4]{shiu@math.ncu.edu.tw}}

\thankstext{T1}{Supported in part by the NSF Grant DMS-07-47758.}
\thankstext{T2}{Supported in part by the NSF Grant DMS-13-07470.}
\thankstext{T3}{Supported in part by the NSF Grant DMS-10-06903.}
\thankstext{T4}{Supported in part by the NSC Grant
101-2115-M-008-10-MY2 and the NCU Grant 102G607-3.}
\runauthor{Georgiou, Joseph, Khoshnevisan and Shiu}

\affiliation{University of Sussex and University of Utah,
University of Sheffield,
University of Utah
and National Central University}

\address[A]{N. Georgiou\\
Department of Mathematics\\
University of Sussex\\
Pevensey 3 5C15\\
Brighton, BN1 9RH\\
United Kingdom\\
\printead{e1}}

\address[B]{M. Joseph\\
Department of Probability and Statistics\\
University of Sheffield\\
Sheffield, England, S3 7RH\\
United Kingdom\\
\printead{e2}\\
\printead{u2}}

\address[C]{D. Khoshnevisan\\
Department of Mathematics\\
University of Utah\\
Salt Lake City, Utah 84112-0090\\
USA\\
\printead{e3}\\
\printead{u3}}

\address[D]{S.-Y. Shiu\\
Department of Mathematics\\
National Central University\\
Jhongli City, Taoyuan County 32001\hspace*{22pt}\\
Taiwan\\
\printead{e4}}
\end{aug}
%
%\begin{aug}
%%\author[A]{\fnms{}~\snm{}\corref{}\ead[label=e1]{}}%,
%%\author[]{\fnms{}~\snm{}\ead[label=]{}}
%% \and
%%\author[]{\fnms{}~\snm{}\ead[label=]{}}
%%\runauthor{}
%%\affiliation{}
%%\dedicated{}
%\address[A]{\\\printead{e1}}
%%\address[]{\\\printead{}}
%\end{aug}

% HISTORY:
\received{\smonth{11} \syear{2013}}
\revised{\smonth{9} \syear{2014}}
%\accepted{\smonth{} \syear{}}

% ABSTRACT
%
\begin{abstract}
Consider an infinite system
\[
\partial_t u_t(x) = (\mathscr{L} u_t) (x)
+ \sigma\bigl(u_t(x)\bigr) \partial_t
B_t(x)
\]
of interacting It\^o diffusions,
started at a nonnegative deterministic bounded initial profile.
%where $\{B(x)\}_{x\in\Z^d}$ is a field of i.i.d.\ Brownian motions and
%$\mathscr{L}$ denotes the generator of a continuous-time random
%walk on the d-dimensional lattice.
We study local and global features of the solution under
standard regularity assumptions on the nonlinearity $\sigma$.
We will show that,
locally in time, the solution behaves as a collection of independent diffusions.
We prove also that the $k$th moment Lyapunov exponent is frequently
of sharp order $k^2$, in contrast to the continuous-space stochastic
heat equation
whose $k$th moment Lyapunov exponent can be of sharp order $k^3$.
When the underlying walk is transient and the noise level is
sufficiently low, we prove
also that the solution is a.s. uniformly dissipative
provided that the initial profile is in $\ell^1(\Z^d)$.
\end{abstract}

% KEYWORDS
\begin{keyword}[class=AMS]
\kwd[Primary ]{60J60}
\kwd{60K35}
\kwd{60K37}
\kwd[; secondary ]{47B80}
\kwd{60H25}
\end{keyword}

\begin{keyword}
\kwd{The stochastic heat equation}
\kwd{semi-discrete stochastic heat equation}
\kwd{discrete space}
\kwd{parabolic Anderson model}
\kwd{Lyapunov exponents}
\kwd{dissipative behavior}
\kwd{comparison principle}
\kwd{interacting diffusions}
\kwd{BDG inequality}
\end{keyword}
%
%\begin{keyword}[class=AMS]
%\kwd[Primary ]{}
%\kwd{}
%\kwd[; secondary ]{}
%\end{keyword}
%\begin{keyword}
%\kwd{}
%\end{keyword}
\end{frontmatter}

%s1 #&#
\section{Introduction}\label{sec1}

\subsection*{Model and motivation}
We propose to study the following system of infinitely-many interacting
diffusions:
{\renewcommand{\theequation}{SHE}
%e1.1 #&#
\begin{equation}\label{SHE}
\frac{\d u_t(x)}{\d t} = (\mathscr{L}u_t) (x) + \sigma
\bigl(u_t(x)\bigr) \frac{\d B_t(x)}{\d t}, %\eqno{}
\end{equation}}
\hspace*{-3pt}where $t>0$ denotes the \emph{time} variable and $x\in\Z^d$
is the \emph{space} variable. In parts of the literature, \eqref{SHE} is
thought of as a \emph{stochastic heat equation} on $(0,\infty)\times
\Z^d$,
viewed as a semi-discrete stochastic partial differential equation.

We interpret \eqref{SHE} as an infinite-dimensional
system of It\^o stochastic differential equations, where
$\{B(x)\}_{x\in\Z^d}$ denotes a field of independent
standard linear Brownian motions and
$\sigma\dvtx\R\to\R$
is a Lipschitz-continuous nonrandom function with
\setcounter{equation}{0}
%e1.1 #&#
\begin{equation}
\label{eq:sigma(0)=0} \sigma(0)=0.
\end{equation}
The drift operator $\mathscr{L}$ acts on the variable $x$ only,
and is the generator of a continuous-time random walk
$X:=\{X_t\}_{t\ge0}:=\{\sum_{j=1}^{N_t}Z_j\}_{t\ge0}$ on $\Z^d$
where $N_t$ is a Poisson process with jump-rate one
and the $Z_j$'s are i.i.d. random variables with values in $\Z^d$.
We consider only initial values $u_0\dvtx\Z^d\to\R$ that are nonrandom
and satisfy
%
%e1.2 #&#
\begin{equation}
\label{init} u_0(x)\ge0\qquad\mbox{for all }x\in\Z^d \quad\mbox{and}\quad
 0<\sup_{x\in\Z^d}u_0(x)<\infty,
\end{equation}
though some of the theory developed here applies to more general
initial profiles.
According to Shiga and Shimizu \cite{ShigaShimizu}, condition
\eqref{init} ensures that the system \eqref{SHE} has an a.s.-unique solution.

Such systems have been studied extensively
\cite{CKM,CM94,CFG,CoxGreven,CranstonMolchanov,CranstonMountfordShiga,CGM,Funaki,GrevenDenHollander,Mueller,Shiga,ShigaShimizu},
most commonly in the context of well-established
models of statistical mechanics, population genetics and related
models of infinitely-many interacting diffusion processes.
One of the central examples of this literature is the
\emph{parabolic Anderson model} \cite{CM94}---also known as
\emph{diffusion in random potential}---which is \eqref{SHE}
when the function $\sigma$ is linear. It is not hard to prove that,
for the parabolic Anderson model, the \emph{$k$th moment Lyapunov
exponent} $\gamma_k(u)$ exists
and is positive and finite for all real numbers $k\ge2$, where
%
%e1.3 #&#
\begin{equation}
\gamma_k(u):= \lim_{t\to\infty} t^{-1}\log\E
\bigl( \bigl|u_t(x)\bigr|^k \bigr).
\end{equation}
[One can prove that, because of \eqref{init}, $\gamma_k(u)$ does not
depend on $x$.]
Jensen's inequality readily implies
that $k\mapsto\gamma_k(u)$ is nondecreasing on
$[2 ,\infty)$. In the case that $\mathscr{L}$ denotes the discrete
Laplacian on $\Z$,
for example, the theory of
Carmona and Molchanov \cite{CM94} implies that
%
%e1.4 #&#
\begin{equation}
\label{CM-I} k\mapsto\gamma_k(u)\qquad\mbox{is strictly increasing on $[2
,\infty)$}.
\end{equation}
This property is referred to as \emph{intermittency} and suggests
that the random function $u$ develops very tall peaks that are distributed
over small space--time ``islands.'' Section~2.4 of Bertini and Cancrini
\cite{BC} and
Section~7.1 of Khoshnevisan \cite{K-CBMS} describe two heuristic derivations
of this ``peaking behavior'' from intermittency condition~\eqref{CM-I}.

In the present nonlinear setup, the Lyapunov exponents do not generally
exist. Therefore, one considers instead
the (maximal) \emph{bottom} and \emph{top Lyapunov exponents} of
the solution $u$ to \eqref{SHE}; those are, respectively, defined as
%
%e1.5 #&#
\begin{eqnarray}
\underline{\gamma}_k(u) &:= &\liminf_{t\to\infty}
\sup_{x\in\Z^d} t^{-1} \log\E \bigl( \bigl|u_t(x)\bigr|^k
\bigr),
\nonumber
\\[-8pt]
\\[-8pt]
\nonumber
\overline{\gamma}_k(u) &:= &\limsup_{t\to\infty} \sup
_{x\in\Z^d} t^{-1} \log\E \bigl( \bigl|u_t(x)\bigr|^k
\bigr). %
\end{eqnarray}
Whenever\vspace*{1pt} $\underline{\gamma}_k(u)=\overline{\gamma}_k(u)$,
we write $\gamma_k(u)$ for their common value and think of
$\gamma_k(u)$ as the \emph{$k$th moment Lyapunov exponent} of $u$.

In the present nonlinear setting, intermittency is defined as the
property that the functions $k \mapsto k^{-1}\underline{\gamma}_k(u)$
and $k\mapsto k^{-1}\overline{\gamma}_k(u)$
are both strictly increasing on $[2 ,\infty)$. Because of convexity,
one can establish intermittency when $\underline{\gamma}_2(u)>0$ and
$\overline{\gamma}_k(u)<\infty$ for all $k\ge2$, \cite{K-CBMS}, Proposition~7.2.

We will prove that $\overline{\gamma}_k(u)$ is
generically finite for all $k\in[2 ,\infty)$;
see Theorem~\ref{th:exist:unique}.
Therefore, as far as matters of intermittency are concerned,
it remains to establish the positivity of the bottom Lyapunov exponent.
This is a nongeneric property. To wit,
when $u_0$ is a constant and $\mathscr{L}$ is the discrete
Laplacian on $\Z^d$, Carmona and Molchanov \cite{CM94} have shown that
$\gamma_2(u)>0$ if and only if $d\in\{1 ,2\}$.

This paper is concerned with various results that surround this general topic.
As a first example of the type of result that we will establish, let
us point out the following, which is related to the mentioned peaking
property of the parabolic Anderson model:
One can apply Theorem~\ref{th:exist:unique} below to the parabolic
Anderson model
in order to see that for all $t>0$ there exist finite and positive
constants $a(t)$
and $A(t)$ such that
%
%e1.6 #&#
\begin{equation}
\label{CM-U} a(t) k^2 \le\log\E \bigl( \bigl|u_t(x)\bigr|^k
\bigr) \le A(t) k^2,
\end{equation}
uniformly for all $x\in\Z^d$ and $k\in[2 ,\infty)$.
It is then possible to combine this bound, together with the method of
Conus et al.
\cite{CJK13}, in order to estimate the size of the peaks of the solution
relative to the spatial variable $x$. For example, if $\mathscr{L}$ is
finite range,
then it is possible to prove that the tall peaks grow at all times
as $\exp\{\operatorname{const}\,\cdot\,
\sqrt{\log\|x\|}\}$ for large values of $\|x\|$ and more precisely, that
%
%e1.7 #&#
\begin{equation}
\label{peaks} 0<\limsup_{\|x\|\to\infty} \frac{\log|u_t(x)|}{\sqrt{\log\|x\|}} <\infty
\qquad\mbox{a.s.}
\end{equation}
We will not establish this fact since it follows fairly readily from
\eqref{CM-U}
and the methods of \cite{CJK13}. Instead let us return to \eqref{SHE} in its
nonlinear form.

In general, properties such as \eqref{CM-U},
whence \eqref{peaks}, can be shown to fail. This is so, for example,
when $\sigma$ is bounded;
see Conus et al. \cite{CJK13} for analogous results.
Therefore, in order to prove \eqref{CM-U}
for \eqref{SHE} we need to impose some growth conditions on the nonlinearity
$\sigma$.

Define
%
%e1.8 #&#
\begin{equation}
\lip:= \sup_{-\infty<a\neq b<\infty} \biggl\llvert \frac{\sigma(a)-\sigma(b)}{a-b}\biggr
\rrvert , \qquad \ell_\sigma:=\inf_{z\in\R} \biggl\llvert
\frac{\sigma(z)}{z} \biggr\rrvert . \label{def:Lipell}
\end{equation}
Because of \eqref{eq:sigma(0)=0},
%
%e1.9 #&#
\begin{equation}
\ell_\sigma\le\biggl\llvert \frac{\sigma(z)}{z}\biggr\rrvert \le\lip\qquad
\mbox{for all $z\in\R\setminus\{0\}$}.
\end{equation}

Suppose that $\ell_{\sigma}>0$ and $\sigma(z)>0$ for all $z>0$.
Then we will prove that \eqref{CM-U} holds, and in particular the
Lyapunov exponents are
always positive and finite; see Theorem~\ref{th:exist:unique}.
Property \eqref{CM-U}
contrasts sharply with the continuous-space analogues of
\eqref{SHE} wherein typically the $k$th moment Lyapunov exponents are
of sharp order $k^3$ as $k\to\infty$ \cite{ACQ,BorCor,BQS,BC,FK} and
implies that,
although the intermittency peaks of semi-discrete stochastic PDEs grow
rapidly with time,
they grow far less rapidly than those of fully-continuous stochastic PDEs.

Not much else seems to be known about the detailed
behavior of the Lyapunov exponents of the solution
to \eqref{SHE}: Ours seem to be the only methods that thus far have succeeded in
analyzing the asymptotics of the Lyapunov exponents of the solution to
\eqref{SHE} when $\sigma$ is nonlinear and/or $\mathscr{L}$ is nonlocal;
Borodin and Corwin have found a few remarkable instances of \eqref{SHE} where
all integer-moment Lyapunov exponents can be computed precisely.

The second nontrivial contribution of this paper concerns the local behavior
of the solution to \eqref{SHE}. The statement is that there frequently
exists a monotone function $S$ such that for all $t>0$ fixed,
$\eta_\tau$ converges to white noise on $\Z^d$ as $\tau\downarrow
0$, where
%
%e1.10 #&#
\begin{equation}
\eta_\tau(x) := \frac{ S(u_{t+\tau}(x)) - S(u_t(x)) }{\sqrt\tau} \qquad\mbox{for all $x\in
\Z^d$}.
\end{equation}
See Theorem~\ref{th:local:RadonNikodym}, and especially
Theorem~\ref{th:CLT}, for more details.
One can interpret our result as saying that the function $S$ is
the infinite-dimensional analogue of the scale function
for a finite-dimensional diffusion. The function $S$
is in fact also an abstract Hopf--Cole transformation for many nonlinear
systems of the form \eqref{SHE};
see \cite{ACQ,BC} for the role of the
latter transformation in the continuous-space parabolic Anderson model.

Finally, we state and prove perhaps the most interesting contribution
of this
paper, namely that the solution to \eqref{SHE} is strongly dissipative
when \eqref{SHE} is weakly disordered. See Theorem~\ref{th:as} for a precise
statement.
The latter theorem implies that \eqref{SHE} is a model with hysterisis;
see Remark~\ref{rem:Hysteresis} for more details. Theorem~\ref{th:as}
is also significant because it
rules out the possibility of the existence of an Anderson mobility edge for
the present model \eqref{SHE}, when \eqref{SHE} is weakly disordered.
We are aware only of one such nonexistence theorem, the original
stationary Anderson model (but on ``tree graphs''); see the recent paper
by Aizenman and Warzel \cite{AizenmanWarzel}. Among other things,
the proof of Theorem~\ref{th:as} relies on a comparison theorem for
renewal processes,
which we state and prove in the \hyperref[app]{Appendix} [see Lemma~\ref{lem:RE:comparison}].
It is likely that our comparison theorem has other uses in applied probability
as well.

%s2 #&#
\section{Results}
In this section we present the main results of this paper.
Let us begin by making the assertions in the \hyperref[sec1]{Introduction} more precise.

%th2.1 #&#
\begin{theorem}
\label{th:exist:unique}
The nonlinear stochastic heat equation \mbox{\eqref{SHE}} has a solution
$u$ that is continuous in the variable $t$, and is
unique among all predictable random fields that satisfy
$\sup_{t\in[0,T]}\sup_{x\in\Z^d}\E(|u_t(x)|^2)<\infty$ for all $T>0$.
Moreover,
%
%e2.1 #&#
\begin{equation}
\label{eq:Upper:LE} \overline{\gamma}_k(u) \le 8\lip^2
k^2 \qquad\mbox{for all integers $k\ge2$}.
\end{equation}
Furthermore, $u_t(x)\ge0$ for all $t\ge0$ and $x\in\Z^d$ a.s., provided
that $u_0(x)\ge0$ for all $x\in\Z^d$.
Finally, if $\ell_\sigma>0$ and
$\sigma(x)> 0$ for all $x>0$, then for all $\varepsilon\in(0 ,1)$,
%
%e2.2 #&#
\begin{equation}
\label{eq:Lower:LE} \underline{\gamma}_k(u) \ge(1-\varepsilon)
\ell_\sigma^2 k^2 \qquad\mbox{for every integer $k\ge
\varepsilon^{-1}+\bigl(\varepsilon\ell _\sigma^2
\bigr)^{-1}$}.
\end{equation}
\end{theorem}

Our next result shows that, at each point $x\in\Z^d$,
the solution behaves locally in time like a Brownian motion.
Standard moment methods---which we will have to reproduce here as well---show
that $t\mapsto u_t(x)$ is almost surely a H\"older-continuous
random function for every H\"older exponent $<\nicefrac12$. The
following proves that the H\"older exponent $\nicefrac12$ is sharp.

%th2.2 #&#
\begin{theorem}[(A Radon--Nikod\'ym property)]\label{th:local:RadonNikodym}
For every $t\ge0$ and $x\in\Z^d$,
%
%e2.3 #&#
\begin{equation}
\label{eq:local:RN} \lim_{\tau\downarrow0}\frac{u_{t+\tau}(x) - u_t(x)}{B_{t+\tau}(x)
- B_t(x)} = \sigma \bigl(
u_t(x) \bigr)\qquad \mbox{in probability}.
\end{equation}
In addition,
%
%e2.4 #&#
\begin{equation}\quad
\label{eq:local:LIL} \limsup_{\tau\downarrow0}\frac{u_{t+\tau}(x) - u_t(x)}{
\sqrt{2\tau\log\log(1/\tau)}} =-\liminf
_{\tau\downarrow0}\frac{u_{t+\tau}(x) - u_t(x)}{
\sqrt{2\tau\log\log(1/\tau)}} = \bigl\llvert \sigma \bigl(
u_t(x) \bigr)\bigr\rrvert ,
\end{equation}
almost surely.
\end{theorem}

Local iterated logarithm laws, such as \eqref{eq:local:LIL}, are well known
in the context of finite-dimensional diffusions; see, for instance, Anderson
\cite{Anderson}, Theorem~4.1. The time-change methods employed in the
finite-dimensional setting will, however, not work effectively in the present
infinite-dimensional context. Here, we obtain~\eqref{eq:local:LIL} as
a ready
consequence of the proof of the ``random Radon--Nikod\'ym property''~\eqref{eq:local:RN}.

%re2.3 #&#
\begin{remark}
Fix an $x\in\Z^d$ and a $t>0$,
and consider the ratio $R(\tau):=[u_{t+\tau}(x)-u_t(x)]/[B_{t+\tau
}(x)-B_t(x)]$;
this is a well-defined random variable for every $\tau>0$, since
$B_{t+\tau}(x)-B_t(x)\neq0$ with probability one for every $\tau>0$.
However, $\{R(\tau)\}_{\tau>0}$ is not a well-defined stochastic process
since there exist random times $\tau>0$ such that
$B_{t+\tau}(x)-B_t(x)=0$ a.s. Thus one does not expect that
the mode of convergence in \eqref{eq:local:RN} can be improved to almost-sure
convergence. This statement can be strengthened further still, but we
will not do so
here.%\qed
\end{remark}

%re2.4 #&#
\begin{remark}
According to \eqref{eq:local:RN}, the solution to the \eqref{SHE}
behaves as the noninteracting system ``$\d u_t(x)\approx\sigma
(u_t(x))\,\d B_t(x)$''
of diffusions, locally to first order. This might seem to
suggest the [false] assertion
that $x\mapsto u_t(x)$ ought to be a sequence of independent
random variables. This is not true, as can be seen by looking more
closely at the time increments of $t\mapsto u_t(x)$. In fact,
our arguments can be extended to show that the spatial correlation
structure of
$u$ appears at second-order approximation levels in the
sense of the following three-term stochastic Taylor expansion:
in the scale $\tau^{1/2}$:
%
%e2.5 #&#
\begin{equation}
\label{eq:Taylor} u_{t+\tau}(x) \simeq u_t(x) +
\tau^{1/2} \sigma\bigl(u_t(x)\bigr) Z_1 +\tau
Z_2 +\tau^{3/2} U(\tau) \qquad\mbox{as $\tau\downarrow0$},
\end{equation}
where: (i) ``$\simeq$'' denotes approximation of
distributions; (ii) $Z_1$ is a standard normal variable
independent of $u_t(x)$;
(iii) $Z_2$ is a nontrivial random variable that depends on
the entire random field $\{u_s(y)\}_{s\in[0,t],y\in\Z^d}$; (iv)
$U(\tau
)=O_\P(1)$ as $\tau\downarrow0$ means that
$\lim_{m\uparrow\infty}\limsup_{\tau\downarrow0}\P\{|U(\tau
)|\ge m\}=0$.
In particular, \eqref{eq:Taylor} tells us that the temporally-local
interactions in the random field $x\mapsto u_t(x)$ are second order in nature.
\end{remark}

Rather than prove these refined assertions, we next turn our attention
to a different
local property of the solution to \eqref{SHE} and show that, after a scale change,
the local-in-time behavior of the solution to \eqref{SHE} is that of spatial
white noise.

%th2.5 #&#
\begin{theorem}\label{th:CLT}
Suppose $\sigma(z)>0$ for all $z\in\R\setminus\{0\}$, and define
%
%e2.6 #&#
\begin{equation}
S(z) := \int_{z_0}^z \frac{\d w}{\sigma(w)} \qquad(z\ge0),
\end{equation}
where $z_0\in\R\setminus\{0\}$ is a fixed number. Then,
$S(u_t(x))<\infty$ a.s. for all $t>0$ and $x\in\Z^d$.
Furthermore, if we choose and fix
$m$ distinct points $x_1,\ldots,x_m\in\Z^d$, then for all $t> 0$
and $q_1,\ldots,q_m\in\R$,
%
%e2.7 #&#
\begin{equation}
\lim_{\tau\downarrow0} \P \Biggl(\bigcap_{j=1}^m
\bigl\{ S \bigl( u_{t+\tau}(x_j) \bigr) - S \bigl(
u_t(x_j) \bigr) \le q_j\sqrt\tau \bigr\}
\Biggr) = \prod_{j=1}^ m
\Phi(q_j),
\end{equation}
where $\Phi(q) :=(2\pi)^{-1/2}\int_{-\infty}^q
\exp(-w^2/2) \,\d w$
denotes the standard Gaussian cumulative distribution function.
\end{theorem}

The preceding manifests itself in interesting ways for different
choices of
the nonlinearity coefficient $\sigma$. Let us mention the following parabolic
Anderson model, which has been a motivating example for us.

%ex2.6 #&#
\begin{example}
Consider the semi-discrete
parabolic Anderson model, which is \eqref{SHE} with
$\sigma(x)\equiv qx$ [for some fixed constant $q> 0$]. In that case,
the solution to \eqref{SHE} is positive and
the ``scale function'' $S$ is
$S(z) = q^{-1}\ln(z/z_0)$ for $z,z_0>0$. As such,
$\sigma(u_t(x))=qu_t(x)$ in \eqref{SHE}, and we find the following
log-normal limit law: For every $t>0$ and
$x_1,\ldots,x_m\in\Z^d$ fixed,
%
%e2.8 #&#
\begin{equation}\qquad
\biggl( \biggl[\frac{u_{t+\tau}(x_1)}{u_t(x_1)} \biggr]^{1/\sqrt{\tau}},\ldots, \biggl[
\frac{u_{t+\tau}(x_m)}{u_t(x_m)} \biggr]^{1/\sqrt{\tau}} \biggr)
\Rightarrow\bigl({\mathrm{e}}^{qN_1},\ldots,{\mathrm{e}}^{qN_m}\bigr)
\qquad\tau \downarrow0,
\end{equation}
where $N_1 ,\ldots,N_m$ are
i.i.d. standard normal variables, and ``$\Rightarrow$'' denotes convergence
in distribution.
\end{example}

Our final main result is a statement about the large-time behavior of the
solution $u$ to \eqref{SHE}. We prove a rigorous version of the
following assertion: ``\emph{If the random walk $X$ is transient and
$\lip$ is sufficiently small---so that \textup{\eqref{SHE}} is not very noisy---then
a decay condition such as $u_0\in\ell^1(\Z^d)$ on the initial
profile is enough to ensure that
$\sup_{x\in\Z^d}|u_t(x)|\to0$ almost surely as $t\to\infty$.}''
This is new even for the parabolic Anderson model,
where $\sigma(x)\propto x$ and $\mathscr{L}:=$ the generator of
the simple walk on $\Z^d$. In fact, this result gives a partial [though
strong] negative answer to an open problem of Carmona and Molchanov
\cite{CM94}, page 122, and rules out the existence of [the analogue of]
a nontrivial ``Anderson mobility edge'' in the present nonstationary setting,
when $u_0\in\ell^1(\Z^d)$.

Recall that $X:=\{X_t\}_{t\ge0}$ is a continuous-time random walk
on $\Z^d$ with generator $\mathscr{L}$. Let $X'$ denote an independent
copy of $X$, and define
%
%e2.9 #&#
\begin{equation}
\label{Upsilon(0)} \Upsilon(0) := \int_0^\infty\P
\bigl\{X_t=X_t'\bigr\} \,\d t = \E\int
_0^\infty \1 _{\{0\}}
\bigl(X_t-X_t'\bigr) \,\d t.
\end{equation}
We can think of $\Upsilon(0)$ as the expected value of
the total occupation time of $\{0\}$, as viewed by the symmetrized random
walk $X-X'$.
Although $\Upsilon(0)$ is always well defined, it is finite
if and only if the symmetrized random walk $X-X'$ is transient~\cite{ChungFuchs}.
We are ready to state our final result.
%\begin{theorem}\label{th:as2}
% Suppose that there exists a real number $k>4$ so that : (i)
% \begin{equation}\label{cond:zk2}
% \lip< \frac{1}{z_k\sqrt{\Upsilon(0)}},
% \end{equation}
% where $1/\infty:=0$; and (ii)
% \begin{equation}\label{cond:alpha2}
% \limsup_{t\to\infty}\
% t^\alpha\P\{X_t=X'_t\}<\infty
% \mbox{for some $\alpha>\max\left(1 ,\frac{2}{k-4}\right);$}
% \end{equation}
% where $X'$ denotes an independent copy of $X$.
% If, in addition, $u_0\in\ell^1(\Z^d)$ is non negative, $\sigma(0)=0$,
% and the underlying probability space is complete, then
% \begin{equation}
% \lim_{t\to\infty}\sum_{x\in\Z^d}|u_t(x)|^2 = 0
% \mbox{almost surely}.
% \end{equation}
%\end{theorem}

%For the special case of walks with zero mean and finite variance, we
%have an optimal result.

%th2.7 #&#
\begin{theorem}[(Dissipation)]\label{th:as}
%Let $\mathscr{L}$ be the generator of a mean zero, finite variance
%random walk in transient dimensions ($d\ge3$) and
Suppose that %there exists a real number $k>4$ so that : (i)
%
%e2.10 #&#
\begin{equation}
\label{cond:zk} \lip< \bigl[\Upsilon(0) \bigr]^{-1/2},
\end{equation}
and that there exists $\alpha\in(1 ,\infty)$ such that
%
%e2.11 #&#
\begin{equation}
\label{cond:alpha} \P\bigl\{X_t=X'_t\bigr\}=
O\bigl(t^{-\alpha}\bigr) \qquad(t\to\infty),
\end{equation}
%
%and
%where $1/\infty:=0$; and (ii)
%\begin{equation}\label{cond:alpha}
% \limsup_{t\to\infty}\
% t^\alpha\P\{X_t=X'_t\}<\infty
% \mbox{for some $\alpha>1$};
%\end{equation}
%where $X'$ denotes an independent copy of $X$.
where $X'$ denotes an independent copy of $X$.
If, in addition, $u_0\in\ell^1(\Z^d)$
and the underlying probability space is complete, then
%
%e2.12 #&#
\begin{equation}
\lim_{t\to\infty}\sup_{x\in\Z^d}\bigl|u_t(x)\bigr|=
\lim_{t\to\infty}\sum_{x\in\Z^d}\bigl|u_t(x)\bigr|^2
= 0\qquad \mbox{almost surely}.
\end{equation}
\end{theorem}

%re2.8 #&#
\begin{remark}[(Hysteresis)]\label{rem:Hysteresis}
Consider the parabolic
Anderson model\break [$\sigma(x)\propto x$], where the underlying
symmetrized walk $X-X'$ is transient, the noise level
is small and $u_0$ is a constant. It is well known that
under these conditions
$u_t(x)$ converges weakly as $t\to\infty$ to a nondegenerate
random variable $u_\infty(x)$ for every $x\in\Z^d$.
See, for example, Greven and den Hollander \cite{GrevenDenHollander}, Theorem~1.4,
Cox and Greven \cite{CoxGreven}, and Shiga \cite{Shiga}. These results
provide a partial affirmative answer to a question of Carmona and Molchanov
\cite{CM94}, page 122, about the existence of long-term invariant laws
in the low-noise regime
of the transient parabolic Anderson model, in particular.
By contrast, Theorem~\ref{th:as} shows
that if $u_0$ is far from stationary (here, it decays at infinity),
then the system is very strongly dissipative
in the low-noise regime. This
result implies that \emph{the parabolic Anderson model
remembers its initial state forever}.%\qed
\end{remark}

%ex2.9 #&#
\begin{example}
Continuous-time walks that have property \eqref{cond:alpha} include
all transient finite-variance centered random walks on $\Z^d$
[$d>2$, necessarily]. For those walks, $\alpha:=d/2$,
thanks to the local central limit theorem. There are more interesting examples
as well. For instance, suppose $t^{-1/p}X_t$ converges in distribution
to a stable random variable $S$ as $t\to\infty$; see
Gnedenko and Kolmogorov~\cite{GK}, Section~35, for necessary and sufficient
conditions. Then $S$ is necessarily stable with index $p$, $p\in(0 ,2]$,
and $t^{-1/p}(X_t-X_t')$ converges in law to a symmetric stable
random variable $S$ with stability index $p$.
If, in addition, the group of all possible values of
$X_t-X_t'$ generates all of $\Z^d$, then a theorem of Gnedenko
\cite{GK}, page 236, ensures that $t^{1/p}\P\{X_t=X_t'\}$ converges
to $f(0)<\infty$, where $f$ denotes the probability density function
of $S$,
as long as $p\in(0 ,1)$.
%\qed
\end{example}

\textit{Organization of the paper}.
In Section~\ref{sec:prelim}
we introduce the mild solution to \eqref{SHE} and state the version of
Burkholder--Davis--Gundy
inequality that we will use throughout the paper.

In Section~\ref{sec:thm1.1}
we prove Theorem~\ref{th:exist:unique}: Section~\ref{subsec:upper lyap}
includes proof of upper bounds for the Lyapunov exponents,
from which existence of the solution follows;
Section~\ref{subsec:lower lyap} contains a comparison principle for \eqref{SHE},
together with proof of lower bounds for the Lyapunov exponents.

Section~\ref{sec:approx}
has some results on the local (in time) behavior of the solution.
These results
are used in Section~\ref{sec:thm1.2} in order to prove Theorem~\ref
{th:local:RadonNikodym}.

Section~\ref{sec:thm1.4} contains a proof of Theorem~\ref{th:CLT}.
That proof hinges on Theorem~\ref{th:local:RadonNikodym} and the fact
that the solution to
\eqref{SHE} immediately becomes strictly positive everywhere (Proposition~\ref
{prop:pos}).

Section~\ref{sec:prelim thm1.6} explains how $\|u_t\|^2_{\ell^{2}(\Z^d)}$
is connected to the intersection local times of two independent
continuous-%
time random walks with common generator $\mathscr{L}$.
The results of Section~\ref{sec:prelim thm1.6} are then used
in Section~9 in order to prove Theorem~\ref{th:as}.

%s3 #&#
\section{Preliminaries} \label{sec:prelim}
In this small section we collect some preliminary facts about SPDEs interacting
diffusion processes and BDG-type martingale inequalities. These facts are
used throughout the rest of the paper.

%s3.1 #&#
\subsection{The mild solution} \label{subsec:mild}

Recall that the convolution on $\Z^d$ is defined by
%
%e3.1 #&#
\begin{equation}
(f*g) (x) := \sum_{y\in\Z^d} f(x-y)g(y) \qquad\bigl(x\in
\Z^d\bigr).
\end{equation}
For every function $h\dvtx\Z^d\to\R$
we define a new function $\tilde{h}$,
%
%e3.2 #&#
\begin{equation}
\tilde{h}(x) := h(-x)\qquad \bigl(x\in\Z^d\bigr),
\end{equation}
as the \emph{reflection} of $h$.

By a ``solution'' to \eqref{SHE} we mean a solution in integrated---%
or ``mild''---form. That is,
a predictable process $t\mapsto u_t$,
with values in $\R^{\Z^d}$, that solves the following infinite system
of It\^o SDEs:
%
%e3.3 #&#
\begin{equation}
\label{mild} u_t(x) = (\tilde{p}_t*u_0)
(x) + \sum_{y\in\Z^d}\int_0^t
p_{t-s}(y-x)\sigma\bigl(u_s(y)\bigr) \,\d
B_s(y),
\end{equation}
where $p_t(x) := \P\{X_t=x\}$.

It might be helpful to note also that
$(P_t\phi)(x) :=  (\tilde{p}_t*\phi )(x)$ defines the semigroup
of the random walk $X$ via the identity $(P_t\phi)(x) =
\E\phi(x+X_t)$. Thus we can write \eqref{mild} in the following,
perhaps more familar, form:
%
%e3.4 #&#
\begin{equation}
u_t(x) = (P_tu_0) (x) + \sum
_{y\in\Z^d}\int_0^t
p_{t-s}(y-x)\sigma\bigl(u_s(y)\bigr) \,\d
B_s(y).
\end{equation}

%s3.2 #&#
\subsection{A BDG inequality} \label{subsec:BDG}

We begin this subsection with some background on Burkholder's constants
which will give us the best constants in Lemma~\ref{lem:BDG}. According
to the Burkholder--Davis--Gundy inequality
\cite{Burkholder,BDG,BG},
%
%e3.5 #&#
\begin{equation}
z_p := \sup_x \sup_{t>0}
\biggl[\frac{\E (|x_t|^p )}{
\E (\langle x\rangle_t^{p/2} )} \biggr]^{1/p}<\infty,
\end{equation}
where the supremum ``$\sup_x$'' is taken over all nonzero martinagles
$x:=\{x_t\}_{t\ge0}$ that have continuous trajectories and are
in $L^2(\P)$ at all times, $\langle x\rangle_t$ denotes the quadratic
variation of $x$ at time $t$ and $0/0:=\infty/\infty:=0$.
Davis \cite{Davis} has computed the numerical value of $z_p$
in terms of zeroes of special functions. When $p\ge2$ an integer,
Davis's theorem implies that $z_p$
is equal to the largest positive root of the modified Hermite polynomial
$\mbox{\sl He}_p$.
Thus, for example, we obtain the following from direct evaluation of
the zeros:
%
%e3.6 #&#
\begin{eqnarray}
z_2 &= &1, \qquad z_3=\sqrt3, \qquad z_4=
\sqrt{3+\sqrt{6}} \approx2.334,
\nonumber
\\[-8pt]
\\[-8pt]
\nonumber
z_5&=&\sqrt{5+\sqrt{10}}\approx2.857,\qquad  z_6
\approx3.324,\ldots.
\end{eqnarray}

It is known that $ z_p\sim2\sqrt p$ as $p\to\infty$, and
$\sup_{p\ge2} ( z_p/\sqrt p)=2$;
see Carlen and Kr{\'e}e \cite{CK}, Appendix.%\qed

Suppose $Z:=\{Z_t(x)\}_{t\ge0,x\in\Z^d}$ is a predictable random field,
with respect to the infinite-dimensional
Brownian motion $\{B_t(\bullet)\}_{t\ge0}$, that satisfies the moment bound
$\E\int_0^t\|Z_s\|_{\ell^2(\Z^d)}^2 \,\d s<\infty$. Then the It\^o integral
process defined by
%
%e3.7 #&#
\begin{equation}
\int_0^t Z_s\cdot\,\d
B_s:=\sum_{y\in\Z^d}\int
_0^t Z_s(y) \,\d
B_s(y)\qquad (t\ge0)
\end{equation}
exists and defines a continuous $L^2(\P)$ martingale. See, for example,
Pr{\'e}v{\^o}t and R{\"o}ckner
\cite{PrevotRoeckner}.
The following variation of the Burkholder--Davis--Gundy
inequality yields moment bounds for this martingale. %that also pay
%special attention
%to the constants in such inequalities.

%le3.1 #&#
\begin{lemma}[(BDG lemma)]\label{lem:BDG}
For all finite real numbers $k\ge2$ and $t\ge0$,
%
%e3.8 #&#
\begin{equation}
\label{BDG} \E \biggl(\biggl\llvert \int_0^t
Z_s\cdot\,\d B_s \biggr\rrvert ^k \biggr)
\le \biggl|4k\sum_{y\in\Z^d}\int_0^t
\bigl\{ \E \bigl( \bigl|Z_s(y)\bigr|^k \bigr) \bigr
\}^{2/k}\,\d s \biggr|^{k/2}.
\end{equation}
\end{lemma}

\begin{pf}
We follow a method of Foondun and Khoshnevisan \cite{FK}.

A standard approximation argument tells us that it suffices to consider
the case where $y\mapsto Z_s(y)$ has finite support.
Let $F\subset\Z^d$ be a finite set of cardinality \mbox{$m\ge1$},
and suppose $Z_s(y)=0$
for all $y\notin F$. Consider the (standard, finite-dimensional)
It\^o integral process $\int_0^t Z_s\cdot\,\d B_s := \sum_{y\in F}\int_0^t
Z_s(y) \,\d B_s(y)$.
According to Davis's \cite{Davis} form
of the Burkholder--Davis--Gundy inequality
$m$-dimensional Brownian motion \cite{Burkholder,BDG,BG},
%
%e3.9 #&#
\begin{equation}
\label{eq:Davis} \E \biggl( \biggl\llvert \int_0^t
Z_s\cdot\,\d B_s\biggr\rrvert ^k \biggr)
\le z_k^k \E \biggl( \biggl| \sum_{y\in F}
\int_0^t \bigl[Z_s(y)
\bigr]^2 \,\d s \biggr|^{k/2} \biggr).
\end{equation}
Finally, we use the Carlen--Kr{\'e}e bound $z_k\le2\sqrt k$ \cite{CK}
together with the Minkowski inequality to finish the proof in
the case where $F$ is finite. A standard finite-dimensional approximation
completes the proof.
\end{pf}

%s4 #&#
\section{Proof of Theorem \texorpdfstring{\protect\ref{th:exist:unique}}{2.1}} \label{sec:thm1.1}
%s4.1 #&#
\subsection{Bounds for the upper Lyapunov exponents} \label
{subsec:upper lyap}
Existence and uniqueness, and also continuity, of the solution are
dealt with extensively in the literature and are well
known; see, for example, Shiga and Shimizu \cite{ShigaShimizu}
and the general theory of Pr{\'e}v{\^o}t and R{\"o}ckner \cite{PrevotRoeckner}
for some of the latest developments. However,
in order to derive our estimates of the Lyapunov exponents we will need
a priori
estimates which will also yield existence and uniqueness. Therefore,
in this section, we hash out some---though not all---of the details.

Let us proceed by applying Picard iteration. Let
$u^{(0)}_t(x):=u_0(x)$, and
then define iteratively for all $n\ge0$,
%
%e4.1 #&#
\begin{equation}
\label{eq:Picard} u^{(n+1)}_t(x) := (\tilde{p}_t*u_0)
(x) + \sum_{y\in\Z^d}\int_0^t
p_{t-s}(y-x)\sigma \bigl( u^{(n)}_s(y) \bigr) \,\d
B_s(y).
\end{equation}
It follows from the properties of the It\^o integral that
%
%e4.2 #&#
\begin{equation}
\label{eq:Mn} M_t^{(n+1)}:=\sup_{x\in\Z^d} \E
\bigl(\bigl| u_t^{(n+1)}(x) \bigr|^k \bigr)
\le2^{k-1}\sup_{x\in\Z^d}(I_x+J_x),
\end{equation}
where
%
%e4.3 #&#
\begin{eqnarray}
\label{eq:IJ} %
I_x&:=& \bigl\llvert (\tilde
p_t*u_0) (x) \bigr\rrvert ^k,
\nonumber
\\[-8pt]
\\[-8pt]
\nonumber
J_x&:=&\E \biggl( \biggl| \sum_{y\in\Z^d}\int
_0^t p_{t-s}(y-x) \sigma \bigl(
u_s^{(n)}(y) \bigr) \,\d B_s(y)
\biggr|^k \biggr).
\end{eqnarray}
The first term is easy to bound:
%
%e4.4 #&#
\begin{equation}
\label{eq:Ix} \sup_{x\in\Z^d} I_x \le
\|u_0\|_{\ell^\infty(\Z^d)}^k,
\end{equation}
since $\sum_xp_t(x)= 1$. Next we bound $J_x$.

Because $\sigma$ is Lipschitz continuous and $\sigma(0)=0$, we can
see that
$|\sigma(z)|\le\lip|z|$ for all $z\in\R$. Thus
we may use the BDG lemma (Lemma~\ref{lem:BDG}) in order to see that
%
%e4.5 #&#
\begin{equation}
\label{eq:Jx:prec} J_x^{2/k} \le4k\lip^2\sum
_{y\in\Z^d}\int_0^t
\bigl[p_{t-s}(y-x)\bigr]^2 \bigl\{ \E \bigl( \bigl\llvert
u^{(n)}_s(y)\bigr\rrvert ^k \bigr) \bigr
\}^{2/k} \,\d s.
\end{equation}
Therefore, we may recall the inductive
definition \eqref{eq:Mn} of $M$ to see that
%
%e4.6 #&#
\begin{eqnarray}
\label{eq:Jx} %
J_x^{2/k} &\le&4k
\lip^2\sum_{y\in\Z^d}\int_0^t
\bigl[p_{t-s}(y-x)\bigr]^2 \bigl( M^{(n)}_s
\bigr)^{2/k} \,\d s
\nonumber
\\[-8pt]
\\[-8pt]
\nonumber
&\le&4k\lip^2\int_0^t \bigl(
M^{(n)}_s \bigr)^{2/k} \,\d s,
\end{eqnarray}
since
%
%e4.7 #&#
\begin{equation}
\label{eq:replica} \sum_{z\in\Z^d}\bigl[p_r(z)
\bigr]^2=\P\bigl\{X_r=X_r'
\bigr\}\le1,
\end{equation}
where $X'$ denotes
an independent copy of $X$. (This last bound might appear to be quite crude,
and it is when $r$ is large. However, it turns out that the behavior of
$r$ near zero matters more to us. Therefore, the inequality is tight in the
regime $r\approx0$ of interest to us.)

We may combine \eqref{eq:Mn}, \eqref{eq:Ix}
and \eqref{eq:Jx} in order to see that for all $\beta,t>0$,
%
%e4.8 #&#
\begin{eqnarray}\quad
&&{\mathrm{e}}^{-\beta t}M^{(n+1)}_t
\nonumber
\\[-8pt]
\\[-8pt]
\nonumber
&&\qquad\le2^{k-1}\|u_0\|_{\ell^\infty(\Z^d)}^k
+\bigl(16k\lip^2\bigr)^{k/2}\biggl\llvert \int
_0^t{\mathrm{e}}^{-2\beta(t-s)/k} \bigl( {\mathrm{e}}^{-\beta s}M^{(n)}_s \bigr)^{2/k} \,\d s
\biggr\rrvert ^{k/2}.
\end{eqnarray}
Consequently, the sequence defined by
%
%e4.9 #&#
\begin{equation}
N_\beta^{(m)}:= \sup_{t\ge0} \bigl( {\mathrm{e}}^{-\beta
t}M^{(m)}_t \bigr)\qquad (m\ge0)
\end{equation}
satisfies the recursive inequality
%
%e4.10 #&#
\begin{eqnarray}
N_\beta^{(n+1)} &\le&2^{k-1}
\|u_0\|_{\ell^\infty(\Z^d)}^k + \bigl(16k\lip^2
\bigr)^{k/2}\biggl\llvert \int_0^t{\mathrm{e}}^{-2\beta s/k} \,\d s
\biggr\rrvert ^{k/2}N_\beta^{(n)}
\nonumber
\\[-8pt]
\\[-8pt]
\nonumber
&\le&2^{k-1}\|u_0\|_{\ell^\infty(\Z^d)}^k +
\biggl(\frac{8 k^2\lip^2}{\beta} \biggr)^{k/2}N_\beta^{(n)}.
\end{eqnarray}
In particular, if we denote (temporarily for this proof)
%
%e4.11 #&#
\begin{equation}
\alpha:= 8(1+\delta)\lip^2,
\end{equation}
where $\delta>0$ is fixed but arbitrary, then
%
%e4.12 #&#
\begin{equation}
N_{\alpha k^2}^{(n+1)} \le2^{k-1}\|u_0
\|_{\ell^\infty(\Z^d)}^k + (1+\delta)^{-k/2}N_{\alpha k^2}^{(n)}.
\end{equation}
We may apply induction on $n$ now in order to see that
$\sup_{n\ge0}N_{\alpha k^2}^{(n)}<\infty$;
equivalently, for all $k\ge2$ there exists $c_k\in(0 ,\infty)$
such that
%
%e4.13 #&#
\begin{equation}
\label{eq:u:UB} \sup_{x\in\Z^d}\E \bigl( \bigl\llvert
u^{(n)}_t(x) \bigr\rrvert ^k \bigr) \le
c_k {\mathrm{e}}^{8(1+\delta)\lip^2 k^2t} \qquad\mbox{for all $t\ge0$}.
\end{equation}
Similarly,
%
%e4.14 #&#
\begin{eqnarray}
&&\E \bigl( \bigl| u^{(n+1)}_t(x) - u^{(n)}_t(x)
\bigr|^k \bigr)\nonumber
\\
&&\qquad=\E \biggl( \biggl| \sum_{y\in\Z^d}\int
_0^t p_{t-s}(y-x) \bigl\{ \sigma
\bigl( u^{(n)}_s(y) \bigr) - \sigma \bigl(
u^{(n-1)}_s(y) \bigr) \bigr\} \,\d B_s(y)
\biggr|^k \biggr)
\\
\nonumber
&&\qquad\le\bigl(4k\lip^2\bigr)^{k/2} \E \biggl( \biggl| \sum
_{y\in\Z^d}\int_0^t
\bigl[ p_{t-s}(y-x) \bigr]^2 \bigl\{ u^{(n)}_s(y)
- u^{(n-1)}_s(y) \bigr\}^2 \,\d s\biggr |^{k/2}
\biggr).\hspace*{-12pt}
\end{eqnarray}
Define
%
%e4.15 #&#
\begin{equation}
L_t^{(n+1)}:= \sup_{x\in\Z^d} \E \bigl( \bigl
\llvert u^{(n+1)}_t(x) - u^{(n)}_t(x)
\bigr\rrvert ^k \bigr)
\end{equation}
to deduce from the preceding, \eqref{eq:replica} and Minkowski's
inequality that
%
%e4.16 #&#
\begin{eqnarray}
L^{(n+1)}_t &\le&\bigl(4k\lip^2
\bigr)^{k/2} \biggl( \sum_{y\in\Z^d}\int
_0^t \bigl[ p_{t-s}(y-x)
\bigr]^2 \bigl( L_s^{(n)} \bigr)^{2/k}
\,\d s \biggr)^{k/2}
\nonumber
\\[-8pt]
\\[-8pt]
\nonumber
&\le&\bigl(4k\lip^2\bigr)^{k/2} \biggl( \int
_0^t \bigl(L_s^{(n)}
\bigr)^{2/k} \,\d s \biggr)^{k/2}. %
\end{eqnarray}
Therefore,
%
%e4.17 #&#
\begin{equation}
K_{\alpha k^2}^{(m)}= \sup_{t\ge0} \bigl( {\mathrm{e}}^{-\alpha k^2 t} L_t^{(m)} \bigr)
\end{equation}
satisfies
%
%e4.18 #&#
\begin{eqnarray}
\nonumber
K_{\alpha k^2}^{(n+1)} &\le&\bigl(4k\lip^2
\bigr)^{k/2} \biggl( \int_0^t{\mathrm{e}}
^{-2\alpha k(t-s)} \,\d s \biggr)^{k/2} K_{\alpha k^2}^{(n)}\\
&\le&
\biggl( \frac{4\lip^2}
{2\alpha} \biggr)^{k/2} K_{\alpha k^2}^{(n)}
\\
&\le&2^{-k}K_{\alpha k^2}^{(n)}.\nonumber
\end{eqnarray}
From this we conclude that $\sum_{n=0}^\infty K_{\alpha k^2}^{(n)}<
\infty$. Therefore, there exists a random field $u_t(x)$ such that
$\lim_{n\to\infty}u_t^{(n)}(x)=u_t(x)$ in $L^k(\P)$.
It follows readily that $u$ solves \eqref{SHE}, and $u$ satisfies \eqref
{eq:Upper:LE} by
\eqref{eq:u:UB} and Fatou's lemma. Uniqueness is proved by similar means,
and we skip the details.
%\qed

%s4.2 #&#
\subsection{Bounds for the lower Lyapunov exponents} \label
{subsec:lower lyap}
We start the section with a truncation error estimate for the
nonlinearity $\sigma$.
This will be needed to use the results of Cox, Fleischmann and Greven
\cite{CFG}
on comparison of moments for interacting diffusions.
We can then reduce our problem to the case of $\sigma(x)=\ell_\sigma x$.

%le4.1 #&#
\begin{lemma}\label{lem:truncate}
Define $\sigma^{(N)}$ by $\sigma^{(N)}:=\sigma$ on $(-N ,N)$,
$\sigma^{(N)}:=0$ on $[-N-1 ,N+1]^c$, and defined by linear interpolation
on $[-N-1 ,-N]\cup[N ,N+1]$.
Let $U^{(N)}_t(x)$ denote the a.s.-unique solution to
\mbox{\eqref{SHE}} where $\sigma$ is replaced by
$\sigma^{(N)}$.
Then,
$\lim_{N\to\infty} U^{(N)}_t(x)=u_t(x)$ a.s. and
in $L^k(\P)$ for all
$k\ge2$, $t\ge0$ and $x\in\Z^d$.
\end{lemma}

\begin{pf}
Since $\sigma^{(N)}$ is Lipschitz continuous, Theorem~\ref{th:exist:unique}
ensures the existence and uniqueness of $U^{(N)}$ for every $N\ge1$.
Then by \eqref{mild}
%
%e4.19 #&#
\begin{equation}
u_t(x) - U^{(N)}_t(x) = T_1 +
T_2,
\end{equation}
where
%
%e4.20 #&#
\begin{eqnarray}
T_1 &:=&\sum_{y\in\Z^d}\int
_0^t p_{t-s}(y-x) \bigl\lbrace
\sigma\bigl( u_s(y)\bigr)-\sigma^{(N)}(u_s(y)
\bigr\rbrace\,\d B_s(y);
\nonumber
\\[-8pt]
\\[-8pt]
\nonumber
T_2 &:=& \sum_{y\in\Z^d}\int
_0^t p_{t-s}(y-x) \bigl\{
\sigma^{(N)} \bigl( u_s(y)\bigr) - \sigma^{(N)}
\bigl( U^{(N)}_s(y) \bigr) \bigr\} \,\d B_s(y).
\end{eqnarray}
Because $|\sigma(z)|\le\lip|z|$, Lemma~\ref{lem:BDG} implies
that $ \{\E(|T_1|^k) \}^{2/k}$ is at most
%
%e4.21 #&#
\begin{equation}
4k\lip^2\sum_{y\in\Z^d}\int
_0^t \bigl[p_{t-s}(y-x)
\bigr]^2 \bigl\{ \E \bigl(\bigl\llvert u_s(y)\bigr
\rrvert ^k; \bigl|u_s(y)\bigr|\ge N \bigr) \bigr\}^{2/k}
\,\d s.
\end{equation}
We have $\E( |Y|^k; |Y|\ge N ) \le N^{-k}\E(Y^{2k} )$,
valid for all $Y\in L^{2k}(\Omega)$. Therefore,
%
%e4.22 #&#
\begin{equation}\quad
\bigl\{\E\bigl(|T_1|^k\bigr) \bigr\}^{2/k}\le
\frac{4k\lip^2}{N^2} \sum_{y\in\Z^d}\int
_0^t \bigl[p_{t-s}(y-x)
\bigr]^2 \bigl\{ \E \bigl(\bigl\llvert u_s(y)\bigr
\rrvert ^{2k} \bigr) \bigr\}^{2/k} \,\d s.
\end{equation}
Because $\sum_{y\in\Z^d}[p_{t-s}(y-x)]^2\le1$---see \eqref
{eq:replica}---%
the already-proved
bound \eqref{eq:Upper:LE} tells us that
%
%e4.23 #&#
\begin{equation}
\label{eq:T_1} \bigl\{\E\bigl(|T_1|^k\bigr) \bigr
\}^{2/k} \le\frac{a_k}{N^2} \int_0^t{\mathrm{e}}^{128\lip^2 k s} \,\d s \le\frac{Aa_k{\mathrm{e}}^{A kt}}{N^2},
\end{equation}
where $a_k$ and $A$ are uninteresting finite and positive constants;
moreover, $a_k$ depends only on $k$. This estimates the norm of
$T_1$.

As for $T_2$, we use the simple inequality
$|\sigma^{(N)}(r)-\sigma^{(N)}(\rho)|\le C |r-\rho|$,
together with the BDG Lemma~\ref{lem:BDG} in order to find that
%
%e4.24 #&#
\begin{equation}
 \bigl\{ \E \bigl(|T_2|^k \bigr) \bigr
\}^{2/k} \le b_k\int_0^t
\sup_{y\in\Z^d} \bigl\{\E \bigl(\bigl\llvert u_s(y)-
U^{(N)}_s(y)\bigr\rrvert ^k \bigr) \bigr
\}^{2/k} \,\d s, %
\end{equation}
where $b_k$ is a constant dependent on $\sigma$ and $k$.
Together, the preceding moment bounds for $T_1$ and $T_2$ imply that
%
%e4.25 #&#
\begin{equation}
D^{(N)}_t := \sup_{x\in\Z^d} \bigl\{\E
\bigl(\bigl\llvert u_t(x)-U^{(N)}_t(x)\bigr
\rrvert ^k \bigr) \bigr\}^{2/k}
\end{equation}
satisfies the recursion
%
%e4.26 #&#
\begin{equation}
D^{(N)}_t \le\frac{\tilde{a}_k{\mathrm{e}}^{\tilde{A}k^2t}}{N^2} + \tilde{b}_k
\int_0^t D^{(N)}_s \,\d
s,
\end{equation}
where $\tilde{a}_k$, $\tilde{b}_k$ and $\tilde{A}$ are positive and
finite constants, and the first two depend only on $k$ (whereas the latter
is universal). An application of the Gronwall
inequality shows that $\sup_{t\in[0,T]}D^{(N)}_t = O(N^{-2})$
as $N\to\infty$, for every fixed value $T\in(0 ,\infty)$. This
is enough to yield the lemma.
\end{pf}

We complete the proof of Theorem~\ref{th:exist:unique} by verifying the
two remaining assertions of that theorem:
(i) The solution is nonnegative because $u_0(x)\ge0$
and $\sigma(0)=0$; and (ii) The lower bound \eqref{eq:Lower:LE}
for the lower Lyapunov exponent holds. We keep the two parts
separate, as they use different ideas.

%th4.2 #&#
\begin{theorem}[(Comparison principle)]\label{th:comparison}
Suppose $u$ and $v$ are the solutions to \mbox{\eqref{SHE}}
with respective initial functions $u_0$ and $v_0$. If $u_0(x)\ge v_0(x)$
for all $x\in\Z^d$, then $u_t(x)\ge v_t(x)$ for all $t\ge0$
and $x\in\Z^d$ a.s.
\end{theorem}

The nonnegativity assertion of Theorem~\ref{th:exist:unique}
is well known \cite{Shiga}, but also
follows from the preceding comparison principle.
This is because condition \eqref{eq:sigma(0)=0} implies
that $v_t(x)\equiv0$ is the unique solution to \eqref{SHE} with initial
condition $v_0(x)\equiv0$. Therefore, the comparison principle yields
$u_t(x)\ge v_t(x)=0$ a.s.

\begin{pf*}{Proof of Theorem~\ref{th:comparison}}
Consider the following infinite dimensional SDE:
%
%e4.27 #&#
\begin{eqnarray}\label{eq:infinitedimSHE}
\qquad w_t(x)=w_0(x)+\int_{0}^{t}(
\mathscr{L}w_s) (x) \,\d s + \int_{0}^{t}
\sigma\bigl(w_s(x)\bigr)\,\d B_s(x)
\nonumber
\\[-8pt]
\\[-8pt]
\eqntext{\bigl(x\in
\Z^d\bigr).}
\end{eqnarray}
It is a well-known fact that the mild solution to \eqref{SHE} is also a
solution in the weak sense. See, for example, Theorem~3.1 of Iwata
\cite{Iwata} and its proof.
Therefore, $u_t(x)$ and $v_t(x)$, respectively, solve \eqref{eq:infinitedimSHE}
with initial conditions $u_0(x)$ and $v_0(x)$.

Let $\{S_n\}_{n=1}^\infty$
denote a growing sequence of finite subsets of $\Z^d$ that exhaust all
of $\Z^d$. Consider, for every $n\ge1$, the stochastic integral equation,
%
%e4.28 #&#
\begin{eqnarray}
\cases{\displaystyle u_t^{(n)}(x)=u_0(x)+\int
_{0}^{t}\bigl(\mathscr{L}u_s^{(n)}
\bigr) (x) \,\d s\vspace*{2pt}\cr
{}\hspace*{42pt} +\displaystyle \int_{0}^{t}\sigma
\bigl(u_s^{(n)}(x)\bigr)\,\d B_s(x),&\quad $\mbox{if
$x\in S_n$}$;\vspace*{2pt}
\cr
u_t^{(n)}(x)=u_0(x),&\quad $
\mbox{if $x\notin S_n$}$.}
\end{eqnarray}
Similarly, we let $v^{(n)}$ solve the same equation, but start it from $v_0(x)$.

Each of these equations is in fact a finite-dimensional SDE, and has a unique
strong solution, by It\^o's theory. Moreover,
Shiga and Shimizu's proof of their Theorem~2.1 \cite{ShigaShimizu}
shows that,
for every $x\in\Z^d$ and $t>0$,
there exists a subsequence $\{n_k\}_{k=1}^\infty$ of increasing integers
such that
%
%e4.29 #&#
\begin{equation}
u^{(n_k)}_t(x) \stackrel{\P} {\longrightarrow}u_t(x)
\quad\mbox{and}\quad v^{(n_k)}_t(x) \stackrel{\P} {
\longrightarrow}v_t(x),
\end{equation}
as $k\to\infty$.
Therefore, we may appeal to a comparison principle for finite-dimensional
SDEs, such as that of Geiss and Manthey \cite{GM}, Theorem~1.2,
in order to conclude the result; the quasi-monotonicity condition of
\cite{GM} is met simply because $\mathscr{L}$ is the generator of a
Markov chain.
The verification of that detail is left to the interested reader.
\end{pf*}

We are now in position to establish the lower bound \eqref{eq:Lower:LE}
on the lower Lyapunov exponent of the solution to \eqref{SHE}.

\begin{pf*}{Proof of Theorem~\ref{th:exist:unique}: Verification of
\eqref{eq:Lower:LE}}
Let $v$ solve the stochastic heat equation
%
%e4.30 #&#
\begin{equation}
\d v_t(x) = (\mathscr{L}v_t) (x) \,\d t +
\ell_\sigma v_t(x) \,\d B_t(x),
\end{equation}
subject to $v_0(x) := u_0(x)$. Also define
$V^{(N)}$ to be the solution to
%
%e4.31 #&#
\begin{equation}
\d V^{(N)}_t(x) = (\mathscr{L}v_t) (x) \,\d t +
\zeta^{(N)} \bigl(V^{(N)}_t(x) \bigr) \,\d
B_t(x),
\end{equation}
where $\zeta^{(N)}(x) :=\ell_\sigma x$ on $(-N ,N)$, $\zeta^{(N)}(x):=0$
when $|x|\ge N+1$, and $\zeta^{(N)}$ is defined by
linear interpolation everywhere else.

Define $\sigma^{(N)}$ and $U^{(N)}$ as in Lemma~\ref{lem:truncate}.
Because $\sigma^{(N)}\ge\zeta^{(N)}$ everywhere on $\R_+$,
and since both $U^{(N)}$ and $V^{(N)}$ are $\ge0$ a.s. and pointwise,
the comparison theorem of Cox, Fleischmann and Greven \cite{CFG}, Theorem~1,
shows us that
%
%e4.32 #&#
\begin{equation}
\E \bigl(\bigl\llvert V^{(N)}_t(x)\bigr\rrvert
^k \bigr) \le\E \bigl( \bigl\llvert U^{(N)}_t(x)
\bigr\rrvert ^k \bigr),
\end{equation}
for all $t\ge0$, $x\in\Z^d$, $k\ge2$ and $N\ge1$. Let $N\to\infty$,
and apply Lemma~\ref{lem:truncate} to find that
$V^{(N)}_t(x)\to v_t(x)$ and $U^{(N)}_t(x)\to u_t(x)$ in $L^k(\P)$ for
all $k\ge2$.
As a result, one can let $N\to\infty$ in the preceding display in
order to deduce
the following:
%
%e4.33 #&#
\begin{equation}
\E \bigl(\bigl\llvert v_t(x)\bigr\rrvert ^k \bigr) \le
\E \bigl( \bigl\llvert u_t(x) \bigr\rrvert ^k \bigr).
\end{equation}
Therefore, it remains to bound $\underline{\gamma}_k(v)$ from below.

Let $\{X^{(i)}\}_{i=1}^k$ denote $k$ independent copies of the random walk
$X$. It is possible to prove that
%
%e4.34 #&#
\begin{equation}
\label{eq:conus} \E \bigl( \bigl|v_t(x)\bigr|^k \bigr) = \E \Biggl(
\prod_{j=1}^k u_0 \bigl(
X^{(j)}_t+x \bigr) \,\cdot\,{\mathrm{e}}^{M_k(t)} \Biggr),
\end{equation}
where $M_k(t)$ denotes the ``multiple collision local time,''
%
%e4.35 #&#
\begin{equation}
M_k(t) := 2\ell_\sigma^2 \mathop{\sum
\sum}_{1\le i<j\le k} \int_0^t
\1_{\{0\}} \bigl( X^{(i)}_s-X^{(j)}_s
\bigr) \,\d s.
\end{equation}
When $X$ is the continuous-time
simple random walk on $\Z^d$, this is a well-known
consequence of a Feynman--Kac formula; see, for instance,
Carmona and Molchanov \cite{CM94}, page 19. When $X$ is replaced by a
L\'evy process, Conus \cite{Conus} has found an elegant derivation of
this formula.
The class of all L\'evy processes includes that of continuous-time
random walks,
whence follows \eqref{eq:conus}.

Note that a.s. on the event that none of the walks $X^{(1)},\ldots,X^{(k)}$
jump in the time interval $[0 ,t]$,
%
%e4.36 #&#
\begin{equation}
\prod_{j=1}^k u_0 \bigl(
X^{(j)}+x \bigr) {\mathrm{e}}^{%
M_k(t)} \ge\bigl[u_0(x)
\bigr]^k {\mathrm{e}}^{k(k-1)\ell_\sigma^2 t}.
\end{equation}
Since the probability is $\exp(-t)$ that $X^{(j)}$ does not jump in
$[0 ,t]$,
it follows from the independence of $X^{(1)},\ldots,X^{(k)}$ that
%
%e4.37 #&#
\begin{equation}
\E \bigl(\bigl|v_t(x)\bigr|^k \bigr) \ge\bigl[u_0(x)
\bigr]^k\exp \bigl\{ \bigl[ k(k-1) \ell_\sigma^2
- k \bigr]t \bigr\}.
\end{equation}
Because $u_0$ is not identically zero, it follows that
%
%e4.38 #&#
\begin{equation}
\underline{\gamma}_k(u) \ge\underline{\gamma}_k(v) \ge
k(k-1) \ell_\sigma^2 - k.
\end{equation}
The preceding is $\ge(1-\varepsilon)
k^2\ell_\sigma^2$ when
$k\ge\varepsilon^{-1}+(\varepsilon\ell_\sigma^2)^{-1}$.
This completes the proof of the theorem.
\end{pf*}

%s5 #&#
\section{A local approximation theorem} \label{sec:approx}
In this section we develop a description of the local dynamics of
the random field $t\mapsto u_t(\bullet)$ in the form of several
approximation results.

Our first approximation lemma is a standard sample-function
continuity result; it states basically that outside a single null set,
%
%e5.1 #&#
\begin{equation}\qquad
\label{eq:local:modulus} u_{t+\tau}(x)= u_t(x) + O \bigl(
\tau^{(1+o(1))/2} \bigr)\qquad \mbox{as $\tau\to0$, for all $t\ge0$ and $x\in
\Z^d$}.
\end{equation}
The result is well known, but we need to be cautious with various constants
that crop up in the proof. Therefore, we include the details to account
for the
dependencies of the implied constants.

%le5.1 #&#
\begin{lemma}\label{lem:modulus}
There exists a version of $u$ that is a.s. continuous in $t$ with
critical H\"older exponent $\ge\nicefrac12$. In fact,
for every $T\ge1$, $\varepsilon\in(0 ,1)$ and $k\ge2$,
%
%e5.2 #&#
\begin{equation}
\sup_{x\in\Z^d}\sup_I \E \biggl( \mathop{
\sup_{s,t\in I}}_{
s\neq t} \biggl[ \frac{|u_t(x) - u_s(x)|}{|t-s|^{(1-\varepsilon)/2}}
\biggr]^k \biggr)<\infty,
\end{equation}
where ``$\sup_I$'' denotes the supremum over all closed subintervals
$I$ of
$[0 ,T]$ that have length $\le1$.
\end{lemma}

\begin{pf}
Minkowski's inequality gives
%
%e5.3 #&#
\begin{equation}
\label{QQQ} \bigl[ \E \bigl( \bigl\llvert u_{t+\tau}(x) -
u_t(x)\bigr\rrvert ^k \bigr) \bigr]^{1/k}
\le|Q_1| + Q_2 + Q_3,
\end{equation}
where
%
%e5.4 #&#
\begin{eqnarray}
\label{eq:Q_1-Q_3}
\nonumber
Q_1 &:=& (\tilde{p}_{t+\tau}*u_0)
(x) - (\tilde{p}_t*u_0) (x),
\\
\qquad Q_2 &:=& \biggl[ \E \biggl(\biggl | \sum
_{y\in\Z^d}\int_0^t \bigl[
p_{t+\tau-s}(y-x)
\nonumber
\\[-10pt]
\\[-10pt]
\nonumber
&&\hspace*{58pt}{} - p_{t-s}(y-x) \bigr] \sigma
\bigl(u_s(y)\bigr) \,\d B_s(y) \biggr|^k \biggr)
\biggr]^{1/k} ,
\\
Q_3 &:=& \biggl[ \E \biggl( \biggl| \sum_{y\in\Z^d}
\int_t^{t+\tau} p_{t+\tau-s}(y-x) \sigma
\bigl(u_s(y)\bigr) \,\d B_s(y)\biggr |^k \biggr)
\biggr]^{1/k}.\nonumber
\end{eqnarray}
We estimate each item in turn.

Let $J_{t,t+\tau}$ denote the event that the random walk $X$ jumps
some time during the time interval $(t ,t+\tau)$. Because
%
%e5.5 #&#
\begin{eqnarray}
\label{eq:p-p} %
\sum_{x\in\Z^d}\bigl\llvert
p_{t+\tau}(x) - p_t(x) \bigr\rrvert &=&\sum
_{x\in\Z^d}\bigl\llvert \E ( \1_{\{X_{t+\tau}=x\}} -
\1_{\{
X_t=x\}
}; J_{t,t+\tau} )\bigr\rrvert
\nonumber
\\[-8pt]
\\[-8pt]
\nonumber
&\le& 2\P(J_{t,t+\tau}) = 2 \bigl(1-{\mathrm{e}}^{-\tau} \bigr) \le2\tau,
\end{eqnarray}
we obtain the following estimate for $|Q_1|$:
%
%e5.6 #&#
\begin{equation}
\label{eq:Q_1} |Q_1| \le2\|u_0\|_{\ell^\infty(\Z^d)}\tau.
\end{equation}

By BDG Lemma~\ref{lem:BDG},
%
%e5.7 #&#
\begin{eqnarray}
\label{eq:Q2:prec}
Q_2^2 &\le&4k\sum
_{y\in\Z^d}\int_0^t \bigl[
p_{t+\tau-s}(y-x) - p_{t-s}(y-x) \bigr]^2 \bigl\{ \E
\bigl( \bigl\llvert \sigma\bigl(u_s(y)\bigr)\bigr\rrvert
^k \bigr) \bigr\}^{2/k} \,\d s
\nonumber
\\[-8pt]
\\[-8pt]
\nonumber
&\le&4k \int_0^t \mathscr{Q}(s)\sup
_{y\in\Z^d} \bigl\{ \E \bigl( \bigl\llvert \sigma
\bigl(u_s(y)\bigr)\bigr\rrvert ^k \bigr) \bigr
\}^{2/k} \,\d s,
\end{eqnarray}
where
%
%e5.8 #&#
\begin{equation}
\mathscr{Q}(s) := \sum_{z\in\Z^d}\bigl|p_{t+\tau-s}(z)-p_{t-s}(z)\bigr|^2
\qquad (0<s<t).\vadjust{\goodbreak}
\end{equation}
%
%It is possible to find a real-variable estimate for $\mathscr{Q}(s)$
%using \eqref{eq:p-p};
%namely, $\mathscr{Q}(s) \le\sum_{z\in\Z^d}|p_{t+\tau-s}(z)-p_{t-s}(z)|
%\le2\tau.$
%Unfortunately, this is not good enough for our present needs;
% we need to do a little better by showing that $\mathscr{Q}(s)\le2
%\tau^2$: Recall
% that we can represent $X_t:= \sum_{j=0}^{N(t)} Y_j$,
% where $Y_0:=0$, $\{Y_j\}_{j=1}^\infty$ is a sequence of i.i.d.\
% random variables and $\{N(t)\}_{t\ge0}$ is an independent
% rate-one Poisson process. Let $\varphi(\xi):=
% \E\exp(i\xi\cdot Y_1)$ denote the characteristic function of
% the increments of the continuous-time random walk $X$. It is an
%elementary
% calculation that
% \begin{equation}\label{eq:CHF}
% \E\e^{i\xi\cdot X_t}=\e^{-t(1-\varphi(\xi))} \mbox{for all
% $\xi\in\R^d$ and $t\ge0$}.
% \end{equation}
% Therefore, we appeal to the Parseval identity and find that
% \begin{equation}
% \mathscr{Q}(s)= (2\pi)^{-d}\int_{[-\pi,\pi]^d}\left| \e^{-(t+
%\tau-s)(1-\varphi(\xi))} -
% \e^{-(t-s)(1-\varphi(\xi))}\right|^2 \,\d\xi\le2\tau^2,
% \end{equation}
% uniformly for all $s\in(0 ,t)$.
Note that $\mathscr{Q}(s) \le [\sum_z | p_{t+\tau-s}(z)-p_{t-s}
| ]^2 \le4\tau^2$ uniformly for $s\in(0,t)$ from~\eqref{eq:p-p}.
This shows that
\[
Q_2^2 \le16k\tau^2 \int
_0^t\sup_{y\in\Z^d} \bigl\{ \E
\bigl( \bigl\llvert \sigma\bigl(u_s(y)\bigr)\bigr\rrvert
^k \bigr) \bigr\}^{2/k} \,\d s.
\]
Because $|\sigma(z)|\le\lip|z|$ for all $z\in\R$,
the already-proved bound \eqref{eq:Upper:LE} tells us that there
exist constants $c,c_k\in(0 ,\infty)$ $[k\ge2]$ such that
%
%e5.9 #&#
\begin{equation}
\label{eq:LE:k} \sup_{y\in\Z^d}\E \bigl( \bigl\llvert \sigma
\bigl(u_s(y)\bigr)\bigr\rrvert ^k \bigr) \le
c_k^k {\mathrm{e}}^{ck^2s}\qquad \mbox{for all integers $k
\ge2$ and $s\ge0$}.
\end{equation}
Therefore,
%
%e5.10 #&#
\begin{equation}
\label{eq:Q_2} Q_2^2 \le8c^{-1}c_k^2{\mathrm{e}}^{2ckt}\tau^2.
\end{equation}

Finally, we apply BDG Lemma~\ref{lem:BDG} to see that
%
%e5.11 #&#
\begin{eqnarray}
\label{eq:Q_3:prec} %
Q_3^2 &\le&4k\sum
_{y\in\Z^d} \int_t^{t+\tau} \bigl[
p_{t+\tau-s}(y-x) \bigr]^2 \bigl\{ \E \bigl( \bigl\llvert \sigma
\bigl(u_s(y)\bigr)\bigr\rrvert ^k \bigr) \bigr
\}^{2/k} \,\d s
\nonumber
\\[-8pt]
\\[-8pt]
\nonumber
&\le&4kc_k^2\sum_{y\in\Z^d} \int
_t^{t+\tau} \bigl[ p_{t+\tau-s}(y-x)
\bigr]^2 {\mathrm{e}}^{2cks} \,\d s,
\end{eqnarray}
owing to \eqref{eq:LE:k}. Because $\sum_{y\in\Z^d}[p_h(y-x)]^2\le1$
for all $h\ge0$, we find that
%
%e5.12 #&#
\begin{equation}
\label{eq:Q_3} Q_3^2 \le4c^{-1}c_k^2{\mathrm{e}}^{2ck(t+\tau)}\tau.
\end{equation}
We combine \eqref{eq:Q_1}, \eqref{eq:Q_2} and \eqref{eq:Q_3} and
find that
for all integers $k\ge2$, there exists
a finite and positive constant $\tilde{a}:=
\tilde{a}(T,k)$ such that for every $\tau\in(0 ,1)$,
%
%e5.13 #&#
\begin{equation}
\label{eq:mod:u} \sup_{x\in\Z^d} \sup_{t\in(0,T)} \E
\bigl(\bigl\llvert u_{t+\tau}(x) - u_t(x)\bigr\rrvert
^k \bigr) \le \tilde{a} {\mathrm{e}}^{\tilde{a}T}\tau^{k/2}.
\end{equation}
The lemma follows from this bound, and an application of
a quantitative form of the Kolmogorov continuity theorem
\cite{RevuzYor}, Theorem~2.1, page 25.
We omit the remaining details, as they are nowadays standard.
\end{pf}

Our next approximation result is
the highlight of this section, and
refines \eqref{eq:local:modulus} by inspecting
more closely the main contribution to the $O(\tau^{(1+o(1))/2})$ error term
in \eqref{eq:local:modulus}. In order to describe the next
approximation result,
we first define for every fixed $t\ge0$ an infinite-dimensional
Brownian motion $B^{(t)}$ as follows:
%
%e5.14 #&#
\begin{equation}
B^{(t)}_\tau(x) := B_{\tau+t}(x)-B_t(x)\qquad
\bigl(x\in\Z^d, \tau\ge0\bigr).
\end{equation}
If we continue to hold $t$ fixed, then it is easy to see that
$\{B^{(t)}_\bullet(x)\}_{x\in\Z^d}$ is a collection of independent
$d$-dimensional
Brownian motions. Furthermore, the entire process $B^{(t)}$ is independent
of the infinite-dimensional random variable $u_t(\bullet)$, since it
is easy
to see from the proof of the first part of Theorem~\ref{th:exist:unique}
that $u_t$ is a measurable function of $\{B_s(y)\}_{s\in[0,t],y\in\Z^d}$,
which is therefore independent of $B^{(t)}$ by the Markov property of $B$.
Now for every fixed $t\ge0$ and
$x\in\Z^d$, consider the solution $u^{(t)}_\bullet(x)$ to the following
(autonomous/noninteracting) It\^o stochastic differential equation:
%
%e5.15 #&#
\begin{equation}
\label{eq:Feller}%
\cases{\displaystyle\frac{\d u^{(t)}_\tau(x)}{\d\tau} = \frac{\d(\tilde{p}_\tau*u_t)(x)}{\d\tau} +
\sigma \bigl( u^{(t)}_\tau(x) \bigr)\frac{\d B^{(t)}_\tau(x)}{\d\tau},\vspace
*{2pt}
\cr
\mbox{subject to\qquad$u^{(t)}_0(x) =
u_t(x)$}.} %
\end{equation}
Note, once again, that $B^{(t)}$ is independent of $u_t(\bullet)$. Moreover,
%
%e5.16 #&#
\begin{equation}
\sup_{\tau>0}\E \bigl(\bigl\llvert (\tilde{p}_\tau*u_t
) (x)\bigr\rrvert ^2 \bigr) \le\sup_{y\in\Z^d}\E
\bigl(\bigl|u_t(y)\bigr|^2 \bigr)<\infty,
\end{equation}
thanks to the already-proved bound \eqref{eq:Upper:LE}
and the Cauchy--Schwarz inequality. Therefore,
\eqref{eq:Feller} is a standard It\^o-type SDE and hence
has a unique strong solution.

%th5.2 #&#
\begin{theorem}[(The local-diffusion property)]\label{th:local:OU}
For every $t\ge0$, the following holds a.s. for all $x\in\Z^d$:
%
%e5.17 #&#
\begin{equation}
\label{eq:local:OU} u_{t+\tau}(x) = u^{(t)}_\tau(x) + O
\bigl( \tau^{(3/2)+o(1)} \bigr)\qquad \mbox{as $\tau\downarrow0$}.
\end{equation}
\end{theorem}

The proof of Theorem~\ref{th:local:OU} hinges on three technical lemmas
that we state next.

%le5.3 #&#
\begin{lemma}\label{lem:local:OU}
Choose and fix $t\ge0$, $\tau\in[0 ,1]$, and $x\in\Z^d$, and define
%
%e5.18 #&#
\begin{eqnarray}
\mathscr{A} &:=&\sum_{y\in\Z^d} \int
_t^{t+\tau} p_{t+\tau-s}(y-x) \sigma
\bigl(u_s(y)\bigr) \,\d B_s(y),
\nonumber
\\[-8pt]
\\[-8pt]
\nonumber
\mathscr{B} &:=& \int_t^{t+\tau} \sigma
\bigl(u_s(x)\bigr) \,\d B_s(x). %
\end{eqnarray}
Then, for all real numbers $k\ge2$ there exist a finite
constant $C_k>0$---depending on $k$ but
not on $(t ,\tau,x)$---and a finite constant $C>0$---not depending
on $(t ,\tau,x ,k)$---such that
%
%e5.19 #&#
\begin{equation}
\E \bigl( \llvert \mathscr{A}-\mathscr{B}\rrvert ^k \bigr) \le
C_k{\mathrm{e}}^{Ck^2(t+1)}\tau^{3k/2}.
\end{equation}
\end{lemma}

%le5.4 #&#
\begin{lemma}\label{lem:local:OU1}
For every $k\ge2$ and $T\ge1$, there exists a finite constant $C(k ,T)$
such that for every $\tau\in(0 ,1]$,
%
%e5.20 #&#
\begin{equation}
\sup_{t\in[0,T]}\sup_{x\in\Z^d} \E \bigl( \bigl
\llvert u_{t+\tau}(x) - u^{(t)}_\tau(x) \bigr\rrvert
^k \bigr) \le C(k ,T)\tau^{3k/2}.
\end{equation}
\end{lemma}

%le5.5 #&#
\begin{lemma}\label{lem:modulus:OU}
There exists a version of $u^{(\bullet)}$ that is a.s. continuous in
$(t ,\tau)$.
Moreover,
for every $T\ge1$, $\varepsilon\in(0 ,1)$ and $k\ge2$,
%
%e5.21 #&#
\begin{equation}
\sup_{t\in[0,T]}\sup_{x\in\Z^d}\sup
_I \E \biggl( \mathop{\sup_{\nu,\mu\in I}}_{
\nu\neq\mu}
\biggl[ \frac{|u^{(t)}_\nu(x) - u^{(t)}_\mu(x)|}{
|\nu-\mu|^{(1-\varepsilon)/2}} \biggr]^k \biggr)<\infty,
\end{equation}
where ``$\sup_I$'' denotes the supremum over all closed subintervals
$I$ of
$[0 ,T]$ that have length $\le1$.
\end{lemma}

In order to maintain the flow of the discussion, we prove Theorem~\ref{th:local:OU} first. Then we conclude this section by establishing
the three supporting lemmas mentioned above.

\begin{pf*}{Proof of Theorem~\ref{th:local:OU}}
Throughout the proof we choose and fix some $t\in[0 ,T]$ and
$x\in\Z^d$.

Our plan is to prove that for all $\delta\in(0 ,\nicefrac12)$,
%
%e5.22 #&#
\begin{equation}
u_{t+\tau}(x) - u^{(t)}_\tau(x) = O \bigl(
\tau^{(3/2)-\delta} \bigr)\qquad \mbox{as $\tau \downarrow0$, a.s.}
\end{equation}
Henceforth, we choose and fix some $\delta\in(0 ,\nicefrac12)$,
and denote by $A_k,A_k',A_k''$, etc. finite constants
that depend only on a parameter $k\ge2$ that will be selected later, during
the course of the proof.

Thanks to Lemma~\ref{lem:local:OU1}, for all $k\ge2$ and $\tau\in[0 ,1]$,
%
%e5.23 #&#
\begin{equation}
\P \bigl\{\bigl\llvert u_{t+\tau}(x) - u^{(t)}_\tau(x)
\bigr\rrvert \ge\tfrac{1}3\tau^{(3/2)-\delta} \bigr\} \le C(k ,T)\tau
^{\delta k}.
\end{equation}
We can choose $k$ large enough
and then apply the Borel--Cantelli lemma in order to deduce that with
probability one,
%
%e5.24 #&#
\begin{equation}\qquad
\bigl\llvert u_{t+\tau_n}(x) - u^{(t)}_{\tau_n}(x) \bigr
\rrvert < \tau_n^{(3/2)-\delta}\qquad \mbox{for all but a finite number of
$n$'s,}
\end{equation}
where $\tau_n := n^{\delta-(1/2)}$.
Because $\tau_n-\tau_{n+1} \sim\operatorname{const} \times
n^{-1}\tau_n$
as $n\to\infty$, H\"older continuity ensures the following
(Lemmas \ref{lem:modulus} and
\ref{lem:modulus:OU}): Uniformly for all $\tau\in[\tau_{n+1},\tau_n]$,
%
%e5.25 #&#
\begin{eqnarray}
\bigl\llvert u_{t+\tau}(x) - u_{t+\tau_n}(x)\bigr\rrvert +
\bigl\llvert u^{(t)}_{\tau_n}(x) - u^{(t)}_\tau(x)
\bigr\rrvert &=& O \bigl( [\tau_n/n ]^{(1/2)-\delta} \bigr)\qquad
\mbox{a.s.}
\nonumber
\\
&=&O \bigl(\tau_n^{(3/2)-\delta} \bigr),
\end{eqnarray}
by the particular choice of the sequence $\{\tau_n\}_{n=1}^\infty$.
The preceding two displays can now be combined to imply \eqref{eq:local:OU}.
\end{pf*}

\begin{pf*}{Proof of Lemma~\ref{lem:local:OU}}
We may rewrite $\mathscr{B}$ as follows:
%
%e5.26 #&#
\begin{equation}
\mathscr{B} =\sum_{y\in\Z^d} \int_t^{t+\tau}
\1_{\{0\}}(y-x) \sigma\bigl(u_s(y)\bigr) \,\d
B_s(y).
\end{equation}
Therefore, BDG Lemma~\ref{lem:BDG} can be used to show that
%
%e5.27 #&#
\begin{eqnarray}
\nonumber
&& \bigl\{ \E \bigl(\llvert \mathscr{A}-\mathscr{B}\rrvert
^k \bigr) \bigr\}^{2/k}
\\
&&\qquad\le4k\sum_{y\in\Z^d}\int_t^{t+\tau}
\bigl[ p_{t+\tau-s}(y-x) - \1_{\{0\}}(y-x) \bigr]^2
\bigl\{ \E \bigl( \bigl\llvert \sigma\bigl(u_s(y)\bigr)\bigr\rrvert
^k \bigr) \bigr\}^{2/k} \,\d s
\nonumber
\\[-8pt]
\\[-8pt]
\nonumber
&&\qquad\le4kc_k^2{\mathrm{e}}^{2ck(t+1)} \biggl(\sum
_{y\in\Z^d\setminus\{0\}
}\int_0^\tau \bigl[
p_{s}(y) \bigr]^2 \,\d s+\int_{0}^{\tau}
\bigl[ 1- p_s(0) \bigr]^2\,\mathrm{d}s \biggr)
\\
\nonumber
&&\qquad\le4kc_k^2{\mathrm{e}}^{2ck(t+1)} 2\int
_{0}^{\tau} \bigl[ 1-p_s(0)
\bigr]^2\,\mathrm{d}s,
\end{eqnarray}
where $c,c_k$ appear in \eqref{eq:LE:k}. Observe that
$p_s(0)= \P\{X_s=0\}\geq\P\{N_s=0\}={\mathrm{e}}^{-s}$, where $\{N_s\}
_{s\ge0}$
denotes the underlying Poisson clock.
Therefore, we obtain
$\int_0^{\tau}  [ 1-p_s(0) ]^2 \,\d s
\le(1/3)\tau^3$, and hence
%
%e5.28 #&#
\begin{equation}
\E \bigl(\llvert \mathscr{A}-\mathscr{B}\rrvert ^k \bigr)
\le(8/3)^{k/2} k^{k/2}c_k^k{\mathrm{e}}^{ck^2(t+1)}\tau^{3k/2}.
\end{equation}
This implies the lemma.
\end{pf*}

\begin{pf*}{Proof of Lemma~\ref{lem:local:OU1}}
In accord with \eqref{mild}, we may write $u_{t+\tau}(x)$ as
%
%e5.29 #&#
\begin{equation}
( \tilde{p}_{t+\tau}*u_0 ) (x) + \sum
_{y\in\Z^d} \int_0^t
p_{t+\tau-s}(y-x)\sigma\bigl(u_s(y)\bigr) \,\d
B_s(y) + \mathscr{A}, \label{eq:post:t:u}
\end{equation}
where $\mathscr{A}$ was defined in Lemma~\ref{lem:local:OU}.

By the Chapman--Kolmogorov
property of the transition functions $\{p_t\}_{t\ge0}$,
%
%e5.30 #&#
\begin{eqnarray}
&&(\tilde{p}_\tau*u_t) (x)
\nonumber
\\[-8pt]
\\[-8pt]
\nonumber
&&\qquad= (\tilde{p}_{t+\tau}*u_0)
(x) + \sum_{y\in\Z^d}\int_0^t
p_{t+\tau-s}(y-x)\sigma \bigl(u_s(y) \bigr) \,\d
B_s(y).
\end{eqnarray}
The exchange of summation with stochastic integration can be justified,
using the already-proved moment bound \eqref{eq:Upper:LE}
of Theorem~\ref{th:exist:unique}; we omit the details.
Instead, let us apply this in \eqref{eq:post:t:u} to see that
%
%e5.31 #&#
\begin{eqnarray}
u_{t+\tau}(x) &=& (\tilde{p}_\tau*u_t)
(x) + \int_t^{t+\tau} \sigma\bigl(u_s(x)
\bigr) \,\d B_s(x) + ( \mathscr{A}-\mathscr{B} )
\nonumber
\\[-8pt]
\\[-8pt]
\nonumber
&=& (\tilde{p}_\tau*u_t) (x) + \int_0^\tau
\sigma\bigl(u_{t+s}(x)\bigr) \,\d_s B^{(t)}_s(x)
+ ( \mathscr{A}-\mathscr{B} ). %
\end{eqnarray}
Lemma~\ref{lem:local:OU} implies that for all $k\ge2$, $t,\tau\ge0$
and $x\in\Z^d$,
%
%e5.32 #&#
\begin{eqnarray}
\label{eq:local:approx1} %
&&\E \biggl( \biggl| u_{t+\tau}(x) - (
\tilde{p}_\tau*u_t) (x) -\int_0^\tau
\sigma\bigl(u_{t+s}(x)\bigr) \,\d_s B^{(t)}_s(x)
\biggr|^k \biggr)
\nonumber
\\[-8pt]
\\[-8pt]
\nonumber
&&\qquad \le a_k {\mathrm{e}}^{ak^2(t+1)}\tau^{3k/2}, %
\end{eqnarray}
where $a\in(0 ,\infty)$ is universal and $a_k\in(0 ,\infty)$
depends only on $k$.
On the other hand,
%
%e5.33 #&#
\begin{equation}
u^{(t)}_\tau(x) - (\tilde{p}_\tau*u_t
) (x) - \int_0^\tau\sigma
\bigl(u^{(t)}_s(x) \bigr) \,\d B^{(t)}_s(x)=0\qquad
\mbox{a.s.,}
\end{equation}
by the very definition of $u^{(t)}$, and thanks to the fact that
$u^{(t)}_0(y)=u_t(y)$.
The preceding two displays and Minkowski's inequality that
%
%e5.34 #&#
\begin{equation}
\psi(\tau) := \bigl\{\E \bigl( \bigl\llvert u_{t+\tau}(x) -
u^{(t)}_\tau(x) \bigr\rrvert ^k \bigr) \bigr
\}^{1/k} \le a_k^{1/k}{\mathrm{e}}^{ak(t+1)}
\tau^{3/2} + Q,
\end{equation}
where
%
%e5.35 #&#
\begin{equation}
Q := \biggl\{\E \biggl(\biggl\llvert \int_0^\tau
\bigl[\sigma\bigl(u_{t+s}(x)\bigr) - \sigma \bigl( u^{(t)}_s(x)
\bigr) \bigr] \,\d_s B^{(t)}_s(x)\biggr\rrvert
^k \biggr) \biggr\}^{1/k}.
\end{equation}
According to BDG Lemma~\ref{lem:BDG} (actually we need a one-dimensional
version of that lemma only), and since $|\sigma(r)-\sigma(\rho)|
\le\lip|r-\rho|$,
%
%e5.36 #&#
\begin{eqnarray}
Q^2 &\le&4k\lip^2\int_0^\tau
\bigl\{ \E \bigl(\bigl\llvert u_{t+s}(x)-u^{(t)}_s(x)
\bigr\rrvert ^k \bigr) \bigr\}^{2/k} \,\d s
\nonumber
\\[-8pt]
\\[-8pt]
\nonumber
&=&4k\lip^2 \int_0^\tau\bigl[
\psi(s)\bigr]^2 \,\d s. %
\end{eqnarray}
Thus we find that
%
%e5.37 #&#
\begin{eqnarray}
\bigl[\psi(\tau)\bigr]^2 \le2a_k^{2/k}{\mathrm{e}}^{2ak(t+1)}\tau^{3} + 8k\lip^2\int
_0^\tau\bigl[\psi(s)\bigr]^2 \,\d s
\nonumber
\\[-8pt]
\\[-8pt]
\eqntext{\mbox{for all $0\le\tau\le1$}.}
\end{eqnarray}
The lemma follows from this and an application of Gronwall's lemma.
\end{pf*}

\begin{pf*}{Proof of Lemma~\ref{lem:modulus:OU}}
One can model closely a proof after that of Lemma~\ref{lem:modulus}.
However, we omit the details, since this is a result about finite-dimensional
diffusions and as such simpler than Lemma~\ref{lem:modulus}.
\end{pf*}

We conclude this section with a final approximation lemma.
The next assertion shows that the solution to \eqref{SHE} depends
continuously on
its initial function (in a suitable topology).

%le5.6 #&#
\begin{lemma}\label{lem:flow}
Let $u$ and $v$ denote the unique solutions to
\mbox{\eqref{SHE}}, corresponding, respectively, to initial functions $u_0$
and $v_0$. Then
%
%e5.38 #&#
\begin{equation}\qquad
\sup_{x\in\Z^d}\E \bigl(\bigl\llvert u_t(x)-v_t(x)
\bigr\rrvert ^2 \bigr) \le \|u_0-v_0
\|_{\ell^\infty(\Z^d)}^2 {\mathrm{e}}^{\lip^2 t} \qquad\mbox{for all $t\ge0$}.
\end{equation}
\end{lemma}

\begin{pf}
Choose and fix $t\ge0$.
The fact that $\sum_{y\in\Z^d} p_t(y)=1$ alone ensures that
%
%e5.39 #&#
\begin{equation}
\sup_{x\in\Z^d}\bigl\llvert ( \tilde{p}_t*u_0
) (x) - ( \tilde{p}_t*v_0 ) (x) \bigr\rrvert \le
\|u_0-v_0\|_{\ell^\infty(\Z^d)}.
\end{equation}
Therefore, \eqref{mild} and It\^o's isometry together imply that
%
%e5.40 #&#
\begin{eqnarray}
&&\E \bigl(\bigl\llvert u_t(x) - v_t(x)\bigr\rrvert
^2 \bigr)
\nonumber
\\
&&\qquad\le\|u_0-v_0\|_{\ell^\infty(\Z^d)}^2
\\
&&\qquad\quad{}+ \lip^2 \int_0^t \llVert
p_s \rrVert _{\ell^2(\Z^d)}^2 \,\cdot\,\sup
_{y\in\Z^d}\E \bigl( \bigl\llvert u_s(y) -
v_s(y)\bigr\rrvert ^2 \bigr) \,\d s.\nonumber
\end{eqnarray}
Since $\|p_s\|_{\ell^2(\Z^d)}^2 = \P\{X_s=X_s'\}\le1$, where $X'$ is
an independent
copy of $X$, we may conclude that $f(t):=\sup_{x\in\Z^d}\E
(|u_t(x)-v_t(x)|^2)$
satisfies
%
%e5.41 #&#
\begin{equation}
f(t) \le\|u_0-v_0\|_{\ell^\infty(\Z^d)}^2 +
\lip^2\int_0^tf(s) \,\d s.
\end{equation}
Therefore, the lemma
follows from Gronwall's inequality.
\end{pf}

%s6 #&#
\section{Proof of Theorem \texorpdfstring{\protect\ref{th:local:RadonNikodym}}{2.2}}
\label{sec:thm1.2}
Theorem~\ref{th:local:RadonNikodym} is a consequence
of the following result.

%pr6.1 #&#
\begin{proposition}\label{pr:local}
For every $t\ge0$, the following holds a.s. for all $x\in\Z^d$:
%
%e6.1 #&#
\begin{eqnarray}
\label{eq:local:OU1} u_{t+\tau}(x) - u_t(x) = \sigma \bigl(
u_t(x) \bigr) \bigl\{ B_{t+\tau}(x)-B_t(x) \bigr
\} + o \bigl( \tau^{1+o(1)} \bigr)
\nonumber
\\[-8pt]
\\[-8pt]
 \eqntext{\mbox{as $\tau\downarrow0$}.}
\end{eqnarray}
\end{proposition}

Indeed,\vspace*{1pt} we obtain \eqref{eq:local:RN} from this proposition,
simply because well-known properties of Brownian motion
imply that for all $\varepsilon\in(0 ,\nicefrac12)$ and $t\ge0$,
%
%e6.2 #&#
\begin{equation}
\lim_{\tau\downarrow0}\frac{\tau^{1-\varepsilon}}{
B_{t+\tau}(x)-B_t(x)}=0 \qquad\mbox{in probability}.
\end{equation}
Moreover, \eqref{eq:local:LIL} follows from the local law of the iterated
logarithm for Brownian motion. It remains to prove Proposition~\ref{pr:local}.

\begin{pf}
According to \eqref{eq:local:approx1}, for every integer $k\ge2$,
and all $t,\tau\ge0$ and $x\in\Z^d$,
%
%e6.3 #&#
\begin{eqnarray}
&&\E \biggl(\biggl\llvert u_{t+\tau}(x) - u_t(x) -
\int_0^\tau\sigma\bigl(u_{t+s}(x)
\bigr) \,\d_s B^{(t)}_s(x)\biggr\rrvert
^k \biggr)
\nonumber
\\[-8pt]
\\[-8pt]
\nonumber
&&\qquad\le2^{k-1} \bigl[ a_k{\mathrm{e}}^{ak^2(t+1)}
\tau^{3k/2} + \E \bigl(\bigl\llvert u_t(x) - (
\tilde{p}_\tau*u_t ) (x) \bigr\rrvert ^k
\bigr) \bigr]. %
\end{eqnarray}
We may write
%
%e6.4 #&#
%e6.5 #&#
\begin{eqnarray}
&&\E \bigl(\bigl\llvert u_t(x) - (\tilde{p}_\tau*u_t
) (x) \bigr\rrvert ^k \bigr)\nonumber
\\
&&\qquad=\E \biggl( \biggl| u_t(x) - \sum_{y\in\Z^d}p_\tau(y-x)u_t(y)
\biggr|^k \biggr)
\\
&&\qquad=\E \biggl( \biggl| u_t(x)\P\{X_\tau\neq0\} - \sum
_{y\in\Z^d\setminus\{x\}} p_\tau(y-x) u_t(y)
\biggr|^k \biggr).\nonumber
\end{eqnarray}
Because $\P\{X_\tau\neq0\}=1-\exp(-\tau)\le\tau$,
Minkowski's inequality shows that
%
%e6.6 #&#
\begin{eqnarray}
&& \bigl\{ \E \bigl(\bigl\llvert u_t(x) - (
\tilde{p}_\tau*u_t ) (x) \bigr\rrvert ^k
\bigr) \bigr\}^{1/k}
\nonumber
\\
&&\qquad\le\tau \bigl\{ \E \bigl(\bigl|u_t(x)\bigr|^k \bigr) \bigr\}
^{1/k} + \sum_{y\in\Z^d\setminus\{x\}} p_\tau(y-x)
\bigl\{ \E \bigl(\bigl|u_t(y)\bigr|^k \bigr) \bigr\}^{1/k}
\\
&&\qquad\le2\tau\sup_{y\in\Z^d} \bigl\{ \E \bigl(\bigl|u_t(y)\bigr|^k
\bigr) \bigr\}^{1/k}. \nonumber
\end{eqnarray}
We can conclude from this development
and from Theorem~\ref{th:exist:unique} that there exists $A_k<\infty
$, depending
only on $k$, and a universal $A<\infty$ such that
%
%e6.7 #&#
\begin{eqnarray}
\label{eq:uuu} %
&&\E \biggl(\biggl\llvert u_{t+\tau}(x) -
u_t(x) -\int_0^\tau\sigma
\bigl(u_{t+s}(x)\bigr) \,\d_s B^{(t)}_s(x)
\biggr\rrvert ^k \biggr)
\nonumber
\\[-8pt]
\\[-8pt]
\nonumber
&&\qquad \le A_k {\mathrm{e}}^{Ak^2(t+1)} \bigl[ \tau^{3k/2} +
\tau^k \bigr] \le A_k {\mathrm{e}}^{Ak^2(t+1)}
\tau^k,
\end{eqnarray}
for all $\tau\in[0 ,1]$. Now, we may apply BDG Lemma~\ref{lem:BDG} in order to see that
%\begin{equation}\begin{split}
%
%e6.8 #&#
\begin{eqnarray}
&&\biggl[\E \biggl(\biggl |\int_0^\tau\sigma
\bigl(u_{t+s}(x)\bigr) \,\d_s B^{(t)}_s(x)
-\sigma\bigl(u_t(x)\bigr) \bigl\{ B_{t+\tau}(x)-B_t(x)
\bigr\}\biggr |^k \biggr) \biggr]^{2/k}
\nonumber
\\
&&\qquad= \biggl[\E \biggl(\biggl\llvert \int_0^\tau
\bigl[ \sigma\bigl(u_{t+s}(x)\bigr) -\sigma\bigl(u_t(x)
\bigr) \bigr] \,\d_s B^{(t)}_s(x) \biggr\rrvert
^k \biggr) \biggr]^{2/k}
\nonumber
\\[-8pt]
\\[-8pt]
\nonumber
&&\qquad\le4k\lip^2\int_0^\tau \bigl[ \E
\bigl(\bigl\llvert u_{t+s}(x)- u_t(x)\bigr\rrvert
^k \bigr) \bigr]^{2/k} \,\d s
\\
&&\qquad\le\tilde{a}_k{\mathrm{e}}^{\tilde{a}t}\int_0^\tau
s \,\d s \le\operatorname{const}\,\cdot\,\tau^2,\nonumber
\end{eqnarray}
%
%\end{split}\end{equation}
using \eqref{eq:mod:u}.
%\begin{equation}\begin{split}
% &\sup_{x\in\Z^d}\left\{\E\left(\left| \int_0^\tau\sigma(u_{t+s}(x))
% \,\d_s B^{(t)}_s(x)
% -\sigma(u_t(x)) B^{(t)}_\tau(x)\right|^k\right)\right\}^{2/k}\\
% &\hskip2.5in\le\tilde{a}_k\e^{\tilde{a}t}\int_0^\tau s \,\d s
% \le\mbox{const}\cdot\tau^2.
%\end{split}\end{equation}
Therefore, we can deduce from \eqref{eq:uuu} that
%
%e6.9 #&#
\begin{equation}
\E \bigl(\bigl\llvert D(\tau, x) \bigr\rrvert ^k \bigr) \le
c_{k,t} \tau^k\qquad (0\le\tau\le1),
\end{equation}
where we defined
%
%e6.10 #&#
\begin{equation}
D(\tau, x):= u_{t+\tau}(x) - u_t(x) -\sigma
\bigl(u_t(x)\bigr) \bigl\{ B_{t+\tau
}(x)-B_t(x)
\bigr\},
\end{equation}
and $c_{k,t}$ is a finite constant that depends only on $k$ and $t$;
in particular, $c_{k,t}$ does not depend on $\tau$.
Choose and fix some $\eta>\xi>0$ such that
$\eta+\xi<\nicefrac12$, and
then apply the Chebyshev inequality, and the preceding with
any choice of integer $k>\xi^{-1}$, in order to see that
$\sum_{n=1}^\infty\P\{ | D(n^{-\eta}, x) | > n^{-(\eta-\xi)} \}
\le c_{k,t}\sum_{n=1}^\infty n^{-\xi k}<\infty$.
Thus
%
%e6.11 #&#
\begin{equation}
\label{eq:D} D \bigl( n^{-\eta}, x \bigr)= O \bigl( n^{-(\eta-\xi)}
\bigr) \qquad\mbox{a.s.}
\end{equation}
by the Borel--Cantelli lemma. Because
$n^{-\eta} - (n+1)^{-\eta} = O ( n^{-1-\eta} )$,
the modulus of continuity
of Brownian motion, together with Lemma~\ref{lem:modulus}, imply that
%
%e6.12 #&#
\begin{equation}\qquad
\sup_{(n+1)^{-\eta}\le\tau\le n^{-\eta}} \bigl\llvert D \bigl( n^{-\eta}, x \bigr)-
D(\tau, x) \bigr\rrvert = O \bigl( n^{-1/2} \bigr)=o \bigl(
n^{-(\eta+\xi)} \bigr) \qquad\mbox{a.s.}
\end{equation}
Therefore a standard monotonicity argument and \eqref{eq:D} together
reveal that $D(t, x) = O ( t^{(\eta-\xi)/\eta} )$
as $t\downarrow0$, a.s.
Since $\eta>\xi$ are arbitrary positive numbers, it follows that
$\limsup_{t\downarrow0}(\log D(t, x)/\log t )\le1$ a.s. This is another
way to state the result.
\end{pf}

%s7 #&#
\section{Proof of Theorem \texorpdfstring{\protect\ref{th:CLT}}{2.5}} \label{sec:thm1.4}
First we prove a preliminary lemma that guarantees strict
positivity of the solution to the \eqref{SHE}. We follow the method
described in Conus, Joseph and Khoshnevisan \cite{CJK11}, Theorem~5.1,
which in turn borrowed heavily from ideas of Mueller \cite{Mueller}
and Mueller and Nualart \cite{MN}.

%le7.1 #&#
\begin{lemma}\label{lem:strict:pos}
$\inf_{0\le t \le T} u_t(x) > 0$ a.s.
for every $T\in(0 ,\infty)$
and all $x\in\Z^d$ that satisfy $u_0(x) > 0$.
\end{lemma}

\begin{pf}
We are going to prove that if $u_0(x_0)>0$ for a fixed $x_0\in\Z^d$, then
there exist finite and positive constants $A$ and $C$ such that
%
%e7.1 #&#
\begin{equation}
\P \Bigl\{\inf_{0<s<t}u_s(x_0) \le
\varepsilon \Bigr\} \le A\varepsilon ^{C\log\vert
\log\varepsilon\vert},
\end{equation}
for that same point $x_0$,
uniformly for all $\varepsilon\in(0 ,1)$. It turns out to be
convenient to prove
the following equivalent formulation of the preceding:
%
%e7.2 #&#
\begin{equation}
\label{eq:goal:pos} \P \Bigl\{\inf_{0<s<t}u_s(x_0)
\le{\mathrm{e}}^{-n} \Bigr\} \le An^{-Cn},
\end{equation}
simultaneously for all $n\ge1$, after a possible relabeling of the
constants $A,C\in(0 ,\infty)$. If so, then we can
simply let $n\to\infty$ and deduce the lemma.

Without loss of generality we assume that $u_0(0)>0$,
and we aim to prove \eqref{eq:goal:pos} with $x_0=0$.
In fact, we will simplify the exposition further and establish \eqref
{eq:goal:pos}
when $u_0(0)=1$; the general case follows from this one
and scaling. Finally, we appeal to a comparison principle
(Theorem~\ref{th:comparison}) in order to reduce our problem further to the
following special
case:
%
%e7.3 #&#
\begin{equation}
u_0(x) = \delta_0(x) \qquad\mbox{for all $x\in
\Z^d$.}
\end{equation}
Thus we consider this case only from now on.

Let $\mathcal{F}_t:=\sigma\{B_s(x)\dvtx x \in\Z^d, 0<s \le t\}$ describe
the filtration generated by time $t$ by all the Brownian motions,
enlarged so that
$t\mapsto u_t$ is a $C(\R)$-valued (strong) Markov process.
Set $T_0 := 0$, and define iteratively for $k\ge0$ the sequence
of $\{ \mathcal{F}_t \}_{t>0}$-stopping times
%
%e7.4 #&#
\begin{equation}
T_{k+1}:= \inf \bigl\{s>T_k\dvtx u_s(0) \le{\mathrm{e}}^{-k-1} \bigr\},
\end{equation}
using the usual convention that $\inf\varnothing:= \infty$.
We may observe that the preceding definitions imply that,
almost surely on $\{T_k<\infty\}$,
%
%e7.5 #&#
\begin{equation}
u_{T_k}(x) \ge{\mathrm{e}}^{-k}\delta_0(x)\qquad \mbox{for all } x\in\Z^d. \label{eq:restart:sol}
\end{equation}
We plan to apply the strong Markov property. In order to do that,
we first define $u^{(k+1)}$ to be the unique continuous solution to
the \eqref{SHE}
(for same Brownian motions, pathwise),
with initial data $u^{(k+1)}_0(x):={\mathrm{e}}^{-k}\delta_0(x)$. Next we
note that,
for every $k\ge0$, the random field
%
%e7.6 #&#
\begin{equation}
w^{(k+1)}_t(x) := {\mathrm{e}}^ku^{(k+1)}_t(x)
\end{equation}
solves the system
%
%e7.7 #&#
\begin{equation}
\cases{ \displaystyle\frac{\d w^{(k+1)}_t(x)}{\d t}=\bigl(\mathscr{L}w^{(k+1)}_t
\bigr) (x) + \sigma_k \bigl( w^{(k+1)}_t(x)
\bigr)\frac{\d B_t(x)}{\d t}, \vspace*{2pt}
\cr
w_0^{(k+1)}(x)=\delta_0(x),}
\label{eq:auxil}
\end{equation}
where $\sigma_k(y):= {\mathrm{e}}^k\sigma({\mathrm{e}}^{-k}y)$.
Because $\sigma(0)=0$, we have
$\operatorname{Lip}_{\sigma_k} = \lip$, uniformly for all $k\ge1$.
Thus we
can keep track of the
constants in the proof of Lemma~\ref{lem:modulus}, in order to deduce
the existence of a finite constant $K:=K(\varepsilon)$ so that for all
$t,s$ with $|t-s|<1$,
%
%e7.8 #&#
\begin{equation}
\E \biggl(\sup_{0<|t-s|<1}\frac{\llvert w^{(k+1)}_{t}(0)
-w^{(k+1)}_s(0)\rrvert ^m}{|t-s|^{m(1-\varepsilon)/2}} \biggr) \le
Km^2{\mathrm{e}}^{Km^2}, \label{eq:sup:const}
\end{equation}
for all real numbers $m\ge2$.

For each $k\ge0$ let us define
%
%e7.9 #&#
\begin{equation}
T^{(k+1)}_1 = \inf \bigl\{ t>0\dvtx w_t^{(k+1)}(0)
\le{\mathrm{e}}^{-1} \bigr\}.
\end{equation}
Equation \eqref{eq:restart:sol}, the strong Markov property and the comparison
principle (Theorem~\ref{th:comparison}) together
imply that outside of a null set, the solution to the revised SPDE
\eqref{eq:auxil} satisfies
%
%e7.10 #&#
\begin{equation}
{\mathrm{e}}^{-k}w^{(k+1)}_t(x) \le
u_{T_k+t}(x).
\end{equation}
Therefore, in particular,
%
%e7.11 #&#
\begin{equation}
T^{(k+1)}_1 \le T_{k+1} - T_k,
\label{eq:time:order}
\end{equation}
and the stopping times $T^{(k+1)}_1$ and $T^{(\ell+1)}_1 $ are
independent if $k \neq\ell$.
For all real numbers $t\in(0 ,1)$ and $m\ge2$,
%
%e7.12 #&#
\begin{eqnarray}
\label{eq:time:bound} \P \bigl\{ T^{(k+1)}_1 \leq t \bigr\} &\le&\P
\Bigl\{ \sup_{0<s< t}\bigl\llvert w_t^{(k+1)}(0)-w_s^{(k+1)}(0)
\bigr\rrvert \ge1- {\mathrm{e}}^{-1} \Bigr\}
\nonumber
\\[-8pt]
\\[-8pt]
\nonumber
&\le& K m^2{\mathrm{e}}^{Km^2} \bigl(1-{\mathrm{e}}^{-1} \bigr)^{-m} t^{(1-\varepsilon)m/2} ,
\end{eqnarray}
where the last inequality follows by Chebyshev's inequality and
\eqref{eq:sup:const} and is valid
for all $0<\varepsilon<1$. Let us emphasize
that the constant of the bound in \eqref{eq:time:bound} does not
depend on
the parameter $k$ which appears in the superscript of the random
variable $T_1^{(k+1)}$.
Now we compute
%
%e7.13 #&#
\begin{eqnarray}
\label{eq:iid:sum} \P \Bigl\{\inf_{0< s \le t} u_s(0)\le{{\mathrm{e}} }^{-n}
\Bigr\} &=& \P\{ T_n \le t\}
\nonumber\\
&=&\P\bigl\{(T_n -T_{n-1})+\cdots+(T_1
- T_0) \le t\bigr\}
\\
\nonumber
&\le&\P \bigl\{T^{(n)}_1+T^{(n-1)}_1+
\cdots+T^{(1)}_1 \le t \bigr\},
\end{eqnarray}
owing to \eqref{eq:time:order}.

The terms $T^{(n)}_1,\ldots,T^{(1)}_1$ that appear in
the ultimate line of \eqref{eq:iid:sum} are independent nonnegative
random variables.
By the triangle inequality, if the sum of those terms
is at most $t$, then certainly
it must be that at least $n/2$ of those terms are at most $2t/n$.
(This application of the triangle inequality is also known as the
pigeon-hole principle.)
If $n$ is an even integer, larger than $t>2$,
then a simple union bound on \eqref{eq:iid:sum} and \eqref{eq:time:bound}
yields
%
%e7.14 #&#
\begin{eqnarray}
\nonumber
&&\P \Bigl\{\inf_{0< s \le t}u_s(0)\le{\mathrm{e}}^{-n} \Bigr\} \label{eq:n-m:rel}
\\
&&\qquad\le\pmatrix{n
\cr
n/2} K^{n/2} m^n{\mathrm{e}}^{Km^2n/2}
\bigl(1-{\mathrm{e}}^{-1} \bigr)^{-mn/2} (2t/n)^{(1-\varepsilon)mn/4}
\\
\nonumber
&&\qquad\le\widetilde{K}^n m^n{\mathrm{e}}^{Km^2n/2} \bigl(1-{\mathrm{e}}^{-1} \bigr)^{-mn/2}
t^{(1-\varepsilon)mn/4} n^{-(1-\varepsilon)mn/4}2^{n(1+m(1-\varepsilon)/4)},
\end{eqnarray}
uniformly for all real numbers $m\ge2$.
Now we set $m := \log n / \log\log n$ in \eqref{eq:n-m:rel}
in order to deduce \eqref{eq:goal:pos} for $x_0=0$ and
every $n\ge1$ sufficiently large. This readily yields
\eqref{eq:goal:pos}.
\end{pf}

Next we show that if we start with an initial profile $u_0$
such that $u_0(x)>0$ \emph{for at least one point} $x\in\Z^d$, then
$u_t(z)>0$
\emph{for all} $z \in\Z^d$ and $t>0$ a.s.
Because we are interested in establishing a lower bound,
we may apply scaling and a comparison theorem (Theorem~\ref{th:comparison})
in order to reduce our problem to the following special case:
%
%e7.15 #&#
\begin{equation}
u_0=\delta_0.
\end{equation}
In this way, we are led to the following representation of the solution:
%
%e7.16 #&#
\begin{equation}
\label{eq:udelta} u_t(x)=p_t(x)+\int_0^t
\sum_{y\in\Z^d} p_{t-s}(y-x) \sigma
\bigl(u_s(y) \bigr) \,\d B_s(y).
\end{equation}

%pr7.2 #&#
\begin{proposition}\label{prop:pos}
If $u_0=\delta_0$, then
$u_t(x)>0$ for all $x \in\Z^d$ and $t>0$ a.s.
\end{proposition}

Proposition~\ref{prop:pos} follows from a few preparatory lemmas.

%le7.3 #&#
\begin{lemma}\label{lem:pos1}
If $u_0=\delta_0$, then
%
%e7.17 #&#
\begin{equation}\quad
\E \bigl(\bigl|u_t(x)\bigr|^2 \bigr) \le\exp \bigl(
\lip^2 t \bigr)\,\cdot\, \bigl[p_t(x)\bigr]^2\qquad
\mbox{for all $t>0$ and $x\in\Z^d$}.
\end{equation}
\end{lemma}

\begin{pf}
We begin with representation \eqref{eq:udelta} of the solution $u$, in integral
form, and appeal to Picard's iteration in order to prove the lemma.

Let $u^{(0)}_t(x) := 1$ for all $t\ge0, x\in\Z^d$,
and then let $\{u^{(n+1)}\}_{n\ge0}$ be defined iteratively by
%
%e7.18 #&#
\begin{equation}
u_t^{(n+1)}(x) := p_t(x)+ \int
_0^t\sum_{y\in\Z^d}
p_{t-s}(y-x) \sigma \bigl(u_s^{(n)}(y) \bigr) \,\d
B_s(y).
\end{equation}

Let us define
%
%e7.19 #&#
\begin{equation}
M^{(k)}_t := \sup_{x\in\Z^d} \E \biggl(
\biggl\llvert \frac{u^{(n+1)}_t(x)}{p_t(x)}\biggr\rrvert ^2 \biggr),
\end{equation}
and apply It\^o's isometry in order to deduce the recursive inequality
for the $M^{(k)}$'s,
%
%e7.20 #&#
\begin{equation}
\label{eq:recur} M^{(n+1)}_t\le1+ \operatorname{Lip}_{\sigma}^2
\,\cdot\, \sup_{x\in\Z^d}\int_0^t
\sum_y \biggl[\frac{p_{t-s}(y-x)p_s(y)}{p_t(x)}
\biggr]^2 M^{(n)}_s \,\d s.
\end{equation}
Because
$\sum_{y\in\Z^d}[f(y)]^2\le[\sum_{y\in\Z^d}f(y)]^2$
for all $f\dvtx\Z^d\to\R_+$, the semigroup property of $\{p_t\}_{t>0}$
yields the bound
%
%e7.21 #&#
\begin{equation}
\label{eq:3p} \sum_{y\in\Z^d} \bigl[p_{t-s}(y-x)p_s(y)
\bigr]^2 \le\bigl[p_t(x)\bigr]^2,
\end{equation}
whence $M^{(n+1)}_t \le1 + \lip^2\,\cdot\,\int_0^t M^{(n)}_s \,\d s$
for all $t>0$ and $n\ge0$.
It follows readily from this that $M^{(n)}_t \le\exp(\lip^2 t)$,
uniformly for all $n\ge0$ and $t>0$; equivalently,
%
%e7.22 #&#
\begin{equation}
\E \bigl(\bigl|u^{(n)}_t(x)\bigr|^2 \bigr)\le{\mathrm{e}}^{\lip^2 t}\bigl[p_t(x)\bigr]^2,
\end{equation}
uniformly for all $n\ge0$, $x\in\Z^d$ and $t>0$. The lemma follows
from this and Fatou's lemma, since
$u^{(n)}_t(x)\to u_t(x)$ in $L^2(\P)$ as $n\to\infty$.
\end{pf}

Our next lemma shows that the random term on
the right-hand side of \eqref{eq:udelta} is small,
for small time, as compared with the nonrandom term in \eqref{eq:udelta}.

%le7.4 #&#
\begin{lemma}\label{lem:pos2}
Assume the conditions of Proposition~\ref{prop:pos}.
Then there exists a finite constant $C>0$ such that for all $t\in(0 ,1)$,
%
%e7.23 #&#
\begin{equation}
\label{eq:randombd} \sup_{x\in\Z^d}\P \biggl\{\biggl\llvert \int
_0^t \sum_{y\in\Z^d}
p_{t-s}(y-x) \sigma \bigl(u_s(y) \bigr)\,\d
B_s(y) \biggr\rrvert >\frac{p_t(x)}{2} \biggr\} \le Ct.
\end{equation}
\end{lemma}
\begin{pf}
By Lemma~\ref{lem:pos1} and It\^o's isometry,
%
%e7.24 #&#
\begin{eqnarray}
&&\E \biggl(\biggl\llvert \int_0^t \sum
_{y\in\Z^d} p_{t-s}(y-x) \sigma
\bigl(u_s(y) \bigr)\,\d B_s(y)\biggr\rrvert
^2 \biggr)
\nonumber
\\
&&\qquad\le\lip^2\,\cdot\,\int_0^t
\sum_{y\in\Z^d} \bigl[p_{t-s}(y-x)
p_s(y) \bigr]^2{\mathrm{e}}^{\lip^2 s} \,\d s \\
&&\qquad\le
\lip^2 \bigl[p_t(x)\bigr]^2\,\cdot\,\int
_0^t{\mathrm{e}}^{\lip^2 s} \,\d s,\nonumber
\end{eqnarray}
where we have used \eqref{eq:3p} in the last inequality.
Because $\int_0^t \exp(\lip^2 s) \,\d s\le ct $
for all $t\in(0 ,1)$ with $c:=\exp(\lip^2)$,
the lemma follows from Chebyshev's inequality.
\end{pf}

Now we can establish Proposition~\ref{prop:pos}.

\begin{pf*}{Proof of Proposition~\ref{prop:pos}}
Let us choose and fix an arbitrary $x\in\Z^d$.
By the strong Markov property of the solution,
and thanks to Lemma~\ref{lem:strict:pos},
we know that once the solution becomes positive at a point,
it remains positive at that point at all future times, almost surely.
Thus it suffices to show that $u_t(x)>0$ for all times of the form $t=2^{-k}$,
when $k$ is a large enough integer.
This is immediate from \eqref{eq:udelta}
and \eqref{eq:randombd}, thanks to the Borel--Cantelli lemma.
\end{pf*}

The preceding lemmas lay the groundwork for the proof of Theorem~\ref{th:CLT}.
We now proceed with the main proof.

\begin{pf*}{Proof of Theorem~\ref{th:CLT}}
Let us first consider the case $m=1$ and, without loss of generality,
$x_1 = 0$. In this
case, we may write
%
%e7.25 #&#
\begin{eqnarray}
&&\lim_{\tau\downarrow0} \P \bigl\{ S\bigl(u_{t+\tau}(0)
\bigr) - S\bigl(u_t(0)\bigr) \le q\sqrt{\tau} \bigr\}
\nonumber\\
&&\qquad=\lim_{\tau\downarrow0} \P \biggl\{ \int_{u_t(0)}^{u_{t+\tau}(0)}
\frac{\d y}{\sigma(y)}\le q\sqrt{\tau} \biggr\}
\\
\nonumber
&&\qquad= \lim_{\tau\downarrow0} \P \biggl\{ \int_{u_t(0)}^{u_{t+\tau}(0)}
\biggl(\frac{1}{\sigma(y)}-\frac{1}{\sigma(u_t(0))} \biggr)\,\d y +\frac{u_{t+\tau}(0)-u_t(0)}{\sigma(u_t(0))} \le
q\sqrt\tau \biggr\}.
\end{eqnarray}
Lemma~\ref{lem:strict:pos} and the positivity condition on
$\sigma$ ensure that $\sigma(u_t(0))>0$ a.s. Therefore,
the theorem follows from Theorem~\ref{th:local:RadonNikodym}
if we were to show that
%
%e7.26 #&#
\begin{equation}
\label{eq:goal:tau} \qquad\frac{1}{\sqrt\tau}\int_{u_t(0)}^{u_{t+\tau}(0)}
\biggl(\frac{1}{\sigma(y)}-\frac{1}{\sigma(u_t(0))} \biggr)\,\d y \to0 \qquad\mbox{almost surely,
as $\tau\downarrow0$}.
\end{equation}
Let $\mathcal{I}(t,t+\tau)$ denote
the random closed interval with endpoints $u_t(0)$ and $u_{t+\tau}(0)$.
Our strict positivity result (Lemma~\ref{lem:strict:pos}) implies that
%
%e7.27 #&#
\begin{equation}
\label{eq:I} \mathcal{I}(t ,t+\tau)\subset(0 ,\infty) \qquad\mbox{for all $t,\tau>0$
a.s.,}
\end{equation}
and thus paves way for the a.s. bounds
%
%e7.28 #&#
\begin{eqnarray*}
\label{eq:tau:bound} %
\biggl\llvert \int_{u_t(0)}^{u_{t+\tau}(0)}
\biggl(\frac{1}{\sigma(y)}-\frac{1}{\sigma(u_t(0))}
 \biggr) \,\d y\biggr\rrvert &\le&\lip
\,\cdot\,\frac{|u_t(0) - u_{t+\tau}(0)|^2}{
\inf_{y\in\mathcal{I}(t,t+\tau)}|\sigma(y)|^2}
\\
&=&O \bigl(\tau\log\vert\log\tau\vert \bigr) \qquad(\tau\downarrow0);
\end{eqnarray*}
see \eqref{eq:local:LIL} for the last part.
This implies \eqref{eq:goal:tau}
and thus completes our proof for $m=1$. The proof for general $m$ is
an easy adaption since $\{B(x_j)\}_{j=1}^m$ are i.i.d. Brownian motions.
\end{pf*}

%s8 #&#
\section{Preliminaries for the proof of Theorem \texorpdfstring{\protect\ref{th:as}}{2.7}}
\label
{sec:prelim thm1.6}
The following function will play a prominent role in the ensuing analysis:
%
%e8.1 #&#
\begin{equation}
\label{eq:P:bar} \bar{P}(\tau) := \llVert p_\tau\rrVert
_{\ell^2(\Z^d)}^2 =\sum_{x\in\Z^d}
\bigl[p_\tau(x)\bigr]^2\qquad \mbox{for all $\tau\ge0$}.
\end{equation}
Because of the Chapman--Kolmogorov property, we can also think of
$\bar{P}$ as
%
%e8.2 #&#
\begin{equation}
\label{eq:P:bar:replica} \bar{P}(\tau) := \P\bigl\{X_\tau- X_\tau'
=0\bigr\},
\end{equation}
where $X'$ is an independent copy of $X$. There is another useful
way to think of $\bar{P}$ as well. Using the fact that
%
%e8.3 #&#
\begin{equation}
\label{eq:CHF} \E{\mathrm{e}}^{i\xi\,\cdot\, X_t}
={\mathrm{e}}^{-t(1-\varphi(\xi))} \qquad\mbox{for all $
\xi\in\R^d$ and $t\ge0$}
\end{equation}
and the Plancherel theorem, we see that
%
%e8.4 #&#
\begin{eqnarray}
\label{eq:P:bar:FT} %
\bar{P}(\tau) &=& (2\pi)^{-d}\int
_{(-\pi,\pi)^d}\bigl\llvert \E\exp (i\xi\,\cdot\, X_\tau) \bigr
\rrvert ^2 \,\d\xi
\nonumber
\\[-8pt]
\\[-8pt]
\nonumber
&=&(2\pi)^{-d}\int_{(-\pi,\pi)^d}{\mathrm{e}}^{-2\tau(1-\Re\varphi(\xi
))} \,\d
\xi,
\end{eqnarray}
where $\varphi(\xi)=\E[\exp(i\xi\,\cdot\, Z_1)]$; recall that $Z_1$ is the
distribution of jump size.
Therefore, in particular, the Laplace transform of $\bar{P}$ is
%
%e8.5 #&#
\begin{eqnarray}
\label{eq:Upsilon} %
\Upsilon(\beta) &:= &\int_0^\infty{\mathrm{e}}^{-\beta\tau}\bar {P}(\tau)
\,\d\tau\qquad (\beta\ge0)
\nonumber
\\[-8pt]
\\[-8pt]
\nonumber
&=&(2\pi)^{-d}\int_{(-\pi,\pi)^d}\frac{\d\xi}{\beta+2(1-\Re
\varphi(\xi))}.
\end{eqnarray}
The interchange of the integrals is justified by Tonelli's theorem, since
$1-\Re\varphi(\xi)\ge0$.

Note that $\Upsilon(0)$ agrees with \eqref{Upsilon(0)}. Also,
the classical theory of random walks tells us that $X-X'$ is transient
if and only if
$\Upsilon(0)=\int_0^\infty\bar{P}(\tau) \,\d\tau<\infty$, which
is in turn
equivalent to the condition
%
%e8.6 #&#
\begin{equation}
\int_{(-\pi,\pi)^d}\frac{\d\xi}{1-\Re\varphi(\xi)}<\infty;
\end{equation}
this is the Chung--Fuchs theorem \cite{ChungFuchs}, transliterated to
the setting of\break continuous-time symmetric random walks, thanks to a standard
Poissonization argument which we feel free to omit.

%le8.1 #&#
\begin{lemma}\label{lem:L2:preserved}
If $u_0\in\ell^2(\Z^d)$, then
$u_t\in\ell^2(\Z^d)$ a.s. for all $t\ge0$. Moreover, for every
$\beta
\ge0$
such that $\lip^2\Upsilon(\beta)<1$,
%
%e8.7 #&#
\begin{equation}
\label{eq:E(E_t)} \E \bigl(\|u_t\|_{\ell^2(\Z^d)}^2 \bigr)
\le\frac{\|u_0\|_{\ell
^2(\Z^d)}^2
{\mathrm{e}}^{\beta t}}{1-\lip^2\Upsilon(\beta)}\qquad \mbox{for all $t\ge0$}.
\end{equation}
\end{lemma}

\begin{pf}
Let $u^{(0)}_t(x) := u_0(x)$ for all $t\ge0$ and $x\in\Z^d$,
and define $u^{(k)}$ to be the resulting $k$th-step approximation
to $u$ via Picard iteration. It follows that
%
%e8.8 #&#
\begin{eqnarray}
\label{eq:mild:L2} %
&&\E \bigl(\bigl|u_t^{(n+1)}(x)\bigr|^2
\bigr)\nonumber
\\
&&\qquad= \bigl\llvert (\tilde{p}_t*u_0 ) (x)\bigr\rrvert
^2 +\sum_{y\in\Z^d}\int_0^t
\bigl[p_{t-s}(y-x)\bigr]^2\E \bigl(\bigl\llvert \sigma
\bigl(u_s^{(n)}(y) \bigr)\bigr\rrvert ^2 \bigr)
\,\d s
\\
&&\qquad\le\bigl\llvert (\tilde{p}_t*u_0 ) (x)\bigr\rrvert
^2 + \lip^2\sum_{y\in\Z^d}\int
_0^t\bigl[p_{t-s}(y-x)
\bigr]^2\E \bigl( \bigl\llvert u_s^{(n)}(y)
\bigr\rrvert ^2 \bigr) \,\d s.\nonumber
\end{eqnarray}
We may add over all $x\in\Z^d$ to deduce from this and Young's
inequality that
%
%e8.9 #&#
\begin{equation}\quad
\label{eq:L2L2} \E \bigl( \bigl\llVert u^{(n+1)}_t\bigr
\rrVert _{\ell^2(\Z^d)}^2 \bigr) \le\|u_0
\|_{\ell^2(\Z^d)}^2+\lip^2\int_0^t
\bar{P}(t-s)\E \bigl(\bigl\llVert u^{(n)}_s\bigr\rrVert
_{\ell^2(\Z^d)}^2 \bigr) \,\d s.
\end{equation}

Since $\Upsilon(\beta)=\beta^{-1}\int_0^\infty\exp(-s)\bar
{P}(s/\beta
) \,\d s
\le\beta^{-1}<\infty$,
we can find $\beta>0$ large enough to guarantee that
$\lip^2\Upsilon(\beta)<1$.

We multiply both sides of \eqref{eq:L2L2} by $\exp(-\beta t)$---for
this choice
of $\beta$---and notice from \eqref{eq:L2L2} that
%
%e8.10 #&#
\begin{equation}
A_k := \sup_{t\ge0} \bigl[ {\mathrm{e}}^{-\beta t}
\E \bigl(\bigl\llVert u^{(k)}_t\bigr\rrVert
_{\ell^2(\Z^d)}^2 \bigr) \bigr] \qquad(k\ge0)
\end{equation}
satisfies
%
%e8.11 #&#
\begin{equation}
A_{n+1}\le\|u_0\|_{\ell^2(\Z^d)}^2 +
\lip^2\Upsilon(\beta) A_n\qquad \mbox{for all $n\ge0$. }
\end{equation}
Since $A_0=\|u_0\|_{\ell^2(\Z^d)}^2$,
the preceding shows that $\sup_{n\ge0}A_n$ is bounded above by
$(1-\lip^2\Upsilon(\beta))^{-1}\|u_0\|_{\ell^2(\Z^d)}^2$.
\end{pf}

%pr8.2 #&#
\begin{proposition}\label{pr:L2:int}
If $u_0\in\ell^1(\Z^d)$, then
for every $\beta\ge0$ such that\break $\lip^2\Upsilon(\beta)<1$,
%
%e8.12 #&#
\begin{equation}
\int_0^\infty{\mathrm{e}}^{-\beta t}\E \bigl(
\|u_t\|_{\ell^2(\Z
^d)}^2 \bigr)\,\d t\le
\frac{\|u_0\|_{\ell^1(\Z^d)}^2\Upsilon(\beta)}{1-\lip^2\Upsilon
(\beta)}.
\end{equation}
Moreover,
%
%e8.13 #&#
\begin{equation}
\int_0^\infty{\mathrm{e}}^{-\beta t}\E \bigl(
\|u_t\|_{\ell^2(\Z
^d)}^2 \bigr)\,\d t =\infty,
\end{equation}
for all $\beta\ge0$ such that $\ell_\sigma^2\Upsilon(\beta)\ge1$.
\end{proposition}

\begin{pf}
We proceed as we did for Lemma \ref{lem:L2:preserved}. But instead
of deducing \eqref{eq:L2L2} from \eqref{eq:mild:L2},
we use a different bound for $\|\tilde{p}_t*u_0\|_{\ell^2(\Z^d)}$
%
%e8.14 #&#
\begin{eqnarray}
\label{eq:ell^2} %
&&\E \bigl( \bigl\llVert u^{(n+1)}_t
\bigr\rrVert _{\ell^2(\Z^d)}^2 \bigr)\nonumber
\\
&&\qquad \le\|p_t\|_{\ell^2(\Z^d)}^2\|u_0
\|_{\ell^1(\Z^d)}^2 +\lip^2\int_0^t
\bar{P}(t-s)\E \bigl(\bigl\llVert u^{(n)}_s\bigr\rrVert
_{\ell^2(\Z^d)}^2 \bigr) \,\d s
\\
&&\qquad = \bar{P}(t)\|u_0\|_{\ell^1(\Z^d)}^2 +
\lip^2\int_0^t \bar{P}(t-s)\E
\bigl(\bigl\llVert u^{(n)}_s\bigr\rrVert
_{\ell^2(\Z^d)}^2 \bigr) \,\d s, \nonumber
\end{eqnarray}
thanks to a slightly different application of Young's inequality.
If we integrate both sides $[\exp(-\beta t)\,\d t]$, then we find that
%
%e8.15 #&#
\begin{equation}
I_k:=\int_0^\infty {\mathrm{e}}^{-\beta t}\E \bigl( \bigl\llVert u^{(k)}_t\bigr
\rrVert _{\ell^2(\Z
^d)}^2 \bigr) \,\d t\qquad (k\ge0)
\end{equation}
satisfies
%
%e8.16 #&#
\begin{eqnarray}
I_{n+1} &\le&\|u_0\|_{\ell^1(\Z^d)}^2
\int_0^\infty{\mathrm{e}}^{-\beta
t}\bar {P}(t) \,\d
t+ I_n\times\lip^2\int_0^\infty{\mathrm{e}}^{-\beta t}\bar{P}(t) \,\d t
\nonumber
\\[-8pt]
\\[-8pt]
\nonumber
&=& \|u_0\|_{\ell^1(\Z^d)}^2\Upsilon(\beta) +
I_n \lip^2\Upsilon(\beta);
\end{eqnarray}
see \eqref{eq:Upsilon}.
The first portion of the lemma follows from this, induction and
Fatou's lemma since
$\lip^2\Upsilon(\beta)<1$.

Next, let us suppose that $\ell_\sigma^2\Upsilon(\beta)\ge1$. The
following complimentary form of~\eqref{eq:ell^2} holds [for
the same reasons that \eqref{eq:ell^2} held]:
%
%e8.17 #&#
\begin{equation}
\label{eq:LB:L2} \quad\E \bigl(\|u_t\|_{\ell^2(\Z^d)}^2
\bigr) \ge\llVert \tilde{p}_t*u_0\rrVert
_{\ell^2(\Z^d)}^2 +\ell_\sigma^2 \int
_0^t\bar{P}(t-s)\E \bigl(\|u_s
\|_{\ell^2(\Z^d)}^2 \bigr)\,\d s.
\end{equation}
It is not hard to verify directly that
%
%e8.18 #&#
\begin{equation}
\label{eq:LB:p*u} \|\tilde{p}_t*u_0\|_{\ell^2(\Z^d)}^2
\ge u_0^2(x_0)\|p_t
\|_{\ell^2(\Z^d)}^2,
\end{equation}
whence, by $u_0(x_0)>0$ for some $x_0>0$, it follows that
%
%e8.19 #&#
\begin{equation}
F(t) := \E \bigl(\|u_t\|_{\ell^2(\Z^d)}^2 \bigr)\qquad (t
\ge0)
\end{equation}
solves the renewal inequality
%
%e8.20 #&#
\begin{equation}
F(t) \ge u_0^2(x_0)\bar{P}(t) +
\ell_\sigma^2 \int_0^t
\bar{P}(t-s)F(s) \,\d s.
\end{equation}
Therefore, $\tilde{F}(\beta):=\int_0^\infty\exp(-\beta t)F(t) \,\d
t$ satisfies
%
%e8.21 #&#
\begin{equation}
\tilde{F}(\beta) \ge u_0^2(x_0)\Upsilon(
\beta) + \ell_\sigma^2\Upsilon(\beta)\tilde{F}(\beta).
\end{equation}
Since $u_0(x_0)>0$ and $\Upsilon(\beta)>0$
for all $\beta\ge0$, it follows that $\tilde{F}(\beta)=\infty$ whenever
$\ell_\sigma^2\Upsilon(\beta)\ge1$.
\end{pf}

%pr8.3 #&#
\begin{proposition}\label{pr:bounded:preserved}
If $u_0\in\ell^1(\Z^d)$, then
%
%e8.22 #&#
\begin{equation}
\label{eq:bounded:preserved} \sup_{t\ge0}\sup_{x\in\Z^d}u_t(x)<
\infty,\qquad \sum_{y\in\Z^d}\int_0^\infty
\bigl\llvert \sigma\bigl(u_s(y)\bigr) \bigr\rrvert ^2\,\d
s<\infty\qquad\mbox{a.s.}
\end{equation}
Moreover: \textup{(i)}
If, in addition, $q:=\lip^2\Upsilon(0)<1$, then
%
%e8.23 #&#
\begin{eqnarray}
\E \Bigl(\sup_{t\ge0}\sup_{x\in\Z^d}
\bigl|u_t(x)\bigr|^2 \Bigr) &\le&\E \Bigl(\sup_{t\ge0}
\|u_t\|_{\ell^1(\Z^d)}^2 \Bigr)
\nonumber
\\[-8pt]
\\[-8pt]
\nonumber
&\le&2\|u_0
\|_{\ell^1(\Z^d)}^2+\frac{8q \,\cdot\,\|u_0\|_{\ell^1(\Z
^d)}^2}{1-q}.
\end{eqnarray}
\textup{(ii)} If, in addition, $\ell_\sigma^2\Upsilon(0)\ge1$,
then
%
%e8.24 #&#
\begin{equation}
\E \Bigl(\sup_{t\ge0}\|u_t\|_{\ell^1(\Z^d)}^2
\Bigr) =\int_0^\infty\E \bigl(
\|u_s\|_{\ell^2(\Z^d)}^2 \bigr) \,\d s= \infty.
\end{equation}
\end{proposition}

%re8.4 #&#
\begin{remark}\label{rem:liminf=0}
Clearly, \eqref{eq:bounded:preserved} implies that if
$u_0\in\ell^1(\Z^d)$, then
%
%e8.25 #&#
\begin{equation}
\liminf_{t\to\infty}\sup_{x\in\Z^d}\bigl|\sigma
\bigl(u_t(x)\bigr)\bigr|^2 \le\liminf_{t\to\infty}
\sum_{x\in\Z^d}\bigl|\sigma\bigl(u_t(x)
\bigr)\bigr|^2=0 \qquad\mbox{a.s.}
\end{equation}
If, in addition, $\ell_\sigma>0$ (say), then we can deduce from the
preceding fact that
$\liminf_{t\to\infty}\sup_{x\in\Z^d} |u_t(x)|=0$
a.s.
%\qed
\end{remark}

Recall that $X-X'$ is transient if and only if $\Upsilon(0)<\infty$.
Therefore,
in order for the condition $\lip^2\Upsilon(0)<1$ to hold, it is
necessary---though
not sufficient---that $X-X'$ be transient.

\begin{pf*}{Proof of Proposition~\ref{pr:bounded:preserved}}
First of all,
Theorem~\ref{th:exist:unique} assures us that $u_t(x)\ge0$ a.s.,
and hence $\|u_t\|_{\ell^1(\Z^d)}=\sum_{x\in\Z^d}u_t(x)$.
Therefore, if we add both sides of \eqref{mild}, then we find that
%
%e8.26 #&#
\begin{equation}
\label{eq:ell^1} \|u_t\|_{\ell^1(\Z^d)} = \|u_0
\|_{\ell^1(\Z^d)} + \sum_{y\in\Z^d}\int
_0^t\sigma\bigl(u_s(y)\bigr) \,\d
B_s(y).
\end{equation}
(It is easy to apply the moment bound of Theorem~\ref{th:exist:unique}
to justify the interchange of the sum and the stochastic integral.)
In particular, it follows that
%
%e8.27 #&#
\begin{equation}
\label{eq:M} M_t := \|u_t\|_{\ell^1(\Z^d)} \qquad(t\ge0)
\end{equation}
defines a nonnegative continuous martingale
with mean $\|u_0\|_{\ell^1(\Z^d)}$. Its\break quadratic variation
satisfies the following relations:
%
%e8.28 #&#
\begin{equation}
\langle M\rangle_t = \sum_{y\in\Z^d}\int
_0^t \bigl|\sigma\bigl(u_s(y)
\bigr)\bigr|^2 \,\d s\le \lip^2 \int_0^t
\|u_s\|_{\ell^2(\Z^d)}^2 \,\d s.
\end{equation}
Bound \eqref{eq:Upper:LE} of Theorem~\ref{th:exist:unique}
is more than enough to show that $M:=\{M_t\}_{t\ge0}$
is a continuous $L^2(\P)$ martingale. Since $M_t\ge0$ a.s.
(Theorem~\ref{th:exist:unique}) it follows from the martingale
convergence theorem that $\lim_{t\to\infty}M_t$ exists a.s.
and is finite a.s., which proves the first part of (\ref
{eq:bounded:preserved}).
And therefore, $\langle M\rangle_\infty=
\sum_{y\in\Z^d}\int_0^t|\sigma(u_s(y))|^2 \,\d s$ has to be
also a.s. finite., since we can realize $M_t$ as $W(\langle M\rangle_t)$
for some Brownian motion $W$, thanks to the Dubins, Dambis-Schwartz
representation theorem \cite{RevuzYor}, page 170.

(i) If we know also that $\lip^2\Upsilon(0)<1$,
then Proposition~\ref{pr:L2:int} guarantees that
$\E\langle M\rangle_\infty$ is bounded
from above by $(1-\lip^2\Upsilon(0))^{-1}
\lip^2\Upsilon(0)\times\break \|u_0\|_{\ell^1(\Z^d)}^2<\infty$,
whence it follows that $M:=\{M_t\}_{t\ge0}$ is a continuous $L^2(\P)$-bounded
martingale with
%
%e8.29 #&#
\begin{equation}
\E \Bigl(\sup_{t\ge0}M_t^2 \Bigr) \le2
\|u_0\|_{\ell^1(\Z
^d)}^2+\frac{8
\lip^2\Upsilon(0) \,\cdot\,\|u_0\|_{\ell^1(\Z^d)}^2
}{1-\lip^2\Upsilon(0)},
\end{equation}
thanks to Doob's maximal inequality. This proves part (i)
because $\|u_t\|_{\ell^\infty(\Z^d)}$
is bounded above by $\|u_t\|_{\ell^1(\Z^d)}$.

(ii) Finally consider the case $\ell_\sigma^2\Upsilon(0)\ge1$.
Since
%
%e8.30 #&#
\begin{eqnarray}
\E \bigl(\|u_t\|_{\ell^1(\Z^d)}^2 \bigr)&=&\E
\bigl( M_t^2 \bigr) =\| u_0
\|_{\ell^1(\Z^d)}^2 + \sum_{y\in\Z^d}\int
_0^t\E \bigl(\bigl\llvert \sigma
\bigl(u_s(y)\bigr)\bigr\rrvert ^2 \bigr) \,\d s
\nonumber
\\[-8pt]
\\[-8pt]
\nonumber
&\ge&\| u_0\|_{\ell^1(\Z^d)}^2 +\ell_\sigma^2
\int_0^t\E \bigl(\|u_s
\|_{\ell^2(\Z^d)}^2 \bigr) \,\d s,
\end{eqnarray}
it suffices to show that this final integral is unbounded (as a
function of $t$),
which follows from the second part of Proposition~\ref{pr:L2:int}.
\end{pf*}

%co8.5 #&#
\begin{corollary}\label{co:iff}
If $u_0\in\ell^1(\Z^d)$, then the following
is a $\P$-null set:
%
%e8.31 #&#
\begin{equation}
\Bigl\{\omega\dvtx\lim_{t\to\infty}\sup_{x\in\Z^d}
\bigl|u_t(x) (\omega)\bigr|=0 \Bigr\} \,\triangle\, \Bigl\{ \omega\dvtx \lim
_{t\to\infty}\|u_t\|_{\ell^2(\Z
^d)}(\omega)=0 \Bigr\}.
\end{equation}
\end{corollary}

\begin{pf}
Let $E_1$ denote the event that $\lim_{t\to\infty}\sup_{x\in\Z^d}
|u_t(x)|=0$
and $E_2$ the event that $\lim_{t\to\infty}\|u_t\|_{\ell^2(\Z^d)}=0$.
Because of the real-variable bounds,
$\|u_t\|_{\ell^\infty
(\Z^d)}^2 \le\|u_t\|_{\ell^2(\Z^d)}^2 \le\|u_t\|_{\ell^\infty(\Z
^d)}\,\cdot\,
\|u_t\|_{\ell^1(\Z^d)}$, we have
%
%e8.32 #&#
\begin{equation}
E_1\,\triangle \, E_2 \subseteq \Bigl\{ \omega\dvtx \limsup
_{t\to\infty}\|u_t\|_{\ell^1(\Z
^d)}(\omega) =\infty
\Bigr\}. \label{eq:subset:null}
\end{equation}
We have already noted, however, that $M_t:=\|u_t\|_{\ell^1(\Z^d)}$ defines
a nonnegative martingale, under the conditions of this corollary.
Therefore, the final event in \eqref{eq:subset:null}
is $\P$-null, thanks to Doob's martingale convergence theorem.
Thus we find that $E_1\,\triangle\,
E_2$ is a measurable subset of a $\P$-null set, and is hence $\P$-null.
\end{pf}

%pr8.6 #&#
\begin{proposition}\label{pr:Liap:L2}
Suppose $u_0\in\ell^1(\Z^d)$ and
the random walk $X$ is transient; that is, $\Upsilon(0)<\infty$. Then
%
%e8.33 #&#
\begin{equation}
\label{eq:Liap:L2:1} \limsup_{t\to\infty}\frac{1}t \log\E \bigl(
\|u_t\|_{\ell^1(\Z
^d)}^2 \bigr) \le\inf \bigl\{\beta>0
\dvtx \lip^2\Upsilon(\beta)<1 \bigr\} <\infty.
\end{equation}
If, in addition, $\ell_\sigma^2\Upsilon(0)>1$, then
%
%e8.34 #&#
\begin{equation}
\label{eq:Liap:L2:2} \liminf_{t\to\infty}\frac{1}t \log\E \bigl(
\llVert u_t\rrVert _{%
\ell^2(\Z^d)}^2 \bigr) \ge\inf
\bigl\{\beta>0\dvtx \ell_\sigma^2\Upsilon(\beta)<1 \bigr\}>0.
\end{equation}
\end{proposition}

\begin{pf}
We have already proved a slightly weaker version of \eqref
{eq:Liap:L2:1}. Indeed,
since $\ell^1(\Z^d)\subset\ell^2(\Z^d)$, \eqref{eq:E(E_t)} implies
that
%
%e8.35 #&#
\begin{equation}
\label{eq:Liap:L2:3} \limsup_{t\to\infty}\frac{1}t \log\E \bigl(
\|u_t\|_{\ell^2(\Z
^d)}^2 \bigr) \le\inf \bigl\{\beta>0
\dvtx \lip^2\Upsilon(\beta)<1 \bigr\}.
\end{equation}
Then
\eqref{eq:ell^1} and \eqref{eq:Liap:L2:3} together tell us that
for every $C>\inf\{\beta>0\dvtx \break \lip^2\Upsilon(\beta)<1\}$,
there exists
$K=K(C)\in(0 ,\infty)$ such that
%
%e8.36 #&#
\begin{eqnarray}\qquad
\E \bigl(\|u_t\|_{\ell^1(\Z^d)}^2 \bigr) &\le&
\|u_0\|_{\ell^1(\Z^d)}^2 +\lip^2\int
_0^t \E \bigl(\|u_s
\|_{\ell^2(\Z^d)}^2 \bigr) \,\d s
\nonumber
\\[-8pt]
\\[-8pt]
\nonumber
&\le&\|u_0\|_{\ell^1(\Z^d)}^2 +K\int
_0^t {\mathrm{e}}^{Cs} \,\d s = O \bigl(
{\mathrm{e}}^{(C+o(1))t} \bigr) \qquad\mbox{as $t\to\infty$}.
\end{eqnarray}
Thus follows the first bound of the proposition.

%Next, we argue in a similar way that led to \eqref{eq:ell^2}, and
%deduce that
%\begin{equation}
% \E\left(\left\| u_t\right\|_{\ell^2(\Z^d)}^2\right)
% \ge\left\| \tilde{p}_t*u_0\right\|_{\ell^2(\Z^d)}^2
% +\ell_\sigma^2\int_0^t\bar{P}(t-s)
% \E\left(\left\| u_s\right\|_{\ell^2(\Z^d)}^2\right) \d s.
%\end{equation}
Because of \eqref{eq:LB:L2} and \eqref{eq:LB:p*u}, we find that
%
%e8.37 #&#
\begin{equation}
F(t) := \E \bigl(\llVert u_t\rrVert _{\ell^2(\Z^d)}^2
\bigr)\qquad (t\ge0)
\end{equation}
solves the renewal inequality
%
%e8.38 #&#
\begin{equation}
F(t) \ge g(t) + \int_0^t h(t-s)F(s) \,\d s\qquad (t
\ge0),
\end{equation}
where
%
%e8.39 #&#
\begin{equation}
g(t) := u_0^2(x_0)\bar{P}(t),\qquad h(t):=
\ell_\sigma^2\bar{P}(t)\qquad (t\ge0).
\end{equation}
A comparison result (Lemma~\ref{lem:RE:comparison}) tells us that
$F(t)\ge f(t)$ for all $t\ge0$, where $f$ is the solution to the renewal
equation
%
%e8.40 #&#
\begin{equation}
f(t) = g(t) + \int_0^t h(t-s)f(s) \,\d s\qquad  (t
\ge0).
\end{equation}
The condition that $\ell_\sigma^2\Upsilon(0)>1$ is equivalent to
$\int_0^\infty h(t) \,\d t>1$. Because of transience $[\Upsilon
(0)<\infty]$ and the fact that $\Upsilon(\beta)$ is strictly decreasing
and continuous,
we can find $\beta^{*}>0$ such that $\int_0^\infty\exp(-\beta^{*}
t)h(t) \,\d t=1$. Note that $f_{\beta^{*}}(t):=\exp(-\beta^{*} t)f(t)$
solves the renewal equation
%
%e8.41 #&#
\begin{equation}
f_{\beta^*} (t) = g_{\beta^*}(t) + \int_0^t
h_{\beta^*}(t-s) f_{\beta
^*}(s) \,\d s\qquad  (t\ge0),
\end{equation}
where
$g_{\beta^*}(t):=\exp(-\beta^* t)g(t)$ and $h_{\beta^*}(t):=\exp
(-{\beta^*} t)h(t)$.
Since $h_{\beta^*}$ is a probability density function, and $g_{\beta
^*}$ is nonincreasing
[see \eqref{eq:P:bar:FT}], Blackwell's key renewal theorem \cite
{FellerVolII} implies
that
%
%e8.42 #&#
\begin{eqnarray}\qquad
\liminf_{t\to\infty}{\mathrm{e}}^{-\beta^* t}F(t) &\ge&\lim
_{t\to
\infty}f_{\beta^*}(t) = \biggl(\int_0^\infty
sh_{\beta^*}(s) \,\d s \biggr)^{-1}\,\cdot\,\int_0^\infty
g_{\beta^*}(s) \,\d s
\nonumber
\\[-8pt]
\\[-8pt]
\nonumber
&=& u_0^2(x_0)\ell_\sigma^{-2}
\biggl(\int_0^\infty s{\mathrm{e}}^{-\beta^* s}
\bar{P}(s) \,\d s \biggr)^{-1}\,\cdot\,\Upsilon\bigl(\beta^*\bigr).
\end{eqnarray}
Since $\bar{P}(s)\le1$, the right-most quantity is at least
$u_0^2(x_0) \ell_\sigma^{-2} (\beta^{*})^2\Upsilon(\beta^*)>0$. This
completes the proof of \eqref{eq:Liap:L2:2}. Note that we have used the
fact that $\Upsilon(\beta)$ is continuous in $\beta$ and strictly
decreasing, so that $\beta^{*}=\inf\{\beta>0\dvtx\ell_{\sigma
}^2\Upsilon
(\beta)<1\}>0$.
%Therefore, we can optimize over all $\beta>0$ such that
%$\ell_\sigma^2\Upsilon(\beta)=1$ and deduce
%\eqref{eq:Liap:L2:2}, whence the result.
\end{pf}

%pr8.7 #&#
\begin{proposition}\label{pr:L2:WeakDisorder}
If $u_0 \in\ell^1(\Z^d)$ and $\lip^2\Upsilon(0)<1$, then\break 
$\lim_{t\to\infty}\E(\|u_t\|_{\ell^2(\Z^d)}^2)=0$.
Furthermore, as $t\to\infty$,
%
%e8.43 #&#
\begin{eqnarray}
\label{ass:assert} %
\bar{P}(t) &=&O \bigl(\E \bigl( \|u_t
\|_{\ell^2(\Z^d)}^2 \bigr) \bigr)\quad \mbox{and}
\nonumber
\\[-8pt]
\\[-8pt]
\nonumber
\qquad\E \bigl(\|u_t\|_{\ell^2(\Z^d)}^2 \bigr)& =&O
\bigl(t^{-\alpha
} \bigr)\qquad \mbox{for all $\alpha\ge0$ such that $\bar{P}(t)=O
\bigl(t^{-\alpha}\bigr)$}. %
\end{eqnarray}
\end{proposition}

\begin{pf}
The first assertion of \eqref{ass:assert} is simple to prove; in fact,\break 
$\E(\|u_t\|_{\ell^2(\Z^d)}^2) \ge
[u_0(x_0)]^2 \bar{P}(t)$ $(t\ge0)$ for any $x_0\in\Z^d$
and all $t>0$;
see \eqref{eq:LB:L2} and \eqref{eq:LB:p*u}. We concentrate our efforts
on the remaining statements.

Thanks to \eqref{eq:ell^2},
%
%e8.44 #&#
\begin{equation}\qquad
\E \bigl( \| u_t\|_{\ell^2(\Z^d)}^2 \bigr) \le
\bar{P}(t) \|u_0\| _{\ell
^1(\Z^d)}^2 +
\lip^2\int_0^t\bar{P}(t-s)\E
\bigl( \|u_s\|_{\ell^2(\Z
^d)}^2 \bigr) \,\d s.
\end{equation}
That is, $F(t):=\E(\|u_t\|_{\ell^2(\Z^d)}^2)$ is a sub solution to a
renewal equation;
namely,
%
%e8.45 #&#
\begin{equation}
\label{eq:subcrit:Fgh} F(t) \le g(t) +\int_0^t
h(t-s)F(s) \,\d s\qquad (t\ge0),
\end{equation}
for
%
%e8.46 #&#
\begin{equation}
g(t) := \bar{P}(t)\|u_0\|_{\ell^1(\Z^d)}^2,\qquad h(t) :=
\lip^2\bar{P}(t).
\end{equation}
A comparison lemma (Lemma~\ref{lem:RE:comparison}) shows that
$0\le F(t)\le f(t)$ for all $t\ge0$, where
%
%e8.47 #&#
\begin{equation}
f(t) = g(t) +\int_0^t h(t-s)f(s) \,\d s\qquad (t
\ge0).
\end{equation}
Therefore, it remains to prove that $f(t)\to0$ as $t\to\infty$. It is
easy, as
well as classical, that we can write $f$ in terms of the renewal
function of $h$; that is,
%
%e8.48 #&#
\begin{equation}
\label{eq:f:g:h} f(t) = g(t) + \sum_{n=0}^\infty
\int_0^t h^{*(n)}(s)g(t-s) \,\d s\qquad (t
\ge0),
\end{equation}
where $h^{*(1)}(t):=\int_0^t h(t-s)h(s) \,\d s$ denotes the convolution
of $h$
with itself, and
$h^{*(k+1)}(t) :=\int_0^t h^{*(k)}(t-s)h(s) \,\d s$ for all $k\ge0$.
We might note that $g(t)\le g(0)=\|u_0\|_{\ell^2(\Z^d)}^2$
because $\bar{P}$ is nonincreasing [see \eqref{eq:P:bar:FT}] and one
at zero.
Therefore,
%
%e8.49 #&#
\begin{eqnarray}
\label{eq:h*h} %
0&\le&\int_0^t
h^{*(n)}(s)g(t-s) \,\d s \le\|u_0\|_{\ell^2(\Z^d)}^2
\int_0^\infty h^{*(n)}(s) \,\d s
\nonumber
\\
&\le&\|u_0\|_{\ell^2(\Z^d)}^2 \biggl(\int
_0^\infty h(s) \,\d s \biggr)^{n+1}\qquad
\mbox{[Young's inequality]}
\\
&=&\|u_0\|_{\ell^2(\Z^d)}^2 \bigl( \lip^2
\Upsilon(0) \bigr)^{n+1}.\nonumber
\end{eqnarray}
\upqed\end{pf}
It is not hard to see that $\lim_{t\to\infty}g(t)=
\lim_{t\to\infty}\bar{P}(t)=0$; this follows from
\eqref{eq:P:bar:FT} and the monotone convergence theorem. Because
$\lip^2\Upsilon(0)<1$, we can deduce from \eqref{eq:h*h} and \eqref
{eq:f:g:h},
in conjunction with the dominated convergence theorem, that $f(t)$---hence
$F(t)=\E(\|u_t\|_{\ell^2(\Z^d)}^2)$---converges to zero as $t\to
\infty$.

It remains to prove the second assertion in \eqref{ass:assert}.
With this in mind, let us suppose $\bar{P}$ satisfies the following:
There exists $c\in(0 ,\infty)$
and $\alpha\in[0 ,\infty)$
such that
%
%e8.50 #&#
\begin{equation}
\bar{P}(t) \le c(1+t)^{-\alpha},
\end{equation}
for there is nothing to consider otherwise. We aim to prove that
%
%e8.51 #&#
\begin{equation}
\label{eq:UB:powerdecay} \E \bigl(\|u_t\|_{\ell^2(\Z^d)}^2
\bigr) \le\operatorname {const}\cdot\, (1+t)^{-\alpha},
\end{equation}
for some finite constant that does not depend on $t$.
This proves the proposition.

Define
$F_k(t):=\E(\|u_t^{(k)}\|_{\ell^2(\Z^d)}^2)$, where $u^{(k)}$ denotes
the $k$th approximation to $u$ via Picard's iteration \eqref
{eq:Picard}, starting at
$u^{(0)}_t(x)\equiv0$. We can write \eqref{eq:ell^2},
in short hand, as follows:
%
%e8.52 #&#
\begin{equation}
\label{eq:ell^2:recur} F_{n+1}(t) \le\bar{P}(t)\|u_0
\|_{\ell^1(\Z^d)}^2 + \lip^2\int_0^t
\bar {P}(t-s) F_n(s) \,\d s.
\end{equation}
Now let us choose and fix $\varepsilon\in(0 ,1)$ and write
%
%e8.53 #&#
\begin{eqnarray}
&&\int_0^t \bar{P}(t-s)
F_n(s) \,\d s\nonumber\\
&&\qquad= \int_{t\varepsilon}^t
\bar{P}(s) F_n(t-s) \,\d s +\int_{t(1-\varepsilon)}^t
\bar{P}(t-s) F_n(s) \,\d s
\nonumber
\\[-8pt]
\\[-8pt]
\nonumber
&&\qquad\le c\int_{t\varepsilon}^t \frac{F_n(t-s)}{(1+s)^\alpha} \,\d s+
\sup_{w\ge0} \bigl[(1+w)^\alpha F_n(w)
\bigr]\int_{t(1-\varepsilon)}^t \frac{\bar{P}(t-s)}{(1+s)^\alpha} \,\d s
\\
&&\qquad\le\frac{c}{\varepsilon^\alpha(1+t)^\alpha}\int_0^\infty
F_n(s) \,\d s+ \sup_{w\ge0} \bigl[(1+w)^\alpha
F_n(w) \bigr]\frac{\Upsilon(0)}{
(1-\varepsilon)^\alpha(1+t)^\alpha}.\nonumber
\end{eqnarray}
The proof of Proposition~\ref{pr:L2:int} shows that
%
%e8.54 #&#
\begin{equation}
\label{eq:F_n:recur} \sup_{n\ge0} \int_0^\infty
F_n(s) \,\d s \le\frac{\|u_0\|_{\ell^1(\Z
^d)}^2\Upsilon(0)}{
1-\lip^2\Upsilon(0)}.
\end{equation}
Consequently,
%
%e8.55 #&#
\begin{equation}
R_k:= \sup_{w\ge0} \bigl[ (1+w)^\alpha
F_k(w) \bigr]\qquad (k\ge0)
\end{equation}
satisfies
%
%e8.56 #&#
\begin{equation}
R_{n+1} \le A + R_n\frac{\lip^2\Upsilon(0)}{(1-\varepsilon)^\alpha} \qquad\mbox{for all $n
\ge0$},
\end{equation}
where
%
%e8.57 #&#
\begin{equation}
A = A(\varepsilon) := c\|u_0\|_{\ell^1(\Z^d)}^2 +
\frac{c\|u_0\|_{\ell^1(\Z^d)}^2\lip^2\Upsilon(0)}{\varepsilon
^\alpha(
1-\lip^2\Upsilon(0))}.
\end{equation}
Since $\lip^2\Upsilon(0)<1$, we can choose $\varepsilon$ sufficiently
close to
zero to ensure that $\lip^2\Upsilon(0)<(1-\varepsilon)^{1+\alpha}$.
For this
particular $\varepsilon$, we find that $R_{n+1}\le A + (1-\varepsilon)R_n$
for all $n$. Since $R_0=0$, this proves that $\sup_{n\ge0} R_n\le
A/\varepsilon$.
Equation \eqref{eq:UB:powerdecay}---whence the proposition---follows
from the
latter inequality and Fatou's lemma.
%\end{pf}

%s9 #&#
\section{Proof of Theorem \texorpdfstring{\protect\ref{th:as}}{2.7}} \label{sec:thm1.6}
Let us begin with an elementary real-variable inequality.

%le9.1 #&#
\begin{lemma}\label{lem:arith}
For all real numbers $k\ge2$ and $x,y,\delta>0$,
%
%e9.1 #&#
\begin{equation}
\label{eq:arith} (x+y)^k \le(1+\delta)^{k-1} x^k
+ \biggl(\frac{1+\delta}{\delta
} \biggr)^{k-1} y^k.
\end{equation}
\end{lemma}

This is a consequence of Jensen's inequality when $\delta=1$.
We are interested in the regime $\delta\ll1$.
\begin{pf}
The function $f(z) := (z+1)^k - (1+\delta)^{k-1} z^k$
$(z>0)$ is maximized at $z_* := \delta^{-1}$,
and $\max_z f(z) = f(z_*)=\{ (1+\delta)/\delta\}^{k-1}$; that is,
$f(x) \le\{(1+\delta)/\delta\}^{k-1}$ for all $x>0$. This
is the desired result when $y=1$.
We can factor the variable $y$ from both sides of \eqref{eq:arith}
in order to reduce
the problem to the previously proved case $y=1$.
\end{pf}

%le9.2 #&#
\begin{lemma}\label{lem:transient:ell^2}
$\int_0^\infty\llVert  p_{s+\tau}-p_s\rrVert _{\ell^2(\Z^d)}^2 \,\d s
\le4\Upsilon(0)\tau^2$ for all $\tau\ge0$.
\end{lemma}

\begin{pf}
We apply the Plancherel theorem and \eqref{eq:CHF} in order to deduce that
%
%e9.2 #&#
\begin{eqnarray}
\llVert p_{s+\tau}-p_s \rrVert
_{\ell^2(\Z^d)}^2 &=& (2\pi)^{-d}\int_{(-\pi,\pi)^d}
\bigl\llvert {\mathrm{e}}^{-(s+\tau)(1-\varphi
(\xi))} -{\mathrm{e}}^{-s(1-\varphi(\xi))}\bigr\rrvert
^2\,\d\xi
\nonumber
\\
&=&(2\pi)^{-d}\int_{(-\pi,\pi)^d}{\mathrm{e}}^{-2s(1-\Re\varphi(\xi))}
\bigl\llvert 1-{\mathrm{e}}^{-\tau(1-\varphi(\xi))}\bigr\rrvert ^2\,\d\xi
\\
&\le&\frac{4\tau^2}{(2\pi)^d}\int_{(-\pi,\pi)^d} {\mathrm{e}}^{-2s(1-\Re\varphi
(\xi))} \,\d
\xi.\nonumber %
\end{eqnarray}
Integrate $[\d s]$ to finish; compare with \eqref{eq:Upsilon}.
\end{pf}

Recall that $z_k$ denotes the optimal constant in BDG inequality
\eqref{eq:Davis}.

%le9.3 #&#
\begin{lemma}\label{lem:ell^k:subcrit}
If $k\in(2 ,\infty)$ satisfies $z_k \lip\sqrt{\Upsilon
(0)}<(1+\delta
)^{-(k-1)/k}$
for some $\delta>0$, then
%
%e9.3 #&#
\begin{equation}
\sup_{t\ge0}\E \Bigl( \sup_{x\in\Z^d}\bigl|u_t(x)\bigr|^k
\Bigr)\le \sup_{t\ge0}\E \bigl( \llVert u_t\rrVert
_{\ell^k(\Z^d)}^k \bigr)<\infty.
\end{equation}
\end{lemma}

\begin{pf}
Let $u^{(0)}_t(x):=u_0(x)$, and define $u^{(n)}$ to be the $n$th step
Picard approximation to $u$, as in \eqref{eq:Picard}.
Define
%
%e9.4 #&#
\begin{equation}
\bar{M}_t^{(n)} := \E \bigl( \bigl\llVert
u^{(n)}_t\bigr\rrVert _{\ell^k(\Z
^d)}^k \bigr)\qquad
\mbox{for all $t\ge0$ and $k\ge1$}.
\end{equation}
Then we can apply Lemma~\ref{lem:arith} and
write%, in analogy with \eqref{eq:M},
%
%e9.5 #&#
\begin{equation}
\label{eq:bar{M}} \bar{M}_t^{(n+1)} \le \biggl(
\frac{1+\delta}{\delta} \biggr)^{k-1}\sum_{x\in\Z^d}
I_x + (1+\delta)^{k-1}\sum_{x\in\Z^d}
J_x,
\end{equation}
where $I_x$ and $J_x$ were defined earlier in \eqref{eq:IJ}. One estimates
$\sum_{x\in\Z^d}I_x$ via Jensen's inequality, using $p_t(\bullet-x)$
as the base measure,
in order to find that
%
%e9.6 #&#
\begin{equation}
\label{eq:bar:Ix} \sum_{x\in\Z^d} I_x \le
\|u_0\|_{\ell^k(\Z^d)}^k.
\end{equation}
In order to estimate $\sum_{x\in\Z^d}J_x$, we define---for all
$(t ,x)\in\R_+\times\Z^d$---a Borel measure $\rho_{t,x}$ on $\R
_+\times\Z^d$
as follows:
%
%e9.7 #&#
\begin{equation}
\label{eq:rho} \rho_{t,x}(\d s \,\d y) := \bigl[p_{t-s}(y-x)
\bigr]^2 \1_{[0,t]}(s) \,\d s \chi(\d y);
\end{equation}
where $\chi$ denotes the counting measure on $\Z^d$. Because of
the transience of $X-X'$, the measure $\rho_{t,x}$ is finite; in fact,
%
%e9.8 #&#
\begin{equation}
\label{mass:rho} \rho_{t,x}\bigl(\R_+\times\Z^d\bigr) = \int
_0^t\|p_s\|_{\ell^2(\Z^d)}^2
\,\d s = \int_0^t\bar{P}(s) \,\d s \le
\Upsilon(0).
\end{equation}
Therefore, we apply \eqref{eq:Jx:prec} and Jensen's
inequality, in conjunction, in order to see that
%
%e9.9 #&#
\begin{eqnarray}
J_x &\le& z_k^k \biggl(
\lip^2\int_{[0,t]\times\Z^d} \bigl\{ \E \bigl(\bigl\llvert
u^{(n)}_s(y)\bigr\rrvert ^k \bigr) \bigr
\}^{2/k} \rho_{t,x}(\d s \,\d y) \biggr)^{k/2}
\nonumber
\\[-8pt]
\\[-8pt]
\nonumber
&\le&(z_k\lip)^k \bigl[\Upsilon(0)\bigr]^{(k-2)/2}
\int_{[0,t]\times\Z^d} \E \bigl( \bigl\llvert u^{(n)}_s(y)
\bigr\rrvert ^k \bigr) \rho_{t,x}(\d s \,\d y).
\end{eqnarray}
Thus
%
%e9.10 #&#
\begin{eqnarray}
\sum_{x\in\Z^d}J_x &
\le&(z_k\lip)^k \bigl[\Upsilon(0)\bigr]^{(k-2)/2}\int
_0^t \bar{P}(t-s)\E \bigl(\bigl\llVert
u^{(n)}_s\bigr\rrVert _{
\ell^k(\Z^d)}^k
\bigr) \,\d s
\nonumber
\\[-8pt]
\\[-8pt]
\nonumber
&\le& \bigl(z_k\lip\sqrt{\Upsilon(0)} \bigr)^k\,\cdot\, \sup
_{r\ge0}\E \bigl(\bigl\llVert u^{(n)}_r
\bigr\rrVert _{
\ell^k(\Z^d)}^k \bigr),
\end{eqnarray}
thanks to \eqref{eq:Upsilon}.

In summary, \eqref{eq:bar{M}} has the following consequence: For all
$n\ge0$,
%
%e9.11 #&#
\begin{eqnarray}\qquad
&&\sup_{t\ge0}\bar{M}^{(n+1)}_t
\nonumber
\\[-8pt]
\\[-8pt]
\nonumber
&&\qquad\le \biggl(\frac{1+\delta}{\delta} \biggr)^{k-1} \|u_0
\|_{\ell^k(\Z^d)}^k + (1+\delta)^{k-1} \bigl(z_k
\lip\sqrt{\Upsilon(0)} \bigr)^k \sup_{t\ge0}
\bar{M}^{(n)}_t.
\end{eqnarray}
Since $(1+\delta)^{k-1}(z_k\lip\sqrt{\Upsilon(0)})^k<1$
and $\sup_{t\ge0}\bar{M}^{(0)}_t=\|u_0\|_{\ell^k(\Z^d)}^k$, this\break 
shows that
$C:=\sup_{n\ge0}\sup_{t\ge0}\bar{M}^{(n)}_t<\infty$. Fatou's lemma
now implies
half of the result, since it shows that $\E(\|u_t\|_{\ell^k(\Z
^d)}^k)\le
\liminf_{n\to\infty}\E(\|u_t^{(n)}\|_{\ell^k(\Z^d)}^k)\le C$. The
remainder of the
proposition follows simply because $\|\bullet\|_{\ell^\infty(\Z
^d)}\le
\|\bullet
\|_{\ell^k(\Z^d)}$.
\end{pf}

%pr9.4 #&#
\begin{proposition}\label{pr:ell^k:modulus}
Assume that $\limsup_{t\to\infty}t^{\alpha}\P\{X_t=X'_t\}<1$
for some \mbox{$\alpha>1$}, where $X$ and $X'$ are two independent
random walks with generator $\mathscr{L}$.
If $k\in(2 ,\infty)$ satisfies $z_k \lip\sqrt{\Upsilon(0)}<
(1+\delta)^{-(k-1)/k}$ for some $\delta>0$, then
there exists a finite constant $A$---depending only on
$\delta$, $\lip$, $\Upsilon(0)$ and $\|u_0\|_{\ell^1(\Z
^d)}$---such that
%
%e9.12 #&#
\begin{equation}
\label{eq:ell^k:modulus} \E \bigl( \llVert u_{t+\tau}-u_t\rrVert
_{\ell^k(\Z^d)}^k \bigr) \le\frac{A\tau^{k/2}}{(1+t)^{\alpha}}\qquad \mbox{for every $t,
\tau\ge0$}.
\end{equation}
%
%where $\beta:=\min(\alpha, d/2)$.
Consequently, there exists a H\"older-continuous modification of
the process $t\mapsto u_t(\bullet)$ with values in $\ell^\infty(\Z
^d)$. Moreover,
for this modification, there is a finite constant $A'$---depending
only on
$\delta$, $\lip$, $\Upsilon(0)$ and $\|u_0\|_{\ell^1(\Z
^d)}$---such that
%
%e9.13 #&#
\begin{equation}
\E \biggl( \sup_{ s\neq r\in[t,t+1]} \sup_{x\in\Z^d}\biggl
\llvert \frac{u_r(x)-u_s(x)}{|r-s|^\eta}\biggr\rrvert ^k \biggr)\le
\frac{A'}{(1+t)^{\alpha}},
\end{equation}
as long as $0\le\eta< (k-2)/(2k)$.
%where `` $\sup_I$'' denotes the supremum over all subintervals $I$ of
%$\R_+$ that have length $1$.
\end{proposition}

\begin{pf}
Thanks to Lemma~\ref{lem:ell^k:subcrit}, $\|u_t\|_{\ell^k(\Z^d)}$
has a finite $k$th moment. This observation justifies the use of these
moments in the ensuing discussion. Now we begin our proof in earnest.

The proof requires us to make a few small adjustments to the derivation
of Lemma~\ref{lem:modulus}; specifically we now incorporate the
fact that $\lip^2\Upsilon(0)<1$ into that proof. Therefore, we mention
only the required
changes.

We use the notation of the proof of Lemma~\ref{lem:modulus} and write
%
%e9.14 #&#
\begin{equation}
\bigl\{ \E \bigl(\bigl|u_{t+\tau}(x)-u_t(x)\bigr|^k \bigr)
\bigr\}^{1/k} \le |Q_1|+Q_2+Q_3,
\end{equation}
whence
%
%e9.15 #&#
\begin{equation}
\E \bigl(\bigl|u_{t+\tau}(x)-u_t(x)\bigr|^k \bigr)
\le3^{k-1} \bigl( |Q_1|^k + Q_2^k
+ Q_3^k \bigr).
\end{equation}
Note that
%
%e9.16 #&#
\begin{eqnarray}
\nonumber
&&\sum_{x\in\Z^d}|Q_1|^k\\[-2pt]
&&\qquad \le\sum_{x\in\Z^d} \biggl( \sum
_{y\in\Z^d} u_0(y) \bigl\llvert p_{t+\tau}(y-x)-p_t(y-x)
\bigr\rrvert \biggr)^k
\nonumber
\\[-8pt]
\\[-8pt]
\nonumber
&&\qquad\le\|u_0\|_{\ell^1(\Z^d)}^{k-1}\,\cdot\, \sum
_{x\in\Z^d}\sum_{y\in\Z^d}
u_0(y) \bigl\llvert p_{t+\tau}(y-x)-p_t(y-x)
\bigr\rrvert ^k
\\
\nonumber
&&\qquad=\|u_0\|_{\ell^1(\Z^d)}^k\,\cdot\,\llVert
p_{t+\tau} - p_t\rrVert _{\ell^k(\Z^d)}^{k} \le
\|u_0\|_{\ell^1(\Z^d)}^k\,\cdot\, \llVert
p_{t+\tau}-p_t\rrVert _{\ell^2(\Z^d)}^{k} ,
\end{eqnarray}
thanks to Jensen's inequality.
We observe that
%
%e9.17 #&#
\begin{eqnarray}
\label{eq:p:diff}
&&
\| p_{t+\tau}-p_t
\|_{\ell^2(\Z^d)}^2 \nonumber\\
&&\qquad= (2\pi)^{-d} \int
_{[-\pi
,\pi]^d} \bigl\llvert {\mathrm{e}}^{-t(1-\varphi(\xi))}\bigr\rrvert
^2 \bigl\llvert {\mathrm{e}}^{-\tau(1-\varphi(\xi))}-1\bigr\rrvert ^2
\,\d\xi
\nonumber
\\[-8pt]
\\[-8pt]
\nonumber
&&\qquad\le\operatorname{const}\,\cdot\,\tau^2 \int
_{[-\pi,\pi]^d} \bigl\llvert {\mathrm{e}}^{-t(1-\varphi(\xi))}\bigr\rrvert
^2\,\d\xi =\operatorname{const}\,\cdot\,\tau^2\P\bigl
\{X_t=X'_t\bigr\}
\\
&&\qquad\le\operatorname{const}\,\cdot\,\frac{\tau^2}{(1+t)^{\alpha} }.\nonumber
\end{eqnarray}
Consequently,
%
%e9.18 #&#
\begin{equation}
\label{eq:Q_1:bis} \sum_{x\in\Z^d}|Q_1|^k
\le\frac{\operatorname{const}\,\cdot\,\tau
^k}{(1+t)^{\alpha k/2}}.
\end{equation}

We estimate $Q_2$ slightly differently from the proof of Lemma~\ref
{lem:modulus} as well.

For every $(t ,x)\in\R_+\times\Z^d$, let us define a similar Borel measure
$R_{t,x}$ to $\rho_{t,x}$ [see \eqref{eq:rho}] as follows:
%
%e9.19 #&#
\begin{equation}
R_{t,x}(\d s \,\d y) := \bigl[p_{t+\tau-s}(y-x)-p_{t-s}(y-x)
\bigr]^2 \1_{[0,t]}(s) \,\d s \chi(\d y).
\end{equation}
Now we reexamine the first line of \eqref{eq:Q2:prec} and note that
%
%e9.20 #&#
\begin{eqnarray}\qquad
Q_2^2&\le&(z_k\lip)^2
\sum_{y\in\Z^d}\int_0^t
\bigl[ p_{t+\tau
-s}(y-x)-p_{t-s} (y-x) \bigr]^2 \bigl
\{ \E \bigl(\bigl|u_s(y)\bigr|^k \bigr) \bigr\}^{2/k}\,\d s
\nonumber
\\[-8pt]
\\[-8pt]
\nonumber
& =&(z_k\lip)^2\int_{\R_+\times\Z^d} \bigl\{ \E
\bigl(\bigl|u_s(y)\bigr|^k \bigr) \bigr\}^{2/k}
R_{t,x}(\d s \,\d y).
\end{eqnarray}
This follows from \eqref{eq:sigma(0)=0} and \eqref{eq:Q2:prec}, but
we use
the optimal constant $z_k$ in place of the slightly weaker
$2\sqrt{k}$ that came from Lemma~\ref{lem:BDG}.

Lemma~\ref{lem:transient:ell^2} implies that
$R_{t,x}(\R_+\times\Z^d) = \int_0^t \llVert  p_{s+\tau}-p_s\rrVert _{
\ell^2(\Z^d)}^2\,\d s \le4\Upsilon(0)\tau^2$. This bound and
Jensen's inequality together show that
$\sum_{x\in\Z^d} Q_2^k$ is bounded from above by
%
%e9.21 #&#
\begin{eqnarray}\label{eq:Q2:bound}
&&(z_k\lip)^k \bigl(4\Upsilon(0)\tau^2
\bigr)^{(k-2)/2} \sum_{x\in\Z^d}\int
_{\R_+\times\Z^d}\E \bigl(\bigl\llvert u_s(y) \bigr\rrvert
^k \bigr) R_{t,x}(\d s \,\d y)
\nonumber
\\[-8pt]
\\[-8pt]
\nonumber
&&\qquad=(z_k\lip)^k \bigl(4\Upsilon(0)
\tau^2\bigr)^{(k-2)/2} \int_0^t
\llVert p_{t+\tau
-s}-p_{t-s} \rrVert _{\ell^2(\Z^d)}^2
\E \bigl( \|u_s\|_{\ell^k(\Z^d)}^k \bigr) \,\d s.
\end{eqnarray}

By an argument similar to the one used in Proposition~\ref{pr:L2:WeakDisorder},
one is able to show that $\E( \|u_s\|_{\ell^k(\Z^d)}^k) \le
\operatorname{const}\,\cdot\,(1+s)^{-\alpha}$. Here is an outline of the proof:
We can follow the proof of Lemma~\ref{lem:ell^k:subcrit},
but derive a better bound on
$\sum_{x\in\Z^d} I_x \le\bar{P}(t) \|u_0\|_{\ell^1(\Z^d)}^k $,
in order to obtain
%
%e9.22 #&#
\begin{eqnarray}
\label{eq:ell^k:prop8.2} &&\E \bigl(\bigl\|u_t^{(n+1)}
\bigr\|_{\ell^k(\Z^d)}^k \bigr) \nonumber\\
&&\qquad \le \biggl(\frac{1+\delta}{\delta}
\biggr)^{k-1} \bar{P}(t) \|u_0\|_{\ell^1(\Z^d)}^k
\\
&&\qquad\quad{} +(1+\delta)^{k-1}(z_k\lip)^k \bigl[\Upsilon
(0)\bigr]^{(k-2)/2}\int_0^t
\bar{P}(t-s)\E \bigl(\bigl\llVert u^{(n)}_s\bigr\rrVert
_{
\ell^k(\Z^d)}^k \bigr) \,\d s.\nonumber
\end{eqnarray}
From here, we proceed along similar lines, as we did from \eqref
{eq:ell^2:recur} onward.
We follow the proof of Proposition~\ref{pr:L2:int}, using \eqref
{eq:ell^k:prop8.2},
in order to derive the following analog of \eqref{eq:F_n:recur}:
%
%e9.23 #&#
\begin{equation}
\sup_{n\ge0} \int_0^{\infty}
F_n(s) \,\d s \le \frac{((1+\delta)/\delta)^{k-1}
\|u_0\|_{\ell^1(\Z^d)}^k \Upsilon(0) }{1- (1+\delta)^{k-1}(z_k\lip
\Upsilon(0))^k },
\end{equation}
where $F_n(t):=\E(\|u_t^{(n+1)}\|_{\ell^k(\Z^d)}^k)$.
In this way, we can obtain the bound
$\E( \|u_s\|_{\ell^k(\Z^d)}^k) \le\operatorname{const}\cdot\,
(1+s)^{-\alpha}$,
as was needed. We use this bound, as well as~\eqref{eq:p:diff} in
\eqref{eq:Q2:bound},
and split the integral into two parts ($0$ to $t/2$ and $t/2$ to $t$),
in order to obtain
the following:
%
%e9.24 #&#
\begin{equation}
\label{eq:Q_2:bis} \sum_{x\in\Z^d} Q_2^k
\le\frac{\operatorname{const}\cdot\,\tau
^k}{(1+t)^{\alpha}}.
\end{equation}

Finally we estimate $\sum_{x\in\Z^d}Q_3^k$ by first modifying \eqref
{eq:Q_3:prec} as follows:
%
%e9.25 #&#
\begin{eqnarray}
Q_3^2 &\le&(z_k\lip)^2
\sum_{y\in\Z^d} \int_t^{t+\tau}
\bigl[ p_{t+\tau-s}(y-x) \bigr]^2 \bigl\{ \E \bigl( \bigl\llvert
u_s(y)\bigr\rrvert ^k \bigr) \bigr\}^{2/k} \,\d
s
\nonumber
\\[-8pt]
\\[-8pt]
\nonumber
&=&(z_k\lip)^2\int_{\R_+\times\Z^d} \bigl\{ \E
\bigl( \bigl\llvert u_s(y)\bigr\rrvert ^k \bigr) \bigr
\}^{2/k} \mathcal{R}_{t,\tau,x}(\d s \,\d y),
\end{eqnarray}
where the Borel measures $\mathcal{R}_{t,\tau,x}$ are defined in a similar
manner as in \eqref{eq:rho}; that is,
%
%e9.26 #&#
\begin{equation}
\mathcal{R}_{t,\tau,x} (\d s \,\d y) := \sum_{y\in\Z^d}
\bigl[ p_{t+\tau
-s}(y-x) \bigr]^2 \1_{[t,t+\tau]}(s) \,\d s
\chi(\d y).
\end{equation}
Because $\mathcal{R}_{t,\tau,x}(\R_+\times\Z^d) = \int_0^\tau
\bar{P}(s)\,\d s\le\tau$, Jensen's inequality assures us that
%
%e9.27 #&#
\begin{eqnarray}
\label{eq:Q_3:bis} %
\sum_{x\in\Z^d}Q_3^k
&\le&(z_k\lip)^k\tau^{(k-2)/2}\sum
_{x\in\Z^d} \int_{\R_+\times\Z^d} \E \bigl(\bigl\llvert
u_s(y)\bigr\rrvert ^k \bigr) \mathcal{R}_{t,\tau,x}(
\,\d s \,\d y)
\nonumber
\\
&=&(z_k\lip)^k\tau^{(k-2)/2}\int
_t^{t+\tau} \bar{P}(t+\tau-s) \E \bigl( \llVert
u_s\rrVert _{\ell^k(\Z^d)}^k \bigr) \,\d s
\\
&\le&\frac{\operatorname{const}\cdot\,\tau^{{k}/{2}}}{(1+t)^{\alpha}}, \nonumber%
\end{eqnarray}
thanks to the bounds
$\E( \|u_s\|_{\ell^k(\Z^d)}^k) \le\operatorname{const}
\,\cdot\,(1+s)^{-\alpha}$ and $\bar{P}(t+\tau-s)\le1$.
Since $\|u_0\|_{\ell^k(\Z^d)}\le\|u_0\|_{\ell^1(\Z^d)}$,
displays \eqref{eq:Q_1:bis}, \eqref{eq:Q_2:bis} and \eqref{eq:Q_3:bis}
together imply~\eqref{eq:ell^k:modulus}. This yields
the first estimate of the proposition. The remaining
assertions follow~\eqref{eq:ell^k:modulus},
using a suitable form of the Kolmogorov continuity
theorem \cite{RevuzYor}, Theorem~2.1, page 25, and the fact that
$\sup_{x\in\Z^d}|u_t(x)-u_s(x)|\le\|u_t-u_s\|_{\ell^k(\Z^d)}$.
\end{pf}

\begin{pf*}{Proof of Theorem~\ref{th:as}}
We apply Proposition~\ref{pr:L2:WeakDisorder} and Chebyshev's
inequality in conjunction in order to see that
%
%e9.28 #&#
\begin{eqnarray}
\sum_{n=1}^\infty\P \Bigl\{ \sup
_{x\in\Z^d} \bigl|u_{n}(x)\bigr| > \varepsilon \Bigr\} &\le&
\frac{1}{\varepsilon^2} \sum_{n=1}^\infty \E
\bigl(\|u_{n}\|_{\ell^2(\Z^d)}^2 \bigr)
\nonumber
\\[-8pt]
\\[-8pt]
\nonumber
&\le&\frac{\operatorname{const}}{\varepsilon^2}\,\cdot\, \sum_{n=1}^\infty
n^{-\alpha}<\infty.
\end{eqnarray}
Therefore, the Borel--Cantelli lemma implies that
%
%e9.29 #&#
\begin{equation}
\label{eq:dissipate:subseq} \lim_{n\to\infty}
\sup_{x\in\Z^d}\bigl|u_{n}(x)\bigr|=0\qquad
\mbox{a.s.}
\end{equation}
We next note that the Burkholder's constants $z_k$ vary
continuously for $k\ge2$, and $z_2=1$ is the minimum;
see Davis \cite{Davis}. Davis \cite{Davis} obtains $z_k$ as
the largest positive zero of the parabolic cylinder function
of parameter $k$ and this varies continuously in $k$; see
Abramowitz and Stegun \cite{AS}.

If $\lip\sqrt{\Upsilon(0)}<1$, we can find $k>2$ and $\delta>0$
such that
%
%e9.30 #&#
\begin{equation}
\label{cond:k} z_k\lip\sqrt{\Upsilon(0)}<(1+\delta)^{-(k-1)/k}.
\end{equation}
We can now use Proposition~\ref{pr:ell^k:modulus} (with $\eta=0$)
along with Chebyshev's inequality to control the spacings
%
%e9.31 #&#
\begin{eqnarray}
&&\P \Bigl\{ \sup_{s\in[n,n+1]}\sup_{x\in\Z^d} \bigl
\llvert u_s(x) - u_{n}(x)\bigr\rrvert >\varepsilon \Bigr
\}
\nonumber
\\[-8pt]
\\[-8pt]
\nonumber
&&\qquad\le\frac{1}{\varepsilon^k} \E \Bigl( \sup_{s\in[n,n+1]} \sup
_{x\in\Z^d}\bigl\llvert u_t(x)-u_s(x)
\bigr\rrvert ^k \Bigr) = O \bigl( n^{-\alpha} \bigr)\qquad \mbox{as $n
\to\infty$}.
\end{eqnarray}
We may use the Borel--Cantelli lemma and \eqref{eq:dissipate:subseq}
in order to deduce that $\lim_{t\to\infty}\sup_{x\in\Z
^d}|u_t(x)|= 0$ a.s.
Thanks to this fact, Corollary~\ref{co:iff} implies the seemingly
stronger assertion that
$\lim_{t\to\infty}\|u_t\|_{\ell^2(\Z^d)}^2= 0$ a.s., and completes the
proof.
\end{pf*}

%sA #&#
\begin{appendix}\label{app}
\section*{Appendix: Some renewal theory}

In this appendix we state and prove a few facts from (linear) renewal
theory. These facts ought to be well known, but we have not succeeded
to find concrete references, and so will describe them in some detail.\vadjust{\goodbreak}

Let us suppose that the functions $h,g\dvtx(0 ,\infty)\to\R_+$ are
locally integrable (say) and pre-defined,
and let us look for a measurable solution $f\dvtx(0 ,\infty)\to\R_+$
to the renewal equation
%
%eA.1 #&#
\setcounter{equation}{0}
\begin{equation}
\label{RE} f(t) = g(t) + \int_0^t
h(t-s)f(s) \,\d s\qquad (t\ge0).
\end{equation}
If $h\in L^1(0 ,\infty)$, then this is a classical subject \cite{FellerVolII}.
For a more general treatment, we may proceed with Picard's iteration:
Let $f^{(0)}(t)\dvtx(0 ,\infty)\to\R_+$ be a fixed measurable
function, and iteratively define
%
%eA.2 #&#
\begin{equation}
f^{(n+1)}(t) := g(t) + \int_0^t
h(t-s)f^{(n)}(s) \,\d s\qquad  (t> 0, n\ge0).
\end{equation}

%leA.1 #&#
\begin{lemmaa}\label{lem:RE1}
Suppose that there exists a constant $\beta\in\R$ that satisfies the
following three conditions:
\textup{(i)} $\gamma:= \sup_{t\ge0}[\exp(-\beta t)g(t)]<\infty$;
\textup{(ii)} $\rho:=\int_0^\infty\exp(-\beta t)h(t) \,\d t<1$;
\textup{(iii)} $\sup_{t\ge0}[\exp(-\beta t)f^{(0)}(t)]<\infty$.
Then \eqref{RE} has a unique nonnegative
solution $f$ that satisfies the following:
%
%eA.3 #&#
\begin{equation}
\label{eq:RE:UB(f)} f(t) \le\frac{\gamma{\mathrm{e}}^{\beta t}}{1-\rho}\qquad (t\ge0).
\end{equation}
Moreover, $\lim_{n\to\infty} \sup_{t\ge0}(
{\mathrm{e}}^{-\beta t}| f^{(n)}(t)-f(t)|)=0$.
\end{lemmaa}

\begin{pf}
Choose such a $\beta\in\R$, and define
%
%eA.4 #&#
\begin{equation}
\gamma:=\sup_{t\ge0} \bigl[ {\mathrm{e}}^{-\beta t}g(t) \bigr],\qquad
\rho:=\int_0^\infty{\mathrm{e}}^{-\beta t}h(t) \,\d
t<1
\end{equation}
and
%
%eA.5 #&#
\begin{equation}\qquad
C_k := \sup_{t\ge0} \bigl({\mathrm{e}}^{-\beta t}
f^{(k)}(t) \bigr), \qquad D_k := \sup_{t\ge0}
\bigl({\mathrm{e}}^{-\beta t}\bigl\llvert f^{(k)}(t) -
f^{(k-1)}(t) \bigr\rrvert \bigr),
\end{equation}
for integers $k\ge1$. Thanks to the definition of the $f^{(k)}$'s,
%
%eA.6 #&#
\begin{equation}
C_{n+1} \le\gamma+ \rho C_n,\qquad D_{n+1} \le\rho
D_n\qquad (n\ge0).
\end{equation}
Consequently, $\sup_{n\ge0}C_n \le\gamma(1-\rho)^{-1}$
and $D_n=O(\rho^n)$. Since $\sum_{n=0}^\infty D_n<\infty$,
it follows that there exists a function $f$ such that
$\sup_{t\ge0}({\mathrm{e}}^{-\beta t}|f^{(n)}(t)-f(t)|)\to0$ as $n\to
\infty$,
and $\sup_{t\ge0}({\mathrm{e}}^{-\beta t}f(t))\le\sup_{n\ge0}C_n$. These
observations together prove the lemma.
\end{pf}

The following is the main result of this appendix.

%leA.2 #&#
\begin{lemmaa}[(Comparison lemma)]\label{lem:RE:comparison}
Suppose there exists $\beta\in\R$ such that:
\textup{(i)}~$\gamma:= \sup_{t\ge0}[\exp(-\beta t)g(t)]<\infty$
and \textup{(ii)} $\rho:=\int_0^\infty\exp(-\beta t)h(t) \,\d t<1$; and
let $f$ denote the unique nonnegative solution to \eqref{RE}
that satisfies \eqref{eq:RE:UB(f)}. If $F\dvtx\R_+\to\R_+$ satisfies:
\textup{(a)} $\sup_{t\ge0}[\exp(-\beta t)F(t)]<\infty$ and \textup{(b)}
%
%eA.7 #&#
\begin{equation}
\label{eq:comp1} F(t) \ge g(t) + \int_0^t
h(t-s)F(s) \,\d s\qquad (t\ge0),
\end{equation}
then $f(t)\le F(t)$ for all $t\ge0$. Finally, if we replace condition
\eqref{eq:comp1}
by
%
%eA.8 #&#
\begin{equation}
\label{eq:comp2} F(t) \le g(t) + \int_0^t
h(t-s)F(s) \,\d s \qquad (t\ge0),
\end{equation}
then $f(t)\ge F(t)$ for all $t\ge0$.
\end{lemmaa}

\begin{pf}
We will prove \eqref{eq:comp1}; \eqref{eq:comp2} is proved similarly.

We apply Picard's iteration with initial function $f^{(0)}:=F$
and note that
%
%eA.9 #&#
\begin{equation}
f^{(1)}(t) = g(t) + \int_0^t
h(t-s) F(s) \,\d s \le F(t) \qquad (t\ge0).
\end{equation}
This and induction together show that $f^{(n+1)}(t)\le f^{(n)}(t)$
for all $t\ge0$ and $n\ge0$.
Let $n\to\infty$ to deduce the lemma from Lemma~\ref{lem:RE1}.
\end{pf}
\end{appendix}
\section*{Acknowledgments}
We thank an anonymous referee and an Associate
Editor whose remarks and suggestions have helped improve
the presentation of this paper.

% imsref loaded by akundreckaite, 2014-10-14 08:32:26
% imsref loaded by akundreckaite, 2014-10-14 08:50:54
%

%\begin{appendix}
%\section{}
%\end{appendix}

%\section*{Acknowledgments}

%\begin{supplement}[id=suppA]
%\slink[doi]{10.1214/00-AAPXXXXSUPP} %[doi,text={...}] - jei reikia
%suskaldyti doi
%\sfilename{aapXXXX\_supp.pdf}
%\end{supplement}

%\begin{thebibliography}{99}
%\bibitem[\protect\citeauthoryear{}{}]{r1}
%\bibitem{r1}
%\end{thebibliography}

\printaddresses

\begin{thebibliography}{39}

%b1 ###
%b1 #&#
\bibitem{AS}
%
\begin{bbook}[mr]
\beditor{\bsnm{Abramowitz},~\bfnm{Milton}\binits{M.}} \AND
\beditor{\bsnm{Stegun},~\bfnm{Irene A.}\binits{I.~A.}}, eds.
(\byear{1992}).
\btitle{Handbook of Mathematical Functions with Formulas, Graphs, and
Mathematical Tables}.
\bpublisher{Dover},
\blocation{New York}.
\bid{mr={1225604}}
\end{bbook}
%

\bptok{imsref}%
% NOT OUTPUTED:
% isbn = 0-486-61272-4
% fpage = xiv+1046
\endbibitem

%b2 ###
%b2 #&#
\bibitem{AizenmanWarzel}
%
\begin{barticle}[auto:parserefs-M02]
\bauthor{\bsnm{Aizenman},~\bfnm{M.}\binits{M.}} \AND
\bauthor{\bsnm{Warzel},~\bfnm{S.}\binits{S.}}
(\byear{2011}).
\btitle{Absence of mobility edge for the Anderson random potential on
tree graphs at weak disorder}.
\bjournal{European Physics Letters}
\bvolume{96}
\bpages{37004}.
\end{barticle}
%

\bptok{imsref}%
\endbibitem

%b3 ###
%b3 #&#
\bibitem{ACQ}
%
\begin{barticle}[mr]
\bauthor{\bsnm{Amir},~\bfnm{Gideon}\binits{G.}},
\bauthor{\bsnm{Corwin},~\bfnm{Ivan}\binits{I.}} \AND
\bauthor{\bsnm{Quastel},~\bfnm{Jeremy}\binits{J.}}
(\byear{2011}).
\btitle{Probability distribution of the free energy of the continuum
directed random polymer in {$1+1$} dimensions}.
\bjournal{Comm. Pure Appl. Math.}
\bvolume{64}
\bpages{466--537}.
\bid{doi={10.1002/cpa.20347}, issn={0010-3640}, mr={2796514}}
\end{barticle}
%

\bptok{imsref}%
% NOT OUTPUTED:
% number = 4
% doi = http://dx.doi.org/10.1002/cpa.20347
% coden = CPAMA
% fjournal = Communications on Pure and Applied Mathematics
\endbibitem

%b4 ###
%b4 #&#
\bibitem{Anderson}
%
\begin{barticle}[mr]
\bauthor{\bsnm{Anderson},~\bfnm{William~J.}\binits{W.~J.}}
(\byear{1972}).
\btitle{Local behaviour of solutions of stochastic integral equations}.
\bjournal{Trans. Amer. Math. Soc.}
\bvolume{164}
\bpages{309--321}.
\bid{issn={0002-9947}, mr={0297031}}
\end{barticle}
%

\bptok{imsref}%
% NOT OUTPUTED:
% fjournal = Transactions of the American Mathematical Society
\endbibitem

%b5 ###
%b5 #&#
\bibitem{BQS}
%
\begin{barticle}[mr]
\bauthor{\bsnm{Bal{\'a}zs},~\bfnm{M.}\binits{M.}},
\bauthor{\bsnm{Quastel},~\bfnm{J.}\binits{J.}} \AND
\bauthor{\bsnm{Sepp{\"a}l{\"a}inen},~\bfnm{T.}\binits{T.}}
(\byear{2011}).
\btitle{Fluctuation exponent of the KPZ/stochastic {B}urgers equation}.
\bjournal{J. Amer. Math. Soc.}
\bvolume{24}
\bpages{683--708}.
\bid{doi={10.1090/S0894-0347-2011-00692-9}, issn={0894-0347}, mr={2784327}}
\end{barticle}
%

\bptok{imsref}%
% NOT OUTPUTED:
% number = 3
% doi = http://dx.doi.org/10.1090/S0894-0347-2011-00692-9
% fjournal = Journal of the American Mathematical Society
\endbibitem

%b6 ###
%b6 #&#
\bibitem{BC}
%
\begin{barticle}[mr]
\bauthor{\bsnm{Bertini},~\bfnm{Lorenzo}\binits{L.}} \AND
\bauthor{\bsnm{Cancrini},~\bfnm{Nicoletta}\binits{N.}}
(\byear{1995}).
\btitle{The stochastic heat equation: {F}eynman--{K}ac formula and
intermittence}.
\bjournal{J. Stat. Phys.}
\bvolume{78}
\bpages{1377--1401}.
\bid{doi={10.1007/BF02180136}, issn={0022-4715}, mr={1316109}}
\bptnote{check year}%
\end{barticle}
%

\bptok{imsref}%
% NOT OUTPUTED:
% number = 5-6
% doi = http://dx.doi.org/10.1007/BF02180136
% coden = JSTPSB
% fjournal = Journal of Statistical Physics
\endbibitem

%b7 ###
%b7 #&#
\bibitem{BorCor}
%
\begin{barticle}[mr]
\bauthor{\bsnm{Borodin},~\bfnm{Alexei}\binits{A.}} \AND
\bauthor{\bsnm{Corwin},~\bfnm{Ivan}\binits{I.}}
(\byear{2014}).
\btitle{Moments and {L}yapunov exponents for the parabolic {A}nderson model}.
\bjournal{Ann. Appl. Probab.}
\bvolume{24}
\bpages{1172--1198}.
\bid{doi={10.1214/13-AAP944}, issn={1050-5164}, mr={3199983}}
\end{barticle}
%

\bptok{imsref}%
% NOT OUTPUTED:
% number = 3
% doi = http://dx.doi.org/10.1214/13-AAP944
% fjournal = The Annals of Applied Probability
\endbibitem

%b8 ###
%b8 #&#
\bibitem{Burkholder}
%
\begin{barticle}[mr]
\bauthor{\bsnm{Burkholder},~\bfnm{D.~L.}\binits{D.~L.}}
(\byear{1966}).
\btitle{Martingale transforms}.
\bjournal{Ann. Math. Statist.}
\bvolume{37}
\bpages{1494--1504}.
\bid{issn={0003-4851}, mr={0208647}}
\end{barticle}
%

\bptok{imsref}%
% NOT OUTPUTED:
% fjournal = Annals of Mathematical Statistics
\endbibitem

%b9 ###
%b9 #&#
\bibitem{BDG}
%
\begin{binproceedings}[mr]
\bauthor{\bsnm{Burkholder},~\bfnm{D.~L.}\binits{D.~L.}},
\bauthor{\bsnm{Davis},~\bfnm{B.~J.}\binits{B.~J.}} \AND
\bauthor{\bsnm{Gundy},~\bfnm{R.~F.}\binits{R.~F.}}
(\byear{1972}).
\btitle{Integral inequalities for convex functions of operators on
martingales}.
In \bbooktitle{Proceedings of the {S}ixth {B}erkeley {S}ymposium on
{M}athematical {S}tatistics and {P}robability ({U}niv. {C}alifornia,
{B}erkeley, {C}alif., 1970/1971), {V}ol. {II}: {P}robability Theory}
\bpages{223--240}.
\bpublisher{Univ. California Press},
\blocation{Berkeley, CA}.
\bid{mr={0400380}}
\end{binproceedings}
%

\bptok{imsref}%
\endbibitem

%b10 ###
%b10 #&#
\bibitem{BG}
%
\begin{barticle}[mr]
\bauthor{\bsnm{Burkholder},~\bfnm{D.~L.}\binits{D.~L.}} \AND
\bauthor{\bsnm{Gundy},~\bfnm{R.~F.}\binits{R.~F.}}
(\byear{1970}).
\btitle{Extrapolation and interpolation of quasi-linear operators on
martingales}.
\bjournal{Acta Math.}
\bvolume{124}
\bpages{249--304}.
\bid{issn={0001-5962}, mr={0440695}}
\end{barticle}
%

\bptok{imsref}%
% NOT OUTPUTED:
% fjournal = Acta Mathematica
\endbibitem

%b11 ###
%b11 #&#
\bibitem{CK}
%
\begin{barticle}[mr]
\bauthor{\bsnm{Carlen},~\bfnm{Eric}\binits{E.}} \AND
\bauthor{\bsnm{Kr{\'e}e},~\bfnm{Paul}\binits{P.}}
(\byear{1991}).
\btitle{{$L^p$} estimates on iterated stochastic integrals}.
\bjournal{Ann. Probab.}
\bvolume{19}
\bpages{354--368}.
\bid{issn={0091-1798}, mr={1085341}}
\end{barticle}
%

\bptok{imsref}%
% NOT OUTPUTED:
% url =
%http://links.jstor.org/sici?sici=0091-1798(199101)19:1<354:EOISI>2.0.CO;2-C&origin=MSN
% number = 1
% coden = APBYAE
% fjournal = The Annals of Probability
\endbibitem

%b12 ###
%b12 #&#
\bibitem{CKM}
%
\begin{barticle}[mr]
\bauthor{\bsnm{Carmona},~\bfnm{Rene}\binits{R.}},
\bauthor{\bsnm{Koralov},~\bfnm{Leonid}\binits{L.}} \AND
\bauthor{\bsnm{Molchanov},~\bfnm{Stanislav}\binits{S.}}
(\byear{2001}).
\btitle{Asymptotics for the almost sure {L}yapunov exponent for the
solution of the parabolic {A}nderson problem}.
\bjournal{Random Oper. Stoch. Equ.}
\bvolume{9}
\bpages{77--86}.
\bid{doi={10.1515/rose.2001.9.1.77}, issn={0926-6364}, mr={1910468}}
\end{barticle}
%

\bptok{imsref}%
% NOT OUTPUTED:
% number = 1
% doi = http://dx.doi.org/10.1515/rose.2001.9.1.77
% fjournal = Random Operators and Stochastic Equations
\endbibitem

%b13 ###
%b13 #&#
\bibitem{CM94}
%
\begin{barticle}[mr]
\bauthor{\bsnm{Carmona},~\bfnm{Ren{\'e}~A.}\binits{R.~A.}} \AND
\bauthor{\bsnm{Molchanov},~\bfnm{S.~A.}\binits{S.~A.}}
(\byear{1994}).
\btitle{Parabolic {A}nderson problem and intermittency}.
\bjournal{Mem. Amer. Math. Soc.}
\bvolume{108}
\bpages{viii+125}.
\bid{doi={10.1090/memo/0518}, issn={0065-9266}, mr={1185878}}
\end{barticle}
%

\bptok{imsref}%
% NOT OUTPUTED:
% number = 518
% doi = http://dx.doi.org/10.1090/memo/0518
% coden = MAMCAU
% fjournal = Memoirs of the American Mathematical Society
\endbibitem

%b14 ###
%b14 #&#
\bibitem{ChungFuchs}
%
\begin{barticle}[mr]
\bauthor{\bsnm{Chung},~\bfnm{K.~L.}\binits{K.~L.}} \AND
\bauthor{\bsnm{Fuchs},~\bfnm{W.~H.~J.}\binits{W.~H.~J.}}
(\byear{1951}).
\btitle{On the distribution of values of sums of random variables}.
\bjournal{Mem. Amer. Math. Soc.}
\bvolume{1951}
\bpages{12}.
\bid{issn={0065-9266}, mr={0040610}}
\end{barticle}
%

\bptok{imsref}%
% NOT OUTPUTED:
% number = 6
% fjournal = Memoirs of the American Mathematical Society
\endbibitem

%b15 ###
%b15 #&#
\bibitem{Conus}
%
\begin{barticle}[mr]
\bauthor{\bsnm{Conus},~\bfnm{Daniel}\binits{D.}}
(\byear{2013}).
\btitle{Moments for the parabolic {A}nderson model: On a result by
{H}u and {N}ualart}.
\bjournal{Commun. Stoch. Anal.}
\bvolume{7}
\bpages{125--152}.
\bid{issn={0973-9599}, mr={3080991}}
\end{barticle}
%

\bptok{imsref}%
% NOT OUTPUTED:
% number = 1
% fjournal = Communications on Stochastic Analysis
\endbibitem

%b16 ###
%b16 #&#
\bibitem{CJK11}
%
\begin{barticle}[mr]
\bauthor{\bsnm{Conus},~\bfnm{Daniel}\binits{D.}},
\bauthor{\bsnm{Joseph},~\bfnm{Mathew}\binits{M.}} \AND
\bauthor{\bsnm{Khoshnevisan},~\bfnm{Davar}\binits{D.}}
(\byear{2012}).
\btitle{Correlation-length bounds, and estimates for intermittent
islands in parabolic {SPDE}s}.
\bjournal{Electron. J. Probab.}
\bvolume{17}
\bpages{15}.
\bid{doi={10.1214/EJP.v17-2429}, issn={1083-6489}, mr={3005720}}
\end{barticle}
%

\bptok{imsref}%
% NOT OUTPUTED:
% doi = http://dx.doi.org/10.1214/EJP.v17-2429
% fjournal = Electronic Journal of Probability
\endbibitem

%b17 ###
%b17 #&#
\bibitem{CJK13}
%
\begin{barticle}[mr]
\bauthor{\bsnm{Conus},~\bfnm{Daniel}\binits{D.}},
\bauthor{\bsnm{Joseph},~\bfnm{Mathew}\binits{M.}} \AND
\bauthor{\bsnm{Khoshnevisan},~\bfnm{Davar}\binits{D.}}
(\byear{2013}).
\btitle{On the chaotic character of the stochastic heat equation,
before the onset of intermitttency}.
\bjournal{Ann. Probab.}
\bvolume{41}
\bpages{2225--2260}.
\bid{doi={10.1214/11-AOP717}, issn={0091-1798}, mr={3098071}}
\end{barticle}
%

\bptok{imsref}%
% NOT OUTPUTED:
% number = 3B
% doi = http://dx.doi.org/10.1214/11-AOP717
% fjournal = The Annals of Probability
\endbibitem

%b18 ###
%b18 #&#
\bibitem{CFG}
%
\begin{barticle}[mr]
\bauthor{\bsnm{Cox},~\bfnm{J.~Theodore}\binits{J.~T.}},
\bauthor{\bsnm{Fleischmann},~\bfnm{Klaus}\binits{K.}} \AND
\bauthor{\bsnm{Greven},~\bfnm{Andreas}\binits{A.}}
(\byear{1996}).
\btitle{Comparison of interacting diffusions and an application to
their ergodic theory}.
\bjournal{Probab. Theory Related Fields}
\bvolume{105}
\bpages{513--528}.
\bid{doi={10.1007/BF01191911}, issn={0178-8051}, mr={1402655}}
\end{barticle}
%

\bptok{imsref}%
% NOT OUTPUTED:
% number = 4
% doi = http://dx.doi.org/10.1007/BF01191911
% coden = PTRFEU
% fjournal = Probability Theory and Related Fields
\endbibitem

%b19 ###
%b19 #&#
\bibitem{CoxGreven}
%
\begin{barticle}[mr]
\bauthor{\bsnm{Cox},~\bfnm{J.~T.}\binits{J.~T.}} \AND
\bauthor{\bsnm{Greven},~\bfnm{Andreas}\binits{A.}}
(\byear{1994}).
\btitle{Ergodic theorems for infinite systems of locally interacting
diffusions}.
\bjournal{Ann. Probab.}
\bvolume{22}
\bpages{833--853}.
\bid{issn={0091-1798}, mr={1288134}}
\end{barticle}
%

\bptok{imsref}%
% NOT OUTPUTED:
% url =
%http://links.jstor.org/sici?sici=0091-1798(199404)22:2<833:ETFISO>2.0.CO;2-I&origin=MSN
% number = 2
% coden = APBYAE
% fjournal = The Annals of Probability
\endbibitem

%b20 ###
%b20 #&#
\bibitem{CGM}
%
\begin{barticle}[mr]
\bauthor{\bsnm{Cranston},~\bfnm{M.}\binits{M.}},
\bauthor{\bsnm{Gauthier},~\bfnm{D.}\binits{D.}} \AND
\bauthor{\bsnm{Mountford},~\bfnm{T.~S.}\binits{T.~S.}}
(\byear{2010}).
\btitle{On large deviations for the parabolic {A}nderson model}.
\bjournal{Probab. Theory Related Fields}
\bvolume{147}
\bpages{349--378}.
\bid{doi={10.1007/s00440-009-0249-z}, issn={0178-8051}, mr={2594357}}
\end{barticle}
%

\bptok{imsref}%
% NOT OUTPUTED:
% number = 1-2
% doi = http://dx.doi.org/10.1007/s00440-009-0249-z
% coden = PTRFEU
% fjournal = Probability Theory and Related Fields
\endbibitem

%b21 ###
%b21 #&#
\bibitem{CranstonMolchanov}
%
\begin{barticle}[mr]
\bauthor{\bsnm{Cranston},~\bfnm{M.}\binits{M.}} \AND
\bauthor{\bsnm{Molchanov},~\bfnm{S.}\binits{S.}}
(\byear{2007}).
\btitle{Quenched to annealed transition in the parabolic {A}nderson problem}.
\bjournal{Probab. Theory Related Fields}
\bvolume{138}
\bpages{177--193}.
\bid{doi={10.1007/s00440-006-0020-7}, issn={0178-8051}, mr={2288068}}
\end{barticle}
%

\bptok{imsref}%
% NOT OUTPUTED:
% number = 1-2
% doi = http://dx.doi.org/10.1007/s00440-006-0020-7
% coden = PTRFEU
% fjournal = Probability Theory and Related Fields
\endbibitem

%b22 ###
%b22 #&#
\bibitem{CranstonMountfordShiga}
%
\begin{barticle}[mr]
\bauthor{\bsnm{Cranston},~\bfnm{M.}\binits{M.}},
\bauthor{\bsnm{Mountford},~\bfnm{T.~S.}\binits{T.~S.}} \AND
\bauthor{\bsnm{Shiga},~\bfnm{T.}\binits{T.}}
(\byear{2002}).
\btitle{Lyapunov exponents for the parabolic {A}nderson model}.
\bjournal{Acta Math. Univ. Comenian. (N.S.)}
\bvolume{71}
\bpages{163--188}.
\bid{issn={0862-9544}, mr={1980378}}
\end{barticle}
%

\bptok{imsref}%
% NOT OUTPUTED:
% number = 2
% fjournal = Acta Mathematica Universitatis Comenianae. New Series
\endbibitem

%b23 ###
%b23 #&#
\bibitem{Davis}
%
\begin{barticle}[mr]
\bauthor{\bsnm{Davis},~\bfnm{Burgess}\binits{B.}}
(\byear{1976}).
\btitle{On the {$L^{p}$} norms of stochastic integrals and other martingales}.
\bjournal{Duke Math. J.}
\bvolume{43}
\bpages{697--704}.
\bid{issn={0012-7094}, mr={0418219}}
\end{barticle}
%

\bptok{imsref}%
% NOT OUTPUTED:
% number = 4
% fjournal = Duke Mathematical Journal
\endbibitem

%b24 ###
%b24 #&#
\bibitem{FellerVolII}
%
\begin{bbook}[mr]
\bauthor{\bsnm{Feller},~\bfnm{William}\binits{W.}}
(\byear{1966}).
\btitle{An Introduction to Probability Theory and Its Applications.
{V}ol. {II}}.
\bpublisher{Wiley},
\blocation{New York}.
\bid{mr={0210154}}
\end{bbook}
%

\bptok{imsref}%
% NOT OUTPUTED:
% fpage = xviii+636
\endbibitem

%b25 ###
%b25 #&#
\bibitem{FK}
%
\begin{barticle}[mr]
\bauthor{\bsnm{Foondun},~\bfnm{Mohammud}\binits{M.}} \AND
\bauthor{\bsnm{Khoshnevisan},~\bfnm{Davar}\binits{D.}}
(\byear{2009}).
\btitle{Intermittence and nonlinear parabolic stochastic partial
differential equations}.
\bjournal{Electron. J. Probab.}
\bvolume{14}
\bpages{548--568}.
\bid{doi={10.1214/EJP.v14-614}, issn={1083-6489}, mr={2480553}}
\end{barticle}
%

\bptok{imsref}%
% NOT OUTPUTED:
% doi = http://dx.doi.org/10.1214/EJP.v14-614
% fjournal = Electronic Journal of Probability
\endbibitem

%b26 ###
%b26 #&#
\bibitem{Funaki}
%
\begin{barticle}[mr]
\bauthor{\bsnm{Funaki},~\bfnm{Tadahisa}\binits{T.}}
(\byear{1983}).
\btitle{Random motion of strings and related stochastic evolution equations}.
\bjournal{Nagoya Math. J.}
\bvolume{89}
\bpages{129--193}.
\bid{issn={0027-7630}, mr={0692348}}
\end{barticle}
%

\bptok{imsref}%
% NOT OUTPUTED:
% url = http://projecteuclid.org/euclid.nmj/1118787110
% coden = NGMJA2
% fjournal = Nagoya Mathematical Journal
\endbibitem

%b27 ###
%b27 #&#
\bibitem{GM}
%
\begin{barticle}[mr]
\bauthor{\bsnm{Gei{\ss}},~\bfnm{Christel}\binits{C.}} \AND
\bauthor{\bsnm{Manthey},~\bfnm{Ralf}\binits{R.}}
(\byear{1994}).
\btitle{Comparison theorems for stochastic differential equations in
finite and infinite dimensions}.
\bjournal{Stochastic Process. Appl.}
\bvolume{53}
\bpages{23--35}.
\bid{doi={10.1016/0304-4149(94)90055-8}, issn={0304-4149}, mr={1290705}}
\end{barticle}
%

\bptok{imsref}%
% NOT OUTPUTED:
% number = 1
% doi = http://dx.doi.org/10.1016/0304-4149(94)90055-8
% coden = STOPB7
% fjournal = Stochastic Processes and their Applications
\endbibitem

%b28 ###
%b28 #&#
\bibitem{GK}
%
\begin{bbook}[mr]
\bauthor{\bsnm{Gnedenko},~\bfnm{B.~V.}\binits{B.~V.}} \AND
\bauthor{\bsnm{Kolmogorov},~\bfnm{A.~N.}\binits{A.~N.}}
(\byear{1968}).
\btitle{Limit Distributions for Sums of Independent Random Variables}.
%\bseries{Translated from the Russian, Annotated,
%and Revised by K. L. Chung.
%With Appendices by J. L. Doob and P. L. Hsu. Revised Edition}.
\bpublisher{Addison-Wesley},
\blocation{Reading, MA}.
\bid{mr={0233400}}
\end{bbook}
%

\bptok{imsref}%
% NOT OUTPUTED:
% fpage = ix+293
\endbibitem

%b29 ###
%b29 #&#
\bibitem{GrevenDenHollander}
%
\begin{barticle}[mr]
\bauthor{\bsnm{Greven},~\bfnm{A.}\binits{A.}} \AND
\bauthor{\bsnm{den Hollander},~\bfnm{F.}\binits{F.}}
(\byear{2007}).
\btitle{Phase transitions for the long-time behavior of interacting
diffusions}.
\bjournal{Ann. Probab.}
\bvolume{35}
\bpages{1250--1306}.
\bid{doi={10.1214/009117906000001060}, issn={0091-1798}, mr={2330971}}
\end{barticle}
%

\bptok{imsref}%
% NOT OUTPUTED:
% number = 4
% doi = http://dx.doi.org/10.1214/009117906000001060
% coden = APBYAE
% fjournal = The Annals of Probability
\endbibitem

%%b30 ###
%%b30 #&#
%\bibitem{IkedaWatanabe}
%%
%\begin{bbook}[mr]
%\bauthor{\bsnm{Ikeda},~\bfnm{Nobuyuki}\binits{N.}} \AND
%\bauthor{\bsnm{Watanabe},~\bfnm{Shinzo}\binits{S.}}
%(\byear{1989}).
%\btitle{Stochastic Differential Equations and Diffusion Processes},
%\bedition{2nd} ed.
%\bseries{North-Holland Mathematical Library}
%\bvolume{24}.
%\bpublisher{North-Holland},
%\blocation{Amsterdam}.
%\bid{mr={1011252}}
%\bptnote{check year}%
%\end{bbook}
%%
%
%\bptok{imsref}%
%% NOT OUTPUTED:
%% isbn = 0-444-87378-3
%% fpage = xvi+555
%\endbibitem

%b31 ###
%b31 #&#
\bibitem{Iwata}
%
\begin{barticle}[mr]
\bauthor{\bsnm{Iwata},~\bfnm{Koichiro}\binits{K.}}
(\byear{1987}).
\btitle{An infinite-dimensional stochastic differential equation with
state space {$C(\mathbf{R})$}}.
\bjournal{Probab. Theory Related Fields}
\bvolume{74}
\bpages{141--159}.
\bid{doi={10.1007/BF01845644}, issn={0178-8051}, mr={0863723}}
\end{barticle}
%

\bptok{imsref}%
% NOT OUTPUTED:
% number = 1
% doi = http://dx.doi.org/10.1007/BF01845644
% coden = PTRFEU
% fjournal = Probability Theory and Related Fields
\endbibitem

%b32 ###
%b32 #&#
\bibitem{K-CBMS}
%
\begin{bbook}[mr]
\bauthor{\bsnm{Khoshnevisan},~\bfnm{Davar}\binits{D.}}
(\byear{2014}).
\btitle{Analysis of Stochastic Partial Differential Equations}.
\bseries{CBMS Regional Conference Series in Mathematics}
\bvolume{119}.
\bpublisher{Amer. Math. Soc.},
\blocation{Providence, RI}.
\bid{mr={3222416}}
\end{bbook}
%

\bptok{imsref}%
% NOT OUTPUTED:
% isbn = 978-1-4704-1547-1
% fpage = viii+116
\endbibitem

%b33 ###
%b33 #&#
\bibitem{Mueller}
%
\begin{barticle}[mr]
\bauthor{\bsnm{Mueller},~\bfnm{Carl}\binits{C.}}
(\byear{1991}).
\btitle{On the support of solutions to the heat equation with noise}.
\bjournal{Stochastics Stochastics Rep.}
\bvolume{37}
\bpages{225--245}.
\bid{issn={1045-1129}, mr={1149348}}
\end{barticle}
%

\bptok{imsref}%
% NOT OUTPUTED:
% number = 4
% coden = STOCBS
% fjournal = Stochastics and Stochastics Reports
\endbibitem

%b34 ###
%b34 #&#
\bibitem{MN}
%
\begin{barticle}[mr]
\bauthor{\bsnm{Mueller},~\bfnm{Carl}\binits{C.}} \AND
\bauthor{\bsnm{Nualart},~\bfnm{David}\binits{D.}}
(\byear{2008}).
\btitle{Regularity of the density for the stochastic heat equation}.
\bjournal{Electron. J. Probab.}
\bvolume{13}
\bpages{2248--2258}.
\bid{doi={10.1214/EJP.v13-589}, issn={1083-6489}, mr={2469610}}
\end{barticle}
%

\bptok{imsref}%
% NOT OUTPUTED:
% doi = http://dx.doi.org/10.1214/EJP.v13-589
% fjournal = Electronic Journal of Probability
\endbibitem

%b35 ###
%b35 #&#
\bibitem{PrevotRoeckner}
%
\begin{bbook}[mr]
\bauthor{\bsnm{Pr{\'e}v{\^o}t},~\bfnm{Claudia}\binits{C.}} \AND
\bauthor{\bsnm{R{\"o}ckner},~\bfnm{Michael}\binits{M.}}
(\byear{2007}).
\btitle{A Concise Course on Stochastic Partial Differential Equations}.
\bseries{Lecture Notes in Math.}
\bvolume{1905}.
\bpublisher{Springer},
\blocation{Berlin}.
\bid{mr={2329435}}
\end{bbook}
%

\bptok{imsref}%
% NOT OUTPUTED:
% isbn = 978-3-540-70780-6; 3-540-70780-8
% fpage = vi+144
\endbibitem

%b36 ###
%b36 #&#
\bibitem{RevuzYor}
%
\begin{bbook}[mr]
\bauthor{\bsnm{Revuz},~\bfnm{Daniel}\binits{D.}} \AND
\bauthor{\bsnm{Yor},~\bfnm{Marc}\binits{M.}}
(\byear{1991}).
\btitle{Continuous Martingales and {B}rownian Motion}.
\bseries{Grundlehren der mathematischen Wissenschaften}
\bvolume{293}.
\bpublisher{Springer},
\blocation{Berlin}.
\bid{doi={10.1007/978-3-662-21726-9}, mr={1083357}}
\end{bbook}
%

\bptok{imsref}%
% NOT OUTPUTED:
% doi = http://dx.doi.org/10.1007/978-3-662-21726-9
% isbn = 3-540-52167-4
% fpage = x+533
\endbibitem

%b37 ###
%b37 #&#
\bibitem{Shiga}
%
\begin{barticle}[mr]
\bauthor{\bsnm{Shiga},~\bfnm{Tokuzo}\binits{T.}}
(\byear{1992}).
\btitle{Ergodic theorems and exponential decay of sample paths for
certain interacting diffusion systems}.
\bjournal{Osaka J. Math.}
\bvolume{29}
\bpages{789--807}.
\bid{issn={0030-6126}, mr={1192741}}
\end{barticle}
%

\bptok{imsref}%
% NOT OUTPUTED:
% url = http://projecteuclid.org/euclid.ojm/1200784090
% number = 4
% coden = OJMAA7
% fjournal = Osaka Journal of Mathematics
\endbibitem

%b38 ###
%b38 #&#
\bibitem{ShigaShimizu}
%
\begin{barticle}[mr]
\bauthor{\bsnm{Shiga},~\bfnm{Tokuzo}\binits{T.}} \AND
\bauthor{\bsnm{Shimizu},~\bfnm{Akinobu}\binits{A.}}
(\byear{1980}).
\btitle{Infinite-dimensional stochastic differential equations and
their applications}.
\bjournal{J. Math. Kyoto Univ.}
\bvolume{20}
\bpages{395--416}.
\bid{issn={0023-608X}, mr={0591802}}
\end{barticle}
%

\bptok{imsref}%
% NOT OUTPUTED:
% number = 3
% coden = JMKYAZ
% fjournal = Journal of Mathematics of Kyoto University
\endbibitem

%%b39 ###
%%b39 #&#
%\bibitem{Walsh}
%%
%\begin{bincollection}[mr]
%\bauthor{\bsnm{Walsh},~\bfnm{John~B.}\binits{J.~B.}}
%(\byear{1986}).
%\btitle{An introduction to stochastic partial differential equations}.
%In \bbooktitle{\'Ecole D'\'et\'e de Probabilit\'es de {S}aint-{F}lour,
%XIV---1984}.
%\bseries{Lecture Notes in Math.}
%\bvolume{1180}
%\bpages{265--439}.
%\bpublisher{Springer},
%\blocation{Berlin}.
%\bid{doi={10.1007/BFb0074920}, mr={0876085}}
%\end{bincollection}
%%
%
%\bptok{imsref}%
%% NOT OUTPUTED:
%% doi = http://dx.doi.org/10.1007/BFb0074920
%\endbibitem
\end{thebibliography}
\end{document}